\documentclass[11pt]{article}
\usepackage{graphicx}
\usepackage{amsmath}
\usepackage{amsthm}
\usepackage{latexsym,bm}
\usepackage{amssymb}

\usepackage{amsfonts}

\usepackage{indentfirst}
\usepackage{CJK}
\newtheorem{thm}{Theorem}[section]
\newtheorem{cor}[thm]{Corollary}
\newtheorem{lem}[thm]{Lemma}
\newtheorem{prop}[thm]{Proposition}
\newtheorem{defn}[thm]{Definition}
\newtheorem{rem}[thm]{Remark}

\begin{document}
\title{{The geometry of $\Phi_{(3)}$-harmonic maps}
\footnotetext{{Mathematics Classification
Primary(2010)}: {58E20, 53C21, 53C25.
}\\
\hspace*{6mm}{Keywords: $\Phi_{(3)}$-harmonic maps, Liouville type results, variation formula, $\Phi_{(3)}$-SSU manifold.\\}
}}
\author{{Shuxiang Feng, Yingbo Han, Kaige Jiang and Shihshu Walter Wei}}

\date{}
\maketitle
\begin{abstract}

In this paper, we motivate and extend the study of harmonic maps or $\Phi_{(1)}$-harmonic maps (cf \cite {es}, Remark 1.3 (iii)), $\Phi$-harmonic maps or $\Phi_{(2)}$-harmonic maps (cf. \cite {hw}, Remark 1.3 (v)), and explore geometric properties of $\Phi_{(3)}$-harmonic maps by unified geometric analytic methods. We define the notion of $\Phi_{(3)}$-harmonic maps and obtain the first variation formula and the second variation formula of the $\Phi_{(3)}$-energy functional $E_{\Phi_{(3)}}$. By using a stress-energy tensor and the asymptotic assumption of maps at infinity, we prove Liouville type results for  $\Phi_{(3)}$-harmonic maps. We introduce the notion of $\Phi_{(3)}$-Superstrongly Unstable ($\Phi_{(3)}$-SSU) manifold and provide many interesting examples. By using an extrinsic average variational method in the calculus of variations (cf. \cite {ww0, W0}), we find $\Phi_{(3)}$-SSU manifold and
prove that any stable $\Phi_{(3)}$-harmonic maps from a compact $\Phi_{(3)}$-SSU manifold (into any compact Riemannian manifold) or (from any compact Riemannian manifold) into a compact $\Phi_{(3)}$-SSU manifold must be constant. We also prove that the homotopy class of any map from a compact $\Phi_{(3)}$-SSU manifold (into any compact Riemannian manifold) or (from any compact Riemannian manifold) into a compact $\Phi_{(3)}$-SSU manifold contains elements of arbitrarily small $\Phi_{(3)}$-energy.
We call a compact Riemannian manifold $M$ to be $\Phi_{(3)}$-strongly unstable ($\Phi_{(3)}$-SU) if it is not the target or domain of a nonconstant stable $\Phi_{(3)}$-harmonic map (from or into any compact Riemannian manifold) and also the homotopy class of any map to or from $M$ (from or into any compact Riemannian manifold) contains elements of arbitrarily small $\Phi_{(3)}$-energy. We prove that every compact $\Phi_{(3)}$-$\operatorname{SSU}$ manifold is $\Phi_{(3)}$-$\operatorname{SU}$.
As consequences, we obtain topological vanishing theorems and sphere theorems by employing a $\Phi_{(3)}$-harmoic map as a catalyst. This is in contrast to the approaches of utilizing a geodesic (\cite {sy}), minimal surface, stable rectifiable current (\cite {ls, HoW1, W01}), $p$-harmonic map (cf. \cite {wwre}), etc., as catalysts. These mysterious phenomena are analogs of harmonic maps or $\Phi_{(1)}$-harmonic maps, $p$-harmonic maps,  $\Phi_{S}$-harmonic maps, $\Phi_{S,p}$-harmonic maps, $\Phi_{(2)}$-harmonic maps, etc., (cf.
Howard and Wei, 1986; Wei, 1992; Wei and Yau, 1994; Wei, 1998; Feng-Han-Li-Wei, 2021; Feng-Han-Wei, 2021).
\end{abstract}

\section{Introduction}
Harmonic maps or $\Phi_{(1)}$-harmonic maps (cf \cite {es}, Remark 1.3 (iii)) which appear in a broad spectrum of contexts in mathematics and physics, have had wide-ranging consequences and influenced developments in other fields (see, e.g., \cite{hy11,hy02,es,hy03,han42}).
From an algebraic invariant point of view, a harmonic map (or $\Phi_{(1)}$-harmonic map) $u : (M, g) \to (N, h)$ between Riemannian manifolds $(M, g)$ and $(N, h)$ can be viewed as a critical point of the energy or $\Phi_{(1)}$-energy
functional, given by the integral of a half $(=\frac {1}{{\bold 1} \cdot 2})$ of
{\it first} elementary symmetric function $\sigma_1$ of eigenvalues of the pullback metric tensor $u^{\ast}h$ relative to the metric $g$.
Similarly, a $\Phi$-harmonic map or $\Phi_{(2)}$-harmonic map (cf. \cite {hw}, Remark 1.3 (v)) is a critical point of the $\Phi$-energy or $\Phi_{(2)}$-energy
functional, given by the integral of a quarter $(=\frac {1}{{\bold 2} \cdot 2})$ of the {\it second} elementary symmetric function $\sigma_2$ of eigenvalues of the pullback metric tensor $u^{\ast}h$ relative to the metric $g$.
Liouville type results for $p$-harmonic maps, $F$-harmonic maps, $F$-stationary  maps, CC-harmonic maps, $\Phi$-harmonic maps, $\Phi_S$-harmonic maps and $\Phi_{S,p}$-harmonic maps are proved by several authors (see \cite{hy09,fhlw,fhw,han61,hw,hy10,T} for more details). In contrast to the usual method of proving Liouville type results by assuming the finiteness of the energy of the map or the smallness of whole image of the domain manifold under the map, Jin \cite{hy04} and Feng-Han \cite{fh} obtain Liouville type results under natural conditions about the asymptotic behavior of maps at infinity.

We propose an extrinsic, average variational method as an approach to confront and resolve problems in global, nonlinear analysis and geometry $($cf. \cite {ww0, W0}$)$.
In contrast to an average method in PDE that we applied in \cite {cw} to obtain sharp growth estimates for warping functions in multiply warped product manifolds, we employ \emph {an extrinsic average variational method} in
the calculus of variations $($\cite {ww0}$)$,
find a large class of manifolds of positive Ricci curvature that enjoy rich properties, and introduce the notions of \emph {superstrongly unstable $(\operatorname{SSU})$ manifolds} $($\cite {ww}$)$ and \emph {$p$-superstrongly unstable $(p$-$\operatorname{SSU})$ manifolds} $($\cite {wy}$)$.
Wei-Yau $($\cite {wy}$)$ discuss Liouville type theorems and the regularity of $p$-minimizers into $p$-$\operatorname{SSU}$ manifolds, extending ground-breaking work of Hardt-Lin (\cite{hl}) and Luckhaus (\cite{lu}).
The method (\cite {ww0}) can be carried over to more general settings: Han and Wei find $\Phi$-SSU manifolds (or $\Phi_{(2)}$-SSU manifolds) \cite{hw} and prove that every compact $
\Phi$-SSU manifold must be $\Phi$-strongly unstable ($\Phi$-SU), or $\Phi_{(2)}$-strongly unstable ($\Phi_{(2)}$-SU), i.e., every compact $\Phi$-SSU manifold can neither be the domain nor the target of any nonconstant smooth stable $\Phi$-harmonic map between two compact Riemannian manifolds and the homotopic class of maps (between two compact Riemannian manifolds) from or into $\Phi$-SSU manifold contains a map of arbitrarily small $\Phi$-energy. These generalize the cases $S^{m},$ for $m > 5$, and compact minimal submanifolds of $S^{m+p}$ with $Ric_{g}>\frac{3}{4}mg$, due to Kawai and Nakauchi \cite{han16, kn1}, in which the nonexistence results of nonconstant stable $\Phi$-harmonic maps are extended to $\Phi$-SU or $\Phi_{(2)}$-SU features
in \cite {hw}.
Cases of SSU manifolds or $\Phi_{(1)}$-SSU manifolds such as spheres, product of spheres with appropriate dimensions, etc. can be found in \cite {W0}.
For some results of the $\Phi$-energy functional, we refer to \cite{han13,han51}. In \cite{fhlw}, Feng-Han-Li-Wei prove that every compact $\Phi_{S}$-SSU manifold is $\Phi_{S}$-SU.
In particular any stable $\Phi_S$-harmonic map between compact Riemannian manifolds is a constant map if the domain manifold or target manifold is a $\Phi_{S}$-SSU manifold. In addition, they obtain properties of weakly conformal $\Phi_S$-harmonic map and horizontally weakly conformal $\Phi_S$-harmonic map. Feng-Han-Wei \cite{fhw} introduce the notions  of $\Phi_{S,p}$-harmonic maps, stable $\Phi_{S,p}$-harmonic maps, and prove that every compact $\Phi_{S,p}$-SSU manifold must be $\Phi_{S,p}$-SU. They obtain several Liouville type theorems by extending the method used in Jin~ \cite{hy04}.

For a symmetric $2$-covariant tensor field $\alpha$ on a hypersurface $M$ in $\mathbb R^{m+1}$,
at any fixed point $x_0 \in M\, ,$ $\alpha$ has the eigenvalues $\lambda$ relative to the metric $g$ of $M$; i.e., the $m$ real roots of the equation $\det (\delta_{ij} \lambda - \alpha_{ij}) = 0\, $ where $\alpha_{ij} = \alpha(e_i,e_j)\, ,$ and $\{e_1, \cdots e_m\}$ is a basis for $T_{x_0}(M)\, .$  The {\it algebraic} invariants --- the $k$-th elementary symmetric function of the eigenvalues of $\alpha$ at $x_0$, denoted by $\sigma_k (\alpha_{x_0}), 1 \le k \le m $ frequently have {\it geometric} meaning of the manifold $M$
with analytic, topological and physical impacts. For example, if we take $\alpha$ to be the second fundamental form of $M\, ,$ in $\mathbb R^{m+1}$,
 then $\frac{1}{m}\sigma _1 (\alpha), \frac {2}{m (m-1)}  \sigma _2 (\alpha)$,  and $\sigma _m (\alpha)$  are the mean curvature, scalar curvature,  and the Gauss-Kronecker curvature
of $M$ respectively and are central themes of Yamabi problem (\cite { A, J, S, Tr}), special Lagrangian graphs (\cite {hla}), geometric aspects of the theory of fully nonlinear elliptic equations (e.g., \cite {Sp}), and conformal geometry (e.g. \cite{CY}, \cite {DLW}), etc. If we take $\alpha$ to be Schoulten tensor, then a study of $\sigma_2 (\alpha)$ leads to a generalized Yamabe problem (\cite {CGY}).
In the study of prescribed curvature problems in PDE, the existence of closed starshaped hypersurfaces of prescribed mean curvature in Euclidean space was proved by A.E. Treibergs and S.W. Wei \cite {TW}, solving a problem of F. Almgren and S.T. Yau \cite {Y}. While the case of prescribed Guass-Kronecker curvature was studied by V.I. Oliker \cite {O} and P. Delano\"e \cite{D}, the case of prescribed $k$-th mean curvature, in particular the intermediate cases, $2 \le k \le m-1$ were treated by L. Caffarelli, L. Nirenberg and J. Spruck \cite{CNS}.

These motivate us from the viewpoint of geometric mapping theory $u: (M^m, g) \to (N^n,h)$, taking $\alpha = u^{\ast} h$, in this paper, to extend the study of harmonic maps or $\Phi_{(1)}$-harmonic maps (cf \cite {es}), $\Phi$-harmonic maps or $\Phi_{(2)}$-harmonic maps (cf. \cite {hw}), to explore geometric properties of $\Phi_{(3)}$-harmonic maps by unified geometric analytic methods. We define the notion of $\Phi_{(3)}$-harmonic maps and obtain the first variation formula and the second variation formula of the $\Phi_{(3)}$-energy functional $E_{\Phi_{(3)}}$.
{\it In fact, $\Phi_{(3)}$-harmonic map $($cf. \cite {hw}) is a critical point of the $\Phi_{(3)}$-energy
functional, given by the integral of a sixth $(=\frac {1}{{\bold 3} \cdot 2})$  of the {\it third} elementary symmetric function $\sigma_3$ of eigenvalues of the pullback metric tensor $u^{\ast}h$ relative to the metric $g$.}
We introduce the notion of $\Phi_{(3)}$-SSU manifold and provide many interesting examples. By using an extrinsic average variational method in the calculus of variations (cf. \cite {ww0, ww}), we find $\Phi_{(3)}$-SSU manifolds, and prove that every compact $\Phi_{(3)}$-SSU manifold is $\Phi_{(3)}$-SSU. As consequences, we prove topological vanishing theorems and sphere theorems by employing a $\Phi_{(3)}$-harmonic map as a catalyst.
This is in contrast to the approaches of utilizing a geodesic (\cite {sy}), minimal surface, stable rectifiable current (\cite {ls, HoW1, W01}), $p$-harmonic map (cf. \cite {wwre}), etc., as catalysts. These mysterious phenomena are analogs of harmonic maps or $\Phi_{(1)}$-harmonic maps, $p$-harmonic maps,  $\Phi_{S}$-harmonic maps, $\Phi_{S,p}$-harmonic maps, $\Phi_{(2)}$-harmonic maps, etc., (cf. \cite{HoW,ww, wy, wwre, fhlw, fhw}).

While we can view the differential of $u$, denoted by $du$, as a $1$-form with values in the pullback bundle $u^{-1}TN$  as $d_{(1)} u$,
in this paper, we introduce the following unified notions.

\begin{defn}
	Let $d_{(1)} u , d_{(2)} u$ and $d_{(3)} u$ be  $1$-forms with values in the pullback bundle $u^{-1}TN$ given by
	\begin{equation}
	\begin{aligned}
	d_{(1)} u (X) &= du (X)\, ,\\
 d_{(2)} u (X)&=\sum_{j=1}^mh\big (du(X),du(e_j)\big )du(e_j)\, ,\quad \operatorname{and} \\
d_{(3)} u (X)&=\sum_{j,k=1}^m h\big (du(X),du(e_j)\big )h\big (du(e_j),du(e_k)\big )du(e_k)\, , \end{aligned} \label{1.1}\end{equation}
	respectively, for any smooth vector field $X$ on $(M,g)$, where $\{e_i\}$ is a local orthonormal frame field on $(M,g)$, with the following {\it corresponding norms}
	\begin{equation*}
	\begin{aligned}
	||d_{(1)} u||^2 &= \sum_{i=1}^mh\big (d_{(1)} u(e_i),du(e_i)\big )= \sum_{i=1}^mh\big (du(e_i),du(e_i)\big )\, ,\\
 ||d_{(2)} u ||^2&= \sum_{i=1}^mh\big (d_{(2)} u(e_i),du(e_i)\big )=\sum_{i,j=1}^mh\big (du(e_i),du(e_j)\big )h\big( du(e_j), du(e_i)\big )\, ,\quad \operatorname{and} \\
  ||d_{(3)} u ||^2&=\sum_{i=1}^mh\big (d_{(3)} u(e_i),du(e_i)\big )
 =\sum_{i,j,k=1}^mh\big (du(e_i),du(e_j)\big )h\big (du(e_j),du(e_k)\big )h\big (du(e_k),du(e_i)\big ) .
 \end{aligned} \end{equation*}

	
	The {\it $\Phi_{(1)}$-energy density $e_{\Phi_{(1)}}(u)$,  $\Phi_{(2)}$-energy density $e_{\Phi_{(2)}}(u)$, and $\Phi_{(3)}$-energy density $e_{\Phi_{(3)}}(u)$ of $u$} are given by
\begin{equation}
\begin{aligned}
	e_{\Phi_{(1)}}(u)&=\frac{||d_{(1)} u||^2}{2}, \\
	e_{\Phi_{(2)}}(u)&=\frac{||d_{(2)} u||^2}{4},\quad \operatorname{and}  \\
	e_{\Phi_{(3)}}(u)&=\frac{||d_{(3)} u||^2}{6}, \quad \operatorname{respectively}.
	\end{aligned} \end{equation}
	
	The $\Phi_{(1)}$-energy $E_{\Phi_{(1)}}(u)$, $\Phi_{(2)}$-energy $E_{\Phi_{(2)}}(u)$, and $\Phi_{(3)}$-energy $E_{\Phi_{(3)}}(u)$ of $u$ are given by
\begin{equation}
\begin{aligned}
E_{\Phi_{(1)}}(u)&=\int_Me_{\Phi_{(1)}}(u)dv_g, \\
	E_{\Phi_{(2)}}(u)&=\int_Me_{\Phi_{(2)}}(u)dv_g,\quad \operatorname{and}  \\
	E_{\Phi_{(3)}}(u)&=\int_Me_{\Phi_{(3)}}(u)dv_g,\quad \operatorname{respectively} .
	\end{aligned} \end{equation}
\end{defn}
\begin{defn}
	For $i=1,2,3$, a smooth map $u$ is said to be a $\Phi_{(i)}$-harmonic map if it is a critical point of the $\Phi_{(i)}$-energy functional $E_{\Phi_{(i)}}$ with respect to any smooth compactly supported variation of $u$, stable $\Phi_{(i)}$-harmonic or simply $\Phi_{(i)}$-stable, if $u$ is a local minimum of $E_{\Phi_{(i)}}(u)$, and $\Phi_{(i)}$-unstable if $u$ is not $\Phi_{(i)}$-stable.\label{D: 1.2}
\end{defn}

\begin{rem} $(\operatorname{i})$ The norm $||d_{(1)} u||$ is the Hibert-Schmid norm of the differential $du$, i.e., $||d_{(1)} u|| = |du|\, .$ $(\operatorname {ii})$ The $\Phi_{(1)}$-energy density $e_{\Phi_{(1)}}(u) = e(u)$ is the energy density of $u$. $(\operatorname {iii})$ $\Phi_{(1)}$-harmonic map is ordinary harmonic map $($cf. \cite {es} $)$.
 $(\operatorname {iv})$ The $\Phi_{(2)}$-energy density $e_{\Phi_{(2)}}(u) = e_{\Phi}(u)$ is the $\Phi$-energy density of $u$. $(\operatorname {v})$  $\Phi_{(2)}$-harmonic map is $\Phi$-harmonic map $($cf. \cite {hw} $)$.  $(\operatorname {vi})$ Definition \ref{D: 1.2} can be extended to $4 \le i \le m=\dim M$. Hence, for any integer $1 \le i \le m$, a smooth map $u$ is said to be a $\Phi_{(i)}$-harmonic map if it is a critical point of the $\Phi_{(i)}$-energy functional $E_{\Phi_{(i)}}$ with respect to any smooth compactly supported variation of $u$, stable $\Phi_{(i)}$-harmonic or simply $\Phi_{(i)}$-stable, if $u$ is a local minimum of $E_{\Phi_{(i)}}(u)$, and $\Phi_{(i)}$-unstable if $u$ is not $\Phi_{(i)}$-stable.
\end{rem}

We recall

\begin{defn}[\cite{wwre}]
	A Riemannian manifold $M^{m}$ is said to be superstrongly
	unstable $(\operatorname{SSU})$ if there exists an isometric immersion
	of $M^{m}$ in $\mathbb{R}^{q}$ with its second fundamental form  $B$ such that for all unit
	tangent vectors $v$ to $M^{m}$ at every point $x\in M^{m},$ the
	following functional is negative valued.
	\begin{equation*}
	\langle Q_x^M(v),v \rangle_M=\sum_{i=1}^{m}\big (2\langle B(v,e_{i}),B(v,e_{i})\rangle-\langle B(v,v),B(e_{i},e_{i})\rangle \big ),
	\end{equation*}
	where $\{e_{1},\cdots,e_{m}\}$ is a local orthonormal
frame field on $M^{m}$ near $x$. A Riemannian manifold $M$ is said to be $p$-superstrongly unstable $(p$-$\operatorname{SSU})$ for $p \geq 2$ if the following functional is negative valued.
\begin{equation*}
F_{p,x} (v)=(p-2)\langle B(v,v),B(v,v)\rangle +\langle Q_x^M(v),v \rangle_M.
\end{equation*}
\end{defn}
The notion of $\Phi_{(3)}$-SSU is defined as follows.

\begin{defn}
	A Riemannian manifold $M^{m}$ is said to be $\Phi_{(3)}$-superstrongly
	unstable $($$\Phi_{(3)}$-$\operatorname{SSU}$$)$ if there exists an isometric immersion
	of $M^{m}$ in $\mathbb{R}^{q}$ with its second fundamental form  $B$ such that for all unit
	tangent vectors $v$ to $M^{m}$ at every point $x\in M^{m},$ the
	following functional is negative valued.
	\begin{equation}
	{F_{\Phi_{(3)},x}} (v)=\sum_{i=1}^{m}\big (6\langle B(v,e_{i}),B(v,e_{i})\rangle _{\mathbb R^q}-\langle B(v,v),B(e_{i},e_{i})\rangle_{\mathbb R^q}\big ),\label{eq:SSU}
	\end{equation}
	where $\{e_{1},\cdots,e_{m}\}$ is a local orthonormal
	frame field on $M^{m}$ near $x$.
\end{defn}

\begin{defn}\label{D: 1.7}
A compact Riemannian manifold $M^{m}$ is $\Phi_{(3)}$-{\it strongly unstable} $(\Phi_{(3)}$-$\operatorname{SU})$ if
it is neither the domain nor the target of any nonconstant smooth $\Phi_{(3)}$-stable harmonic map $($into or from any compact Riemannian manifold$)$, and the homotopic class of maps from or into $M$ $($into or from any compact Riemannian manifold$)$ contains a map of arbitrarily small $\Phi_{(3)}$-energy $E_{\Phi_{(3)}}$.
\end{defn}
This leads to the study of the identity map on a Riemannian manifold. In particular, if $M$ is $\Phi_{(3)}$-$\operatorname{SU}$, then the identity map of $M$ is $\Phi_{(3)}$-unstable. For convenience, we make the following
\begin{defn}
A Riemannian manifold $M^{m}$ is $\Phi_{(3)}$-unstable $(\Phi_{(3)}$-$\operatorname{U})$ if the identity map $\text{Id}_M$ on $M^{m}$ is $\Phi_{(3)}$-unstable.
\end{defn}

This gives a natural analog in the setting of $\Phi_{(1)}$-harmonic maps and $\Phi_{(2)}$-harmonic maps:
\begin{prop}\label{P:1.8}

 Let $M$ be a compact manifold. Then

\[ M\,  \operatorname{is}\, \Phi_{(2)}\operatorname{-SSU}\quad \Rightarrow\quad M \operatorname{is}\, \Phi_{(2)}\operatorname{-SU}\quad  \Rightarrow \quad M \,  \operatorname{is}\, \Phi_{(2)}\operatorname{-U}.
\]
\end{prop}
We organize this paper in the following way. In Section \ref{sec2}, we obtain some fundamental results which will be used in the subsequent sections,
such as the first variation formula $($I$)$ and $($II$)$ in two different settings and the stress energy tensor $S_{\Phi_{(3)}}$ with respect to the functional $E_{\Phi_{(3)}}$. In Section \ref{sec4}, we extend Jin's method for solving a uniqueness problem  about $\Phi_{(3)}$-harmonic maps. To begin with, we obtain the lower energy growth rates of $\Phi_{(3)}$-harmonic maps by using the monotonicity formulas given in \cite{han} (cf. Proposition \ref{prop:L2}). Furthermore, we obtain the upper energy growth rates of these maps by using the asymptotic assumption of the maps at infinity (cf. Proposition \ref{prop:L3}). The two bounds are contradictory if $\Phi_{(3)}$-harmonic map is not a constant map (cf. Theorem \ref{thm:L4}). In Section \ref{sec5}, we calculate the second variation formula of the functional $E_{\Phi_{(3)}}$ (cf. Theorem \ref{thm:v1}) and give the concept of stable $\Phi_{(3)}$-harmonic maps with respect to the $\Phi_{(3)}$-energy $E_{\Phi_{(3)}}(u)$. In Section \ref{sec6}, we describe the motivating examples of the  $\Phi_{(3)}$-SSU manifolds and obtain the relation between  $\Phi_{(3)}$-SSU manifold and $p$-SSU manifold. Furthermore, we prove topological vanishing theorems and sphere theorems by the approach of $\Phi_{(3)}$-harmonic map
as a tool. In Section \ref{sec7}, by using an extrinsic average variational method in the calculus of variations \cite{ww0}, we prove that every stable $\Phi_{(3)}$-harmonic map from a compact $\Phi_{(3)}$-SSU manifold into any compact Riemannian manifold is constant (cf. Theorem \ref{thm:fss1}). In Section \ref{sec8}, by using the similar method to Theorem \ref{thm:fss1}, we obtain that every stable $\Phi_{(3)}$-harmonic maps from any compact Riemannian manifold into a compact $\Phi_{(3)}$-SSU manifolds is constant
(cf. Theorem \ref{thm:iss1}). In Section \ref{sec9}, we prove that the homotopic class of any map from every compact Riemannian manifold into a compact $\Phi_{(3)}$-SSU manifold contains elements of arbitrarily small $\Phi_{(3)}$-energy (cf. Theorem \ref{T: 9.2}). Finally, in Section \ref{sec10}, we prove that the homotopic class of any map from a compact $\Phi_{(3)}$-SSU manifold into every compact Riemannian manifold contains elements of arbitrarily small $\Phi_{(3)}$-energy (cf. Theorem \ref{T:10.2}). For $i=1,2,3$, we prove that every compact $\Phi_{(i)}$-$\operatorname{SSU}$ manifold is $\Phi_{(i)}$-$\operatorname{SU}$, and hence is $\Phi_{(i)}$-$\operatorname{U}$ (cf. Theorem \ref{T: 10.3}).
This generalizes Proposition \ref{P:1.8}, in which $i=2$.

\section{Preliminaries} \label{sec2}

In this section, we give some results which will be used in this paper, and extend the following concept from harmonic maps and $\Phi$-harmonic maps.

\begin{defn}	
The divergence of a $1$-form  with values in the pull-back bundle $u^{-1}(TN)\, ,$ $\tau_{\Phi_{(3)}} (u)$ $($resp. $\tau_{\Phi_{(2)}} (u)$, $\tau_{\Phi_{(1)}} (u))$
	is said to be the $\Phi_{(3)}$-tension field $($resp $\Phi_{(2)}$-tension field, $\Phi_{(1)}$-tension field$)$ of $u$ , if
	\begin{equation}
	\begin{aligned}
	\tau_{\Phi_{(3)}}(u) & = \operatorname {div} \big ( d_{(3)} u \big )\\
	& =\sum_{i,j,k=1}^m \bigg (\widetilde{\nabla}_{e_i} \big (h\big (du(e_i),du(e_j)\big )h\big (du(e_j),du(e_k)\big )du(e_k)\big )\bigg ) \\
	\bigg ( resp. \qquad \tau_{\Phi_{(2)}}(u) & = \operatorname {div} \big ( d_{(2)} u \big )
	 =\sum_{i,j =1}^m \bigg (\widetilde{\nabla}_{e_i} \big (h\big (du(e_i),du(e_j)\big ) du(e_j)\big )\bigg ),\\
	\quad \tau_{\Phi_{(1)}}(u) & = \operatorname {div} \big ( d_{(1)} u \big )
 =\sum_{i=1}^m \big ( \widetilde{\nabla}_{e_i} du(e_i)\big ) \qquad
	\bigg )\, .
	\end{aligned}
	\end{equation}
\end{defn}
We note that $\tau_{\Phi_{(1)}}(u) $ is the same as the tension filed $\tau (u) $ of $u$ (cf. \cite {es}) and $\tau_{\Phi_{(2)}}(u)$ is the same as  the $\Phi$-tension filed of $u$, denoted by $\tau_{\Phi}(u)$ (cf. \cite {hw}).
\begin{thm}[The first variation formula $(I)$]\label{Ft1}
	Let $u:(M^{m},g)\rightarrow(N,h)$ be a smooth map and let $u_t:(M^{m},g)\rightarrow(N,h)$, $(-\delta<t<\delta)$ be a family of compact supported variations such that $u_0=u$ and $v=\frac{\partial u_t}{\partial t}_{\big {|}_{t=0}}$.
	Then we have
	\begin{eqnarray}\label{FV}
	\frac{d E_{\Phi_{(3)}}(u_t)}{dt}_{\big {|}_{t=0}}=-\int_Mh\big (v,\tau_{\Phi_{(3)}}(u)\big )\, dv_g.
	\end{eqnarray}
\end{thm}
\begin{proof}
	We extend the vector field $\frac{\partial}{\partial t}$ on $(-\delta,\delta),$
	$X$ on $M$ naturally on $(-\delta,\delta)\times M,$ and denote
	those also by $\frac{\partial}{\partial t}, X$. Let
     $\nabla$ and $\widetilde{\nabla}$ be the Levi-Civita connection
	on $(-\delta,\delta)\times M$ and the induced connection on $u_{t}^{-1}TN$
	respectively.
	
	Note that
		 \begin{eqnarray}
		  && \quad \frac{\partial}{\partial t} \frac{|| d_{(3)} u_{t}||^2}{6}\nonumber\\
		 && =\frac{1}{6}\frac{\partial}{\partial t}\sum_{i=1}^mh\big (d_{(3) }u_{t}(e_{i}),du_{t}(e_{i})\big )\nonumber\\
		 && =\sum_{i,j,k=1}^mh\bigg (\widetilde{\nabla}_{\frac{\partial}{\partial t}}(du_{t}(e_{i})),du_{t}(e_{k})\bigg )h\big (du_{t}(e_{i}),du_{t}(e_{j})\big )h\big (du_{t}(e_{j}),du_{t}(e_{k})\big )\nonumber\\
		 && =\sum_{i,j,k=1}^m h\bigg (\widetilde{\nabla}_{e_{i}}(du_{t}(\frac{\partial}{\partial t})),du_{t}(e_{k})\bigg )h\big (du_{t}(e_{i}),du_{t}(e_{j})\big )h\big (du_{t}(e_{j}),du_{t}(e_{k})\big )\nonumber\\
		 && =\sum_{i,j,k=1}^m \bigg (e_{i}\, h\big (du_{t}(\frac{\partial}{\partial t}),du_{t}(e_{k})\big ) - h\big (du_{t}(\frac{\partial}{\partial t}),\nabla _{e_i} du_{t}(e_{k})\big ) \bigg )h\big (du_{t}(e_{i}),du_{t}(e_{j})\big )h\big (du_{t}(e_{j}),du_{t}(e_{k}) \big )\nonumber
\end{eqnarray}
		\begin{eqnarray}
		&& =\sum_{i,j,k=1}^m e_{i}\, h\big (du_{t}(\frac{\partial}{\partial t}),du_{t}(e_{k})\big )h\big (du_{t}(e_{i}),du_{t}(e_{j})\big )h\big (du_{t}(e_{j}),du_{t}(e_{k}) \big )\nonumber\\
		&& \qquad -h\bigg (du_{t}(\frac{\partial}{\partial t}), \widetilde{\nabla} _{e_i} du_{t}(e_{k}) \big )\, h\big (du_{t}(e_{i}),du_{t}(e_{j})\big )h\big (du_{t}(e_{j}),du_{t}(e_{k}) \bigg )\nonumber\\
		&& =\sum_{i,j,k=1}^m e_{i}\, h\big (du_{t}(\frac{\partial}{\partial t}),du_{t}(e_{k})\big )h\big (du_{t}(e_{i}),du_{t}(e_{j})\big )h\big (du_{t}(e_{j}),du_{t}(e_{k}) \big ),\nonumber\\
		&& \qquad -h\bigg (du_{t}(\frac{\partial}{\partial t}), h\big (du_{t}(e_{i}),du_{t}(e_{j})\big )h\big (du_{t}(e_{j}),du_{t}(e_{k}) \big )\, \widetilde{\nabla} _{e_i} du_{t}(e_{k}) \bigg )\nonumber\\
		&& =\sum_{i,j,k=1}^m e_{i}\, h\big (du_{t}(\frac{\partial}{\partial t}),du_{t}(e_{k})\big )h\big (du_{t}(e_{i}),du_{t}(e_{j})\big )h\big (du_{t}(e_{j}),du_{t}(e_{k}) \big )\nonumber\\
	&& \qquad -h\bigg (du_{t}(\frac{\partial}{\partial t}), \widetilde{\nabla} _{e_i}  \big ( h\big (du_{t}(e_{i}),du_{t}(e_{j})\big )h\big (du_{t}(e_{j}),du_{t}(e_{k}) \big )\, du_{t}(e_{k}) \big )\bigg )\nonumber\\
		&& \qquad + h\bigg (du_{t}(\frac{\partial}{\partial t}), e_i  \big ( h\big (du_{t}(e_{i}),du_{t}(e_{j})\big )\, h\big (du_{t}(e_{j}),du_{t}(e_{k}) \big )\, du_{t}(e_{k}) \big )\bigg )\nonumber\\	
		&& \qquad +h\bigg (du_{t}(\frac{\partial}{\partial t}), h\big (du_{t}(e_{i}),du_{t}(e_{j}) )\, e_i  \big (h\big (du_{t}(e_{j}), du_{t}(e_{k})\big ) \, du_{t}(e_{k}) \bigg )\nonumber\\
&& =\sum_{i,j,k=1}^m e_{i} \bigg (h\big (du_{t}(\frac{\partial}{\partial t}),du_{t}(e_{k})\big )h\big (du_{t}(e_{i}),du_{t}(e_{j})\big )h\big (du_{t}(e_{j}),du_{t}(e_{k}) \big )\bigg )\nonumber \\
&& \qquad-h\bigg (du_{t}(\frac{\partial}{\partial t}),\widetilde{\nabla}_{e_i} \big (d_{(3)}{u_t}\, {(e_i)}\big ) \bigg ),\nonumber
\end{eqnarray}
	where we use	
	\begin{equation*}
	\widetilde{\nabla}_{\frac{\partial}{\partial t}}\big (du_{t}(e_{i})\big )-\widetilde{\nabla}_{e_{i}}\bigg (du_{t}(\frac{\partial}{\partial t})\bigg )=du_{t}[\frac{\partial}{\partial t},e_{i}]=0
	\end{equation*}
	for the third equality, and \eqref{1.1}	for the last step. Let $X_{t}$ be a compactly supported vector field on $M$ such that
	$g(X_{t},Y)=\sum _{j,k=1} ^m h\big (du_{t}(\frac{\partial}{\partial t}),du_{t}(e_{k})\big )h\big (du_{t}(Y),du_{t}(e_{j})\big )h\big (du_{t}(e_{j}),du_{t}(e_{k})\big )$
	for any vector field $Y$ on $M.$ Then	
    \begin{equation*}
	\begin{aligned}  \frac{\partial}{\partial t}
	  \frac{|| d_{(3)} u_{t}||^2}{6}
	= & \sum_{i=1}^m\big (e_{i}g(X_{t},e_{i})\big )-h\bigg (du_{t}(\frac{\partial}{\partial t}),\widetilde{\nabla}_{e_i} \big (d_{(3)}{u_t}\, {(e_i)}\big )\bigg )\\
	= & \sum_{i=1}^m\big (g(\nabla_{e_{i}}X_{t},e_{i})+g(X_{t},\nabla_{e_{i}}e_{i})\big )-h\bigg (du_{t}(\frac{\partial}{\partial t}),\widetilde{\nabla}_{e_i} \big (d_{(3)}{u_t}\, {(e_i)}\big )\bigg )\\
	= & \operatorname{div}_{g}(X_{t})-h\bigg (du_{t}(\frac{\partial}{\partial t}),\tau_{\Phi_{(3)}} (u_t)\bigg ).
	\end{aligned}
	\end{equation*}
	
	By the Stokes' theorem, we have
	\begin{equation*}
	\begin{aligned}\frac{dE_{\Phi_{(3)}}(u_{t})}{dt}_{\big {|}_{t=0}} & =\int_{M}\begin{aligned}\frac{\partial}{\partial t} & \frac{|| d_{(3)} u_{t}||^2}{6}\end{aligned}_{\big {|}_{t=0}}\, dv_{g}\\
	& =-\int_{M}h\big (v,\tau_{\Phi_{(3)}} (u)\big )dv_{g}.
	\end{aligned}
	\end{equation*}
\end{proof}
\begin{prop}\label{prop:F1}
	A smooth map $u:(M^{m},g)\rightarrow(N,h)$ is $\Phi_{(3)}$-harmonic if and only if $u$ is a solution of the Euler-Lagrange equation for the $\Phi_{(3)}$-energy functional $E_{\Phi_{(3)}}$	
	\begin{equation}\label{eq:EL}
	\tau_{\Phi_{(3)}} (u) = 0.
	\end{equation}
\end{prop}

We introduce
\begin{defn}
	The stress-energy tensor $S_{\Phi_{(3)}}$ of $u$ with respect to the functional
	$E_{\Phi_{(3)}}(u)$ is the symmetric $2$-tensor on $M^{m}$ given by
	\begin{equation}
	S_{\Phi_{(3)}} =  e_{\Phi_{(3)}} g - (d_{(3)} u) ^{-1} h.
	\end{equation}
	That is,
	\begin{eqnarray}
	S_{\Phi_{(3)}}(X,Y)&=&\frac{||d_{(3)} u||^2}{6}g(X,Y) - h\big (d_{(3)} u(X), d u(Y)\big )\nonumber\\
	&=&\frac{1}{6}\sum_{i,j,k=1}^mh\big (du(e_i),du(e_j)\big )h\big (du(e_j),du(e_k)\big )h\big (du(e_k),du(e_i)\big )g(X,Y)\nonumber\\
	&&-\sum_{j,k=1}^mh\big (du(X),du(e_j)\big )h\big (du(e_j),du(e_k)\big )h\big (du(e_k),du(Y)\big )
	\end{eqnarray}
	for every smooth vector fields $X,Y$ on $M$.
\end{defn}
\begin{defn}\label{claw}
	A map $u:(M^m,g)\rightarrow(N,h)$ is said to satisfy the $\Phi_{(3)}$-conservation law
	if $S_{\Phi_{(3)}}$ is divergence free, i.e.
	\begin{equation}
	\operatorname {div}\, S_{\Phi_{(3)}} =0.\label{eq:F10}
	\end{equation}
\end{defn}
\begin{prop} \label{prop:01}
	For any smooth vector fields on $M$, we have	
	\begin{eqnarray}
	(\operatorname{div}\, S_{\Phi_{(3)}})(X)=-h\big (\tau_{\Phi_{(3)}} (u),du(X)\big ).
	\end{eqnarray}
\end{prop}
\begin{proof}
	We choose a local orthonormal frame field $\{e_{i}\}_{i=1}^{m}$ around a point
	$x_0$ in $M$ with $\nabla_{e_{i}} {e_{j}} _{\big {|}_{x_0}} = 0.$ Let $X$ be a vector
	field on $M.$ At $x_0,$ we have	\\
	\begin{equation*}
	\begin{aligned}
	&\quad\big (\operatorname{div}\, S_{\Phi_{(3)}}\big )(X)\\ & =\sum_{i=1}^{m}(\nabla_{e_{i}}S_{\Phi_{(3)}})(e_{i},X)\\
	& =\sum_{i=1}^{m}\bigg(e_{i}S_{\Phi_{(3)}}(e_{i},X)-S_{\Phi_{(3)}}\big (e_{i},\nabla_{e_{i}}X\big )\bigg)\\
	& =\sum_{i=1}^{m}\bigg(e_{i}\big(\frac{\Vert d_{(3)} u\Vert^{2}}{6}g(e_{i},X)\big)-e_{i}\big(h\big( d_{(3)} u(e_{i}),du(X)\big)\big)\\
	& \quad-\frac{\Vert d_{(3)} u\Vert^{2}}{6}g\big (e_{i},\nabla_{e_{i}}X\big )+h\big( d_{(3)} u(e_{i}),du(\nabla_{e_{i}}X)\big)\bigg)
\end{aligned}
\end{equation*}
\begin{equation*}
\begin{aligned}
	& =\sum_{i=1}^{m}\bigg(e_{i}\big(\frac{\Vert d_{(3)} u\Vert^{2}}{6}\big)g(e_{i},X)-h\big(\widetilde{\nabla}_{e_{i}} d_{(3)} u(e_{i}),du(X)\big)\\
	& \quad-h\big( d_{(3)} u(e_{i}),\widetilde{\nabla}_{e_{i}}du(X)\big)+h\big( d_{(3)} u(e_{i}),du(\nabla_{e_{i}}X)\big)\bigg)\\
	& =\sum_{i,j,k=1}^{m}\bigg(\frac{1}{6}X\big(h\big(du(e_{i}),du(e_{j})\big)h\big(du(e_{j}),du(e_{k})\big)h\big(du(e_{k}),du(e_{i})\big)\big)\\
	& \quad-h\big (\tau_{\Phi_{(3)}} (u),du(X)\big )-h\big( d_{(3)} u(e_{i}),(\nabla_{e_{i}}du)(X)\big)\bigg)\\
	& =\sum_{i,j,k=1}^{m}\bigg(h\big (du(e_{i}),du(e_{j})\big )h\big (du(e_{j}),du(e_{k})\big )h\big(du(e_{k}),\nabla_{X}du(e_{i})\big)\\
	& \quad-h\big (\tau_{\Phi_{(3)}} (u),du(X)\big )-h\big ( d_{(3)} u(e_{i}),(\nabla_{e_{i}}du)(X)\big )\bigg)\\
	& =\sum_{i=1}^{m}\bigg (h\big ( d_{(3)} u(e_{i}),(\nabla_{X}du)(e_{i})\big )-h\big ( d_{(3)} u(e_{i}),(\nabla_{e_{i}}du)(X)\big )\bigg )
	-h\big (\tau_{\Phi_{(3)}} (u),du(X)\big )\\
	& =-h\big (\tau_{\Phi_{(3)}} (u),du(X)\big ).
	\end{aligned}
	\end{equation*}
	The last equality holds due to
	\begin{equation*}
	(\nabla_{X}du)(e_{i})=(\nabla_{e_{i}}du)(X).
	\end{equation*}
This concludes the proof of Proposition \ref{prop:01}.	
\end{proof}
Let $X, Y\in\Gamma(TM)$. We denote
the dual one form of $X$ by $ X^{\flat}$, that is, $ X^{\flat} (Y)=g(X,Y)$. The covariant
derivative of $X^{\flat}$ gives a 2-tensor field $\nabla X^{\flat}:$
\begin{equation}
(\nabla X^{\flat})(Y,Z)=(\nabla_{Y} X^{\flat})(Z)=g(\nabla_{Y} X,Z).\label{eq:F6}
\end{equation}
If $X=\nabla f$ is the gradient field of some $C^{2}$ function
$f$ on $M,$ then $X^{\flat}=df$ and $\nabla X^{\flat}=\operatorname{Hess}\, f.$
\begin{lem}[\cite{han23,han12}]
	Let $X$ be a vector field and $T$ be a symmetric $(0,2)$-type
	tensor field. Then
	\begin{equation}
	\operatorname{div}(i_{X}T)=(\operatorname{div}T)(X)+\frac{1}{2}\langle T,L_{X}g\rangle,\label{eq:F5}
	\end{equation}
	where $L_{X}$ is the Lie derivative with respect to the direction
	$X.$
\end{lem}
Let $D$ be any bounded domain of $M$ with $C^{1}$ boundary and $\nu$
be the unit outward normal vector field along $\partial D.$ By using
the Stoke's theorem, we get
\begin{equation}
\int_{\partial D}T(X,\nu)ds_{g}=\int_{D}\big (\langle T,\frac{1}{2}L_{X}g\rangle+(\operatorname{div}T)(X)\big )dv_{g}.\label{eq:F1}
\end{equation}
According to (\ref{eq:F10}) and (\ref{eq:F1}), we have
\begin{equation}
\int_{\partial D}S_{\Phi_{(3)}}(X,\nu)ds_{g}=\int_{D}\langle S_{\Phi_{(3)}},\frac{1}{2}L_{X}g\rangle dv_{g}.\label{eq:F2}
\end{equation}

Let $X$ be a smooth vector field on $M$. We see that $u_{t}=u\circ\varphi_{t}$
is a deformation in Theorem \ref{Ft1}. Next we show the first variation
formula by 1-parameter families of diffeomorphisms.
\begin{thm}[The first variation formula $(II)$]
	\begin{equation*}
	\frac{d}{dt}E_{\Phi_{(3)}}(u_{t})_{\big {|}_{t=0}}=-\int_{M}\langle S_{\Phi_{(3)}},\frac{1}{2}L_{X}g\rangle dv_g,
	\end{equation*}
	where $L_{X}$ is the Lie derivative with respect to the direction
	$X.$
\end{thm}
\begin{proof}
	By the result of Theorem \ref{Ft1} and $u_{t}=u\circ\varphi_{t}$,
	we know that $du(X)$ is the variation vector field for the deformation
	$u_{t}.$
	\begin{equation}
	\begin{aligned}\frac{d}{dt}E_{\Phi_{(3)}}(u_{t})_{\big {|}_{t=0}} & =-\int_{M}\langle du(X),\tau_{\Phi_{(3)}} (u)\rangle dv_{g}\\
	& =\int_{M}\sum_{i=1}^{m}\langle\widetilde{\nabla}_{e_{i}}du(X), d_{(3)} u(e_{i})\rangle dv_{g}.
	\end{aligned}
	\end{equation}
	
	We denote a locally orthonormal frame around a fixed point $x_0$ on
	$M$ by $\{e_{i}\}$, such that  $\nabla_{e_{i}} {e_{j}} _{\big {|}_{x_0}} = 0.$ We have
	\begin{equation*}
	\begin{aligned}\langle S_{\Phi_{(3)}},\frac{1}{2}L_{X}g\rangle & =\frac{1}{2}\sum_{i,j}S_{\Phi_{(3)}}(e_{i},e_{j})L_{X}g(e_{i},e_{j})\\
	& =\sum_{i,j}S_{\Phi_{(3)}}(e_{i},e_{j})g(\nabla_{e_{i}}X,e_{j}).
	\end{aligned}
	\end{equation*}
	
	Hence, at $x_0$, we get
	\begin{equation*}
	\begin{aligned}
	\frac{d}{dt}E_{\Phi_{(3)}}(u_{t})_{\big {|}_{t=0}}  &=\int_{M}\sum_{i=1}^{m}\langle\widetilde{\nabla}_{e_{i}}du(X), d_{(3)} u(e_{i})\rangle dv_{g}\\
	&=\int_{M}\sum_{i=1}^{m}\bigg(h\big((\nabla_{e_{i}}du)(X), d_{(3)} u(e_{i})\big)+h\big(du(\nabla_{e_{i}}X), d_{(3)} u(e_{i})\big)\bigg)dv_{g}\\
	& =\int_{M}\sum_{i=1}^{m}\bigg(h\big ((\nabla_{X}du)(e_{i}), d_{(3)} u(e_{i})\big )+h\big (du(\nabla_{e_{i}}X), d_{(3)} u(e_{i})\big )\bigg)dv_{g}\\
	& =\int_{M}\sum_{i=1}^{m}\bigg(h\big (\widetilde{\nabla}_{X}du(e_{i}), d_{(3)} u(e_{i})\big )+h\big (du(\nabla_{e_{i}}X), d_{(3)} u(e_{i})\big )\bigg )dv_{g}\\
	& =\int_{M}\sum_{i=1}^{m}\bigg(\nabla_{X}(\frac{\Vert d_{(3)} u\Vert^{2}}{6})+h\big (du(\nabla_{e_{i}}X), d_{(3)} u(e_{i})\big )\bigg) dv_{g}\\
	& =\int_{M}\sum_{i=1}^{m}\bigg(L_{X}(\frac{\Vert d_{(3)} u\Vert^{2}}{6})+h\big (du(\nabla_{e_{i}}X), d_{(3)} u(e_{i})\big )\bigg )dv_{g}\\
	& =-\int_{M}\frac{\Vert d_{(3)} u\Vert^{2}}{6}L_{X}(dv_{g})+\int_{M}\sum_{i=1}^{m}h\big (du(\nabla_{e_{i}}X), d_{(3)} u(e_{i})\big )dv_{g}\\
	& =-\int_{M}\frac{\Vert d_{(3)} u\Vert^{2}}{6}\operatorname{div}Xdv_{g}+\int_{M}\sum_{i=1}^{m}h\big (du(\nabla_{e_{i}}X), d_{(3)} u(e_{i})\big )dv_{g}
\end{aligned}
\end{equation*}
\begin{equation*}
\begin{aligned}
	& =-\int_{M}\frac{\Vert d_{(3)} u\Vert^{2}}{6}\sum_{i,j=1}^{m}g(e_{i},e_{j})g\big (\nabla_{e_{i}}X,e_{j}\big )dv_{g}\\
	& \quad+\int_{M}\sum_{i,j=1}^{m}h\big (du(e_{j}), d_{(3)} u(e_{i})\big )g\big (\nabla_{e_{i}}X,e_{j}\big )dv_{g}\\
		& =-\int_{M}\sum_{i,j=1}^{m}S_{\Phi_{(3)}}(e_{i},e_{j})g\big (\nabla_{e_{i}}X,e_{j}\big )dv_{g}\\
	& =-\int_{M}\langle S_{\Phi_{(3)}},\frac{1}{2}L_{X}g\rangle dv_{g}.
	\end{aligned}
	\end{equation*}

Therefore, we have the desired result.
\end{proof}
\begin{thm}
	$(i)\quad$ If $u:(M^{m},g)\rightarrow(N,h)$ is a smooth $\Phi_{(3)}$-harmonic map, then $u$ satisfies the $\Phi_{(3)}$-conservation law.
	
	\noindent
	$(ii)\quad$ Conversely, if $u: (M^{m},g)\rightarrow(N,h)$ is a smooth map satisfying the $\Phi_{(3)}$-conservation law and $u$ is an immersion, then  $u$ is a smooth $\Phi_{(3)}$-harmonic map.
\end{thm}

\begin{rem}
	This is an analog and an extension of harmonic maps in which the stress-energy tensor unifies theory of harmonic maps. (cf. \cite {hy11})
\end{rem}

\begin{proof}
	According to Definition \ref{claw}, Proposition \ref{prop:F1} and Proposition \ref{prop:01}, we can obtain the desired results.
\end{proof}

\section{Liouville type results} \label{sec4}

We denote the $g$-distance function relative to the pole $x_{0}$
by $r(x),$ that is, $r(x)=dist_{g}(x,x_{0}).$ Denote $B(r)=\{x\in M^{m}:r(x)\leq  r\}.$
Obviously, $\frac{\partial}{\partial r}$ is an eigenvector of $\operatorname{Hess}_{g}(r^{2})$
with respect to eigenvalue 2. Let $\lambda_{\max}$$($resp. $\lambda_{\min}$
$)$ be the maximum $($resp. minimal$)$ eigenvalue of $\operatorname{Hess}_{g}(r^{2})-2dr\otimes dr$
at every point of $M\backslash\{x_{0}\}.$
\begin{thm}[Monotonicity Formula]\label{thm:L1}
	Let $u:(M^{m},g)\rightarrow(N^{n},h)$ be a smooth $\Phi_{(3)}$- harmonic map. If
	\begin{equation}
	1+\frac{m-1}{2}\lambda_{\min}-3\max\{2,\lambda_{\max}\}\geq\zeta,\label{eq:F7}
	\end{equation}
where $\zeta > 0$ is a constant, then we have
	\begin{equation}
	\frac{1}{\rho_{1}^{\zeta}} \int_{B(\rho_{1})} \frac {\Vert d_{(3)} u\Vert^{2}}{6} dv_{g}\leq\frac{1}{\rho_{2}^{\zeta}}\int_{B(\rho_{2})} \frac {\Vert d_{(3)} u\Vert^{2}}{6} dv_{g} \label{eq:F28}
	\end{equation}
	for any $0<\rho_{1}\leq\rho_{2}.$
\end{thm}
\begin{proof}
	We choose $D=B(r)$ and $X=r\frac{\partial}{\partial r}=\frac{1}{2}\nabla r^{2}$
	in (\ref{eq:F2}). Hence, we have
	\begin{equation*}
	\int_{\partial B(r)}S_{\Phi_{(3)}}\big (r\frac{\partial}{\partial r},\nu\big )ds_{g}=\int_{B(r)}\langle S_{\Phi_{(3)}},\frac{1}{2}L_{r\frac{\partial}{\partial r}}g\rangle dv_{g}.
	\end{equation*}
	
	From the Coarea formula, we get
	\begin{eqnarray}\label{eq:F3}
	&& \quad\int_{\partial B(r)}S_{\Phi_{(3)}}\big (r\frac{\partial}{\partial r},\nu\big )ds_{g}\nonumber\\
	&&= r\int_{\partial B(r)}\frac{\Vert d_{(3)} u\Vert^{2}}{6}ds_{g}-\int_{\partial B(r)}rh\big ( d_{(3)} u(\frac{\partial}{\partial r}),d u(\frac{\partial}{\partial r})\big )ds_{g}\nonumber\\
	&&= r\int_{\partial B(r)}\frac{\Vert d_{(3)} u\Vert^{2}}{6}ds_{g}-r\int_{\partial B(r)}\Vert d_{(3)} u(\frac{\partial}{\partial r})\Vert^{2}ds_{g}\nonumber\\
	 &&\leq  r\int_{\partial B(r)}\frac{\Vert d_{(3)} u\Vert^{2}}{6}ds_{g}\\
	&&= r\frac{d}{dr}\int_{B(r)}\frac{\Vert d_{(3)} u\Vert^{2}}{6}dv_{g}.\nonumber
	\end{eqnarray}
	
	Suppose $\{e_{i}\}_{i=1}^{m}$ is an orthonormal basis with $g$, such that $\operatorname{Hess}_{g}(r^{2})$ becomes a diagonal matrix and $e_{m}=\frac{\partial}{\partial r}$. By ($\text{\ref{eq:F7}}$), we have
	\begin{equation}\label{eq:F8}
	\begin{aligned}
	& \frac{1}{2}\langle S_{\Phi_{(3)}},L_{r\frac{\partial}{\partial r}}g\rangle\\
	= & \frac{1}{2}\sum_{i,j=1}^{m}S_{\Phi_{(3)}}(e_{i},e_{j})\big (L_{r\frac{\partial}{\partial r}}g\big )(e_{i},e_{j})\\
	= & \frac{1}{2}\sum_{i,j=1}^{m}\bigg( \frac{\Vert d_{(3)} u\Vert^{2}}{6}g(e_{i},e_{j})-h\big (d_{(3)} u(e_{i}), du(e_{j})\big )\bigg ) \big (L_{r\frac{\partial}{\partial r}}g\big )(e_{i},e_{j})\\
%
	= & \frac{1}{2}\bigg( \sum_{i=1}^{m}\frac{\Vert d_{(3)} u\Vert^{2}}{6}\operatorname{Hess}_{g}(r^{2})(e_{i},e_{i})-\sum_{i,j=1}^{m}h\big (d_{(3)} u(e_{i}), du(e_{j})\big )\operatorname{Hess}_{g}(r^{2})(e_{i},e_{j})\bigg ) \\
	\geq & \frac{1}{2}\bigg( \frac{\Vert d_{(3)} u\Vert^{2}}{6}\bigg( 2+(m-1)\lambda_{\min}\bigg ) -\sum_{i=1}^{m}\max\{2,\lambda_{\max}\}h\big (d_{(3)} u(e_{i}), du(e_{i})\big )\bigg ) \\
	= & \frac{\Vert d_{(3)} u\Vert^{2}}{6}\bigg( 1+\frac{m-1}{2}\lambda_{\min}-3\max\{2,\lambda_{\max}\}\bigg ) \\
	\geq & \zeta\frac{\Vert d_{(3)} u\Vert^{2}}{6}.
	\end{aligned}
	\end{equation}
	Using (\ref{eq:F3}) and (\ref{eq:F8}), we get
	\begin{equation*}
	r\frac{d}{dr}\int_{B(r)}\frac{\Vert d_{(3)} u\Vert^{2}}{6}dv_{g}\geq\zeta\int_{B(r)}\frac{\Vert d_{(3)} u\Vert^{2}}{6}dv_{g}.
	\end{equation*}
	
	Thus, integrating with respect to $r$ over the interval $[\rho_1, \rho_2]$
\begin{equation*}
	\frac {\frac{d}{dr}\int_{B(r)}\frac{\Vert d_{(3)} u\Vert^{2}}{6}dv_{g}}{\int_{B(r)}\frac{\Vert d_{(3)} u\Vert^{2}}{6}dv_{g}}\geq \frac {\zeta}{r}
	\end{equation*}
	by the fundamental of theorem of Calculus, we obtain the desired monotonicity formula (\ref{eq:F28}).
	\end{proof}
	
From Theorem \ref{thm:L1}, we have
\begin{prop}
	\label{prop:L2}Suppose $u:(M^{m},g)\rightarrow(N^{n},h)$ is a smooth
	$\Phi_{(3)}$-harmonic map and $r(x)$ satisfies the condition
	(\ref{eq:F7}). If $R\rightarrow\infty$ and $u$ is a nonconstant map, then
	we have
	\begin{equation}
	\int_{B(R)}\frac{\Vert d_{(3)} u\Vert^{2}}{6}dv_{g}\geq c(u)R^{\zeta},\label{eq:L1}
	\end{equation}
	where $c(u)>0$ is a constant.
\end{prop}
In local coordinates, we obtain
	\begin{equation*}\label{eq:L3}
	\begin{aligned}
	&E_{\Phi_{(3)}}^{R}(u)=\int_{B(R)}\frac{\Vert d_{(3)} u\Vert^{2}}{6}dv_{g}\nonumber\\
	&\quad\quad\quad~~=\int_{B(R)}\frac{1}{6}\sum_{i,j,k}h\big (du(e_{i}),du(e_{j})\big )h\big (du(e_{j}),du(e_{k})\big )h\big (du(e_{k}),du(e_{i})\big )dv_{g}\nonumber\\
	&\quad\quad\quad~~ =\int_{B(R)}\frac{1}{6}\sum_{i,j,k,a,b,c,\alpha,\beta,\gamma,\delta,\mu,\nu} g^{ia}g^{jb}g^{kc}h_{\alpha\beta}\frac{\partial u^{\alpha}}{\partial x_{i}}\frac{\partial u^{\beta}}{\partial x_{j}} h_{\gamma\delta}\frac{\partial u^{\gamma}}{\partial x_{b}}\frac{\partial u^{\delta}}{\partial x_{k}}h_{\mu\nu}\frac{\partial u^{\mu}}{\partial x_{c}}\frac{\partial u^{\nu}}{\partial x_{a}}dv_{g}.\nonumber
	\end{aligned}
	\end{equation*}

Extending the method  employed in Jin \cite{hy04} for $\Phi_{(1)}$-harmonic maps (i.e. harmonic maps) on conformally Euclidean domain to $\Phi_{(3)}$-harmonic maps on Riemannian manifolds, we obtain
\begin{prop}
	\label{prop:L3}Suppose $u:(M^{m},g)\rightarrow(N,h)$ is a smooth
	$\Phi_{(3)}$-harmonic map and $r(x)$ satisfies the condition
	(\ref{eq:F7}). When $R$ is large enough,
	\begin{equation*}
	\int_{R}^{\infty}\frac{1}{\operatorname{Vol}_{g}^{\frac{1}{5}}\big (\partial B(r)\big )}dr\geq CR^{-\frac{\zeta}{5}}.
	\end{equation*}
	If $u(x)\rightarrow p_{0}\in N$ as $r(x)\rightarrow\infty$,
	then we have
	\begin{equation}
	E_{\Phi_{(3)}}^{R}(u) = o(R^{\zeta}), \quad \operatorname{as}\quad R \to \infty\label{eq:L0}
	\end{equation}
	
\end{prop}
\begin{proof}
	Suppose $u$ is a nonconstant map. From Proposition $\ref{prop:L2}$, we get
	$E_{\Phi_{(3)}}^{R}(u)\rightarrow\infty$ as $R\rightarrow\infty.$
	We choose a local coordinate neighbourhood $(U, \varphi)$ of $p_{0}$ in $N$,
	such that $\varphi(p_{0})=0$ and
	\begin{equation*}
	h=h_{\alpha\beta}(y)dy^{\alpha}\otimes dy^{\beta},\quad y\in U,
	\end{equation*}
    satisfies
	\begin{equation*}
	\bigg(\frac{\partial h_{\alpha\beta}(y)}{\partial y^{\gamma}}y^{\gamma}+2h_{\alpha\beta}(y)\bigg )\geq\big (h_{\alpha\beta}(y)\big )\quad \text{on} \quad U
	\end{equation*}
	 in the matrice sence (that is, for two $n\times n$ matrices $A,B,$ by $A \geq B$, we mean that  $A-B$
	is a positive semi-definite matrix). Since $u(x)\rightarrow0$ as $r(x)\rightarrow\infty,$ there
	exists a $R_{1}$ such that for $r(x)>R_{1},$ $u(x)\in U$ and
	\begin{equation}
	\bigg(\frac{\partial h_{\alpha\beta}(u)}{\partial y^{\gamma}}u^{\gamma}+2h_{\alpha\beta}(u)\bigg )\geq\big (h_{\alpha\beta}(u)\big ).\label{eq:L3}
	\end{equation}
	For $w\in C_{0}^{2}(M\backslash B(R_{1}),\varphi(U))$ and sufficiently small $t$,
	the variation $u+tw:M^{m}\rightarrow N$ is defined by
	\begin{equation*}
	(u+tw)(q)=\begin{cases}
	u(q) & \quad q\in B(R_{1});\\
	\varphi^{-1}\big (\varphi(u)+tw\big )(q) & \quad q\in M\backslash B(R_{1}).
	\end{cases}
	\end{equation*}
	 By \eqref{FV} and (\ref{eq:EL}), we have
	\begin{equation*}
	\frac{d}{dt}E_{\Phi_{(3)}}(u+tw)_{\big {|}_{t=0}}=0,
	\end{equation*}
	that is, using Einstein notation, we obtain,
	\begin{eqnarray}\label{eq:L4}
	 \int_{M\backslash B(R_{1})}g^{is}g^{jk}g^{lr}\bigg (\frac{\partial h_{\alpha\beta}}{\partial y^{\xi}}w^{\xi}\frac{\partial u^{\alpha}}{\partial x_{i}}\frac{\partial u^{\beta}}{\partial x_{j}}+2h_{\alpha\beta}\frac{\partial u^{\alpha}}{\partial x_{i}}\frac{\partial w^{\beta}}{\partial x_{j}}\bigg )
	 h_{\gamma\delta}\frac{\partial u^{\gamma}}{\partial x_{k}}\frac{\partial u^{\delta}}{\partial x_{l}}h_{\mu\nu}\frac{\partial u^{\mu}}{\partial x_{r}}\frac{\partial u^{\nu}}{\partial x_{s}}dv_{g}=0.
	\end{eqnarray}
    Let $\phi(t)$ be a smooth function on $(R_{1},\infty)$. Choosing $w(x)=\phi\big (r(x)\big )u(x)$ in ($\ref{eq:L4}$), we have
	\begin{eqnarray}\label{eq:L5}
	&&\quad\int_{M\backslash B(R_{1})}g^{is}g^{jk}g^{lr}\bigg (\frac{\partial h_{\alpha\beta}}{\partial y^{\xi}}u^{\xi}+2h_{\alpha\beta}\bigg )\frac{\partial u^{\alpha}}{\partial x_{i}}\frac{\partial u^{\beta}}{\partial x_{j}}\phi\big (r(x)\big )  h_{\gamma\delta}\frac{\partial u^{\gamma}}{\partial x_{k}}\frac{\partial u^{\delta}}{\partial x_{l}}h_{\mu\nu}\frac{\partial u^{\mu}}{\partial x_{r}}\frac{\partial u^{\nu}}{\partial x_{s}}dv_{g}\nonumber\\
	&&=-2\int_{M\backslash B(R_{1})}g^{is}g^{jk}g^{lr}h_{\alpha\beta}\frac{\partial u^{\alpha}}{\partial x_{i}}\frac{\partial\phi\big (r(x)\big )}{\partial x_{j}}u^{\beta}  h_{\gamma\delta}\frac{\partial u^{\gamma}}{\partial x_{k}}\frac{\partial u^{\delta}}{\partial x_{l}}h_{\mu\nu}\frac{\partial u^{\mu}}{\partial x_{r}}\frac{\partial u^{\nu}}{\partial x_{s}}dv_{g}.
	\end{eqnarray}
	The above equation holds for Lipschitz function $\phi$ with compact
	support by an approximation argument.
	
	For $0<\varepsilon\leq 1,$ we define
	\begin{equation*}
	\varphi_{\varepsilon}(t)=\begin{cases}
	1 & t\leq 1;\\
	1+\frac{1-t}{\varepsilon} & 1<t<1+\varepsilon;\\
	0 & t\geq1+\varepsilon
	\end{cases}
	\end{equation*}
	and choose
	\begin{equation*}
	\phi\big (r(x)\big )=\varphi_{\varepsilon}\bigg(\frac{r(x)}{R}\bigg ) \bigg(1-\varphi_{1}\bigg(\frac{r(x)}{R_{1}}\bigg )\bigg ),
   \end{equation*}
	where $R>2R_{1}.$ Let $R_{2}=2R_{1}$. In (\ref{eq:L5}), we obtain
	\begin{eqnarray}\label{eq:L6}
	 &&\quad\int_{M\backslash B(R_{1})}g^{is}g^{jk}g^{lr}\bigg (\frac{\partial h_{\alpha\beta}}{\partial y^{\xi}}u^{\xi}+2h_{\alpha\beta}\bigg )
	 \frac{\partial u^{\alpha}}{\partial x_{i}}\frac{\partial u^{\beta}}{\partial x_{j}}\phi\big (r(x)\big )
     h_{\gamma\delta}\frac{\partial u^{\gamma}}{\partial x_{k}}\frac{\partial u^{\delta}}{\partial x_{l}}h_{\mu\nu}\frac{\partial u^{\mu}}{\partial x_{r}}\frac{\partial u^{\nu}}{\partial x_{s}}dv_{g}\nonumber\\
	&&=\int_{B(R_{2})\backslash B(R_{1})}g^{is}g^{jk}g^{lr}\bigg (\frac{\partial h_{\alpha\beta}}{\partial y^{\xi}}u^{\xi}+2h_{\alpha\beta}\bigg )\frac{\partial u^{\alpha}}{\partial x_{i}}\frac{\partial u^{\beta}}{\partial x_{j}}\bigg(1-\varphi_{1}\bigg(\frac{r(x)}{R_{1}}\bigg )\bigg )\nonumber\\
	&&\quad\times \quad h_{\gamma\delta}\frac{\partial u^{\gamma}}{\partial x_{k}}\frac{\partial u^{\delta}}{\partial x_{l}}h_{\mu\nu}\frac{\partial u^{\mu}}{\partial x_{r}}\frac{\partial u^{\nu}}{\partial x_{s}}dv_{g}\\
	&&+\int_{B(R)\backslash B(R_{2})}g^{is}g^{jk}g^{lr}\bigg (\frac{\partial h_{\alpha\beta}}{\partial y^{\xi}}u^{\xi}+2h_{\alpha\beta}\bigg )\frac{\partial u^{\alpha}}{\partial x_{i}}\frac{\partial u^{\beta}}{\partial x_{j}}
	h_{\gamma\delta}\frac{\partial u^{\gamma}}{\partial x_{k}}\frac{\partial u^{\delta}}{\partial x_{l}}h_{\mu\nu}\frac{\partial u^{\mu}}{\partial x_{r}}\frac{\partial u^{\nu}}{\partial x_{s}}dv_{g}\nonumber\\
	&&+\int_{B((1+\varepsilon)R)\backslash B(R)}g^{is}g^{jk}g^{lr}\bigg (
	\frac{\partial h_{\alpha\beta}}{\partial y^{\xi}}u^{\xi}+2h_{\alpha\beta}\bigg )\frac{\partial u^{\alpha}}{\partial x_{i}}\frac{\partial u^{\beta}}{\partial x_{j}}\varphi_{\varepsilon}\bigg(\frac{r(x)}{R}\bigg )
	h_{\gamma\delta}\frac{\partial u^{\gamma}}{\partial x_{k}}\frac{\partial u^{\delta}}{\partial x_{l}}h_{\mu\nu}\frac{\partial u^{\mu}}{\partial x_{r}}\frac{\partial u^{\nu}}{\partial x_{s}}dv_{g}\nonumber
	\end{eqnarray}
	and
	\begin{eqnarray}\label{eq:L7}
	&&\quad-2\int_{M\backslash B(R_{1})}g^{is}g^{jk}g^{lr}h_{\alpha\beta}\frac{\partial u^{\alpha}}{\partial x_{i}}\frac{\partial\phi\big (r(x)\big )}{\partial x_{j}}u^{\beta}
	h_{\gamma\delta}\frac{\partial u^{\gamma}}{\partial x_{k}}\frac{\partial u^{\delta}}{\partial x_{l}}h_{\mu\nu}\frac{\partial u^{\mu}}{\partial x_{r}}\frac{\partial u^{\nu}}{\partial x_{s}}dv_{g}\nonumber\\
	&&= 2\int_{B(R_{2})\backslash B(R_{1})}g^{is}g^{jk}g^{lr}h_{\alpha\beta}\frac{\partial u^{\alpha}}{\partial x_{i}}\frac{\partial\varphi_{1}\bigg(\frac{r(x)}{R_{1}}\bigg )}{\partial x_{j}}u^{\beta}
	h_{\gamma\delta}\frac{\partial u^{\gamma}}{\partial x_{k}}\frac{\partial u^{\delta}}{\partial x_{l}}h_{\mu\nu}\frac{\partial u^{\mu}}{\partial x_{r}}\frac{\partial u^{\nu}}{\partial x_{s}}dv_{g}\nonumber\\
	&& \quad-2\int_{B((1+\varepsilon)R)\backslash B(R)}g^{is}g^{jk}g^{lr}h_{\alpha\beta}\frac{\partial u^{\alpha}}{\partial x_{i}}\frac{\partial\varphi_{\varepsilon}\bigg(\frac{r(x)}{R}\bigg )}{\partial x_{j}}u^{\beta}
	h_{\gamma\delta}\frac{\partial u^{\gamma}}{\partial x_{k}}\frac{\partial u^{\delta}}{\partial x_{l}}h_{\mu\nu}\frac{\partial u^{\mu}}{\partial x_{r}}\frac{\partial u^{\nu}}{\partial x_{s}}dv_{g}\nonumber\\
	&&= 2\int_{B(R_{2})\backslash B(R_{1})}g^{is}g^{jk}g^{lr}h_{\alpha\beta}\frac{\partial u^{\alpha}}{\partial x_{i}}\frac{\partial\varphi_{1}\bigg(\frac{r(x)}{R_{1}}\bigg )}{\partial x_{j}}u^{\beta}
	h_{\gamma\delta}\frac{\partial u^{\gamma}}{\partial x_{k}}\frac{\partial u^{\delta}}{\partial x_{l}}h_{\mu\nu}\frac{\partial u^{\mu}}{\partial x_{r}}\frac{\partial u^{\nu}}{\partial x_{s}}dv_{g}\\
	&& \quad+2\frac{1}{R\varepsilon}\int_{B((1+\varepsilon)R)\backslash B(R)}g^{is}g^{jk}g^{lr}h_{\alpha\beta}\frac{\partial u^{\alpha}}{\partial x_{i}}\frac{\partial r(x)}{\partial x_{j}}u^{\beta}
	h_{\gamma\delta}\frac{\partial u^{\gamma}}{\partial x_{k}}\frac{\partial u^{\delta}}{\partial x_{l}}h_{\mu\nu}\frac{\partial u^{\mu}}{\partial x_{r}}\frac{\partial u^{\nu}}{\partial x_{s}}dv_{g}.\nonumber
	\end{eqnarray}
	As $\varepsilon\rightarrow0,$
	\begin{equation}
	\begin{aligned} & 2\frac{1}{R\varepsilon}\int_{B((1+\varepsilon)R)\backslash B(R)}g^{is}g^{jk}g^{lr}h_{\alpha\beta}\frac{\partial u^{\alpha}}{\partial x_{i}}\frac{\partial r(x)}{\partial x_{j}}u^{\beta}
	 h_{\gamma\delta}\frac{\partial u^{\gamma}}{\partial x_{k}}\frac{\partial u^{\delta}}{\partial x_{l}}h_{\mu\nu}\frac{\partial u^{\mu}}{\partial x_{r}}\frac{\partial u^{\nu}}{\partial x_{s}}dv_{g}\\
	\to & 2\int_{\partial B(R)}g^{is}g^{jk}g^{lr}h_{\alpha\beta}\frac{\partial u^{\alpha}}{\partial x_{i}}\frac{\partial r(x)}{\partial x_{j}}u^{\beta}
	h_{\gamma\delta}\frac{\partial u^{\gamma}}{\partial x_{k}}\frac{\partial u^{\delta}}{\partial x_{l}}h_{\mu\nu}\frac{\partial u^{\mu}}{\partial x_{r}}\frac{\partial u^{\nu}}{\partial x_{s}}ds_{g}.
	\end{aligned}
	\label{eq:L8}
	\end{equation}
	From ($\ref{eq:L5}$)-($\ref{eq:L8}$), we get\\
	\begin{eqnarray}\label{eq:L9}
	 && \quad\int_{B(R)\backslash B(R_{2})}g^{is}g^{jk}g^{lr}\bigg (\frac{\partial h_{\alpha\beta}}{\partial y^{\xi}}u^{\xi}+2h_{\alpha\beta}\bigg )\frac{\partial u^{\alpha}}{\partial x_{i}}\frac{\partial u^{\beta}}{\partial x_{j}}
	h_{\gamma\delta}\frac{\partial u^{\gamma}}{\partial x_{k}}\frac{\partial u^{\delta}}{\partial x_{l}}h_{\mu\nu}\frac{\partial u^{\mu}}{\partial x_{r}}\frac{\partial u^{\nu}}{\partial x_{s}}dv_{g}\nonumber\\
	&&= 2\int_{\partial B(R)}g^{is}g^{jk}g^{lr}h_{\alpha\beta}\frac{\partial u^{\alpha}}{\partial x_{i}}\frac{\partial r(x)}{\partial x_{j}}u^{\beta}
	h_{\gamma\delta}\frac{\partial u^{\gamma}}{\partial x_{k}}\frac{\partial u^{\delta}}{\partial x_{l}}h_{\mu\nu}\frac{\partial u^{\mu}}{\partial x_{r}}\frac{\partial u^{\nu}}{\partial x_{s}}ds_{g}-D(R_{1}),
	\end{eqnarray}
	where
	\begin{equation*}
	\begin{aligned}D(R_{1}) &=-2\int_{B(R_{2})\backslash B(R_{1})}g^{is}g^{jk}g^{lr}h_{\alpha\beta}\frac{\partial u^{\alpha}}{\partial x_{i}}\frac{\partial\varphi_{1}(\frac{r(x)}{R_{1}})}{\partial x_{j}}u^{\beta}
	h_{\gamma\delta}\frac{\partial u^{\gamma}}{\partial x_{k}}\frac{\partial u^{\delta}}{\partial x_{l}}h_{\mu\nu}\frac{\partial u^{\mu}}{\partial x_{r}}\frac{\partial u^{\nu}}{\partial x_{s}}dv_{g}\\
	&\quad+\int_{B(R_{2})\backslash B(R_{1})}g^{is}g^{jk}g^{lr}\bigg (
	\frac{\partial h_{\alpha\beta}}{\partial y^{\xi}}u^{\xi}+2h_{\alpha\beta}\bigg )
	\frac{\partial u^{\alpha}}{\partial x_{i}}\frac{\partial u^{\beta}}{\partial x_{j}}\bigg(1-\varphi_{1}\bigg(\frac{r(x)}{R_{1}}\bigg )\bigg )\\
	&\quad\times\quad h_{\gamma\delta}\frac{\partial u^{\gamma}}{\partial x_{k}}\frac{\partial u^{\delta}}{\partial x_{l}}h_{\mu\nu}\frac{\partial u^{\mu}}{\partial x_{r}}\frac{\partial u^{\nu}}{\partial x_{s}}dv_{g}.
	\end{aligned}
    \end{equation*}
	
	Next we estimate the term
	
	\begin{equation}\label{eq:L11}
	 2\int_{\partial B(R)}g^{is}g^{jk}g^{lr}h_{\alpha\beta}\frac{\partial u^{\alpha}}{\partial x_{i}}\frac{\partial r(x)}{\partial x_{j}}u^{\beta} h_{\gamma\delta}\frac{\partial u^{\gamma}}{\partial x_{k}}\frac{\partial u^{\delta}}{\partial x_{l}}h_{\mu\nu}\frac{\partial u^{\mu}}{\partial x_{r}}\frac{\partial u^{\nu}}{\partial x_{s}}ds_{g}.
	\end{equation}
	Since the integrand in ($\ref{eq:L11}$) does not depend on the coordinate
	systems on $M^{m}$ and $N^{n}$, at any point $p\in\partial B(R)$
	and $u(p)$, we take the appropriate coordinate systems such that $g_{ij}(p)=\delta_{ij},$
	$g^{ij}(p)=\delta^{ij}$ and $h_{\alpha\beta}\big (u(p)\big )=\delta_{\alpha\beta}.$
	At $p,$ by H\"older's inequality, we have

\begin{equation*}
\begin{aligned}
&\quad g^{is}g^{jk}g^{lr}h_{\alpha\beta}\frac{\partial u^{\alpha}}{\partial x_{i}}\frac{\partial r(x)}{\partial x_{j}}u^{\beta}h_{\gamma\delta}\frac{\partial u^{\gamma}}{\partial x_{k}}\frac{\partial u^{\delta}}{\partial x_{l}}h_{\mu\nu}\frac{\partial u^{\mu}}{\partial x_{r}}\frac{\partial u^{\nu}}{\partial x_{s}}\nonumber \\
&=\sum_{i,j,l}\bigg(\sum_{\alpha}\frac{\partial u^{\alpha}}{\partial x_{i}}u^{\alpha}\frac{\partial r(x)}{\partial x_{j}}\bigg )\bigg(\sum_{\gamma,\mu}\frac{\partial u^{\gamma}}{\partial x_{j}}\frac{\partial u^{\gamma}}{\partial x_{l}}\frac{\partial u^{\mu}}{\partial x_{l}}\frac{\partial u^{\mu}}{\partial x_{i}}\bigg )\nonumber \\
&=\sum_{i,j}\bigg(\sum_{\alpha}\frac{\partial u^{\alpha}}{\partial x_{i}}u^{\alpha}\frac{\partial r(x)}{\partial x_{j}}\bigg )\bigg(\sum_{l} \bigg(\sum_{\gamma}\frac{\partial u^{\gamma}}{\partial x_{j}}\frac{\partial u^{\gamma}}{\partial x_{l}} \bigg )  \bigg(\sum_{\mu}\frac{\partial u^{\mu}}{\partial x_{l}}\frac{\partial u^{\mu}}{\partial x_{i}}\bigg )\bigg )\nonumber \\
&\leq \sum_{i,j}\bigg(\sum_{\alpha}\frac{\partial u^{\alpha}}{\partial x_{i}}u^{\alpha}\frac{\partial r(x)}{\partial x_{j}}\bigg ) \bigg( \sum_{l} \bigg(\sum_{\gamma} \frac{\partial u^{\gamma}}{\partial x_{j}}\frac{\partial u^{\gamma}}{\partial x_{l}} \bigg )^2
\bigg )^{\frac{1}{2}}
\bigg( \sum_{l} \bigg( \sum_{\mu}\frac{\partial u^{\mu}}{\partial x_{l}}\frac{\partial u^{\mu}}{\partial x_{i}} \bigg )^2
\bigg )^{\frac{1}{2}}\nonumber \\
&\leq  \bigg (\sum_{i,j}\bigg(\sum_{\alpha}\frac{\partial u^{\alpha}}{\partial x_{i}}u^{\alpha}\frac{\partial r(x)}{\partial x_{j}}\bigg )^{2}\bigg )^{\frac{1}{2}}
\bigg (\sum_{j,l}\bigg(\sum_{\gamma}\frac{\partial u^{\gamma}}{\partial x_{j}}\frac{\partial u^{\gamma}}{\partial x_{l}}\bigg )^{2}\bigg )^{\frac{1}{2}}\bigg (\sum_{i,l}\bigg(\sum_{\mu}\frac{\partial u^{\mu}}{\partial x_{l}}\frac{\partial u^{\mu}}{\partial x_{i}}\bigg )^{2}\bigg )^{\frac{1}{2}}\nonumber \\
&=\bigg(\sum_{i,j}\bigg (\sum_{\alpha}\frac{\partial u^{\alpha}}{\partial x_{i}}u^{\alpha}\bigg )^{2}\bigg (\frac{\partial r(x)}{\partial x_{j}}\bigg )^{2}\bigg )^{\frac{1}{2}}\bigg (\sum_{j,l}\bigg(\sum_{\gamma}\frac{\partial u^{\gamma}}{\partial x_{j}}\frac{\partial u^{\gamma}}{\partial x_{l}}\bigg )^{2}\bigg )^{\frac{1}{2}}\bigg (\sum_{i,l}\bigg(\sum_{\mu}\frac{\partial u^{\mu}}{\partial x_{l}}\frac{\partial u^{\mu}}{\partial x_{i}}\bigg )^{2}\bigg )^{\frac{1}{2}}\nonumber \\
&\leq \bigg (\sum_{i}\bigg(\sum_{\alpha}\frac{\partial u^{\alpha}}{\partial x_{i}}u^{\alpha}\bigg )^{2}\bigg )^{\frac{1}{2}}\bigg (\sum_{j,l}\bigg(\sum_{\gamma}\frac{\partial u^{\gamma}}{\partial x_{j}}\frac{\partial u^{\gamma}}{\partial x_{l}}\bigg )^{2}\bigg )^{\frac{1}{2}}\bigg (\sum_{i,l}\bigg(\sum_{\mu}\frac{\partial u^{\mu}}{\partial x_{l}}\frac{\partial u^{\mu}}{\partial x_{i}}\bigg )^{2}\bigg )^{\frac{1}{2}}\nonumber \\
&\leq \bigg (\sum_{i}\bigg(\sum_{\alpha}\frac{\partial u^{\alpha}}{\partial x_{i}}\frac{\partial u^{\alpha}}{\partial x_{i}}\bigg )\bigg )^{\frac{1}{2}}\bigg (\sum_{\alpha}\big (u^{\alpha}\big )^{2}\bigg )^{\frac{1}{2}}\bigg (\sum_{j,l}\bigg(\sum_{\gamma}\frac{\partial u^{\gamma}}{\partial x_{j}}\frac{\partial u^{\gamma}}{\partial x_{l}}\bigg )^{2}\bigg )^{\frac{1}{2}}\bigg (\sum_{i,l}\bigg(\sum_{\mu}\frac{\partial u^{\mu}}{\partial x_{l}}\frac{\partial u^{\mu}}{\partial x_{i}}\bigg )^{2}\bigg )^{\frac{1}{2}}\nonumber
\end{aligned}
\end{equation*}
\begin{equation*}
\begin{aligned}
&\leq  m^{\frac{1}{4}}\bigg (\sum_{i}\bigg(\sum_{\alpha}\frac{\partial u^{\alpha}}{\partial x_{i}}\frac{\partial u^{\alpha}}{\partial x_{i}}\bigg )^{2}\bigg )^{\frac{1}{4}}\bigg (\sum_{\alpha}\big (u^{\alpha}\big )^{2}\bigg )^{\frac{1}{2}}\bigg (\sum_{j,l}\bigg(\sum_{\gamma}\frac{\partial u^{\gamma}}{\partial x_{j}}\frac{\partial u^{\gamma}}{\partial x_{l}}\bigg )^{2}\bigg )^{\frac{1}{2}}\nonumber \\
&\quad\times \bigg (\sum_{i,l}\bigg(\sum_{\mu}\frac{\partial u^{\mu}}{\partial x_{l}}\frac{\partial u^{\mu}}{\partial x_{i}}\bigg )^{2}\bigg )^{\frac{1}{2}}\nonumber\\
&=m^{\frac{1}{4}} \bigg (\sum_{i}\bigg(h\big (du(\frac{\partial}{\partial x_{i}}),du(\frac{\partial}{\partial x_{i}})\big )\bigg )^{2}\bigg )^{\frac{1}{4}}\bigg (\sum_{\alpha}\big (u^{\alpha}\big )^{2}\bigg )^{\frac{1}{2}}\bigg (\sum_{j,l}\bigg(\sum_{\gamma}\frac{\partial u^{\gamma}}{\partial x_{j}}\frac{\partial u^{\gamma}}{\partial x_{l}}\bigg )^{2}\bigg )^{\frac{1}{2}}\nonumber\\
&\quad\times\bigg (\sum_{i,l}\bigg(\sum_{\mu}\frac{\partial u^{\mu}}{\partial x_{l}}\frac{\partial u^{\mu}}{\partial x_{i}}\bigg )^{2}\bigg )^{\frac{1}{2}}\nonumber\\
&\leq  m^{\frac{1}{4}}\bigg (\sum_{i,j}\bigg(h\big (du(\frac{\partial}{\partial x_{i}}),du(\frac{\partial}{\partial x_{j}})\big )\bigg )^{2}\bigg )^{\frac{1}{4}}\bigg (\sum_{\alpha}\big (u^{\alpha}\big )^{2}\bigg )^{\frac{1}{2}}\nonumber \\
&\quad\times\bigg (\sum_{j,l}\bigg(h\big (du(\frac{\partial}{\partial x_{j}}),du(\frac{\partial}{\partial x_{l}})\big )\bigg )^{2}\bigg )^{\frac{1}{2}}\bigg (\sum_{i,l}\bigg(h\big (du(\frac{\partial}{\partial x_{l}}),du(\frac{\partial}{\partial x_{i}})\big )\bigg )^{2}\bigg )^{\frac{1}{2}}\nonumber\\
&\leq  m^{\frac{5}{12}}\bigg (\sum_{\alpha}\big (u^{\alpha}\big )^{2}\bigg )^{\frac{1}{2}}\bigg (\sum_{i,j}\bigg(h\big (du(\frac{\partial}{\partial x_{i}}),du(\frac{\partial}{\partial x_{j}})\big )\bigg )^{3}\bigg )^{\frac{1}{6}}\nonumber \\
&\quad\times\bigg (\sum_{j,l}\bigg(h\big (du(\frac{\partial}{\partial x_{j}}),du(\frac{\partial}{\partial x_{l}})\big )\bigg )^{2}\bigg )^{\frac{1}{2}}\bigg (\sum_{i,l}\bigg(h\big (du(\frac{\partial}{\partial x_{l}}),du(\frac{\partial}{\partial x_{i}})\big )\bigg )^{2}\bigg )^{\frac{1}{2}}\nonumber \\
&\leq  m^{\frac{13}{12}}\bigg (\sum_{\alpha}\big (u^{\alpha}\big )^{2}\bigg )^{\frac{1}{2}}\bigg (\sum_{i,j}\bigg(h\big (du(\frac{\partial}{\partial x_{i}}),du(\frac{\partial}{\partial x_{j}})\big )\bigg )^{3}\bigg )^{\frac{1}{6}}\nonumber \\
&\quad\times\bigg (\sum_{j,l}\bigg(h\big (du(\frac{\partial}{\partial x_{j}}),du(\frac{\partial}{\partial x_{l}})\big )\bigg )^{3}\bigg )^{\frac{1}{3}} \bigg (\sum_{i,l}\bigg(h\big (du(\frac{\partial}{\partial x_{l}}),du(\frac{\partial}{\partial x_{i}})\big )\bigg )^{3}\bigg )^{\frac{1}{3}}\nonumber \\
&\leq  m^{\frac{13}{12}}\bigg (\sum_{\alpha}\big (u^{\alpha}\big )^{2}\bigg )^{\frac{1}{2}}\bigg
(\sum_{i,j,l}h\big (du(\frac{\partial}{\partial x_{i}}),du(\frac{\partial}{\partial x_{j}})\big )h\big (du(\frac{\partial}{\partial x_{j}}),du(\frac{\partial}{\partial x_{l}})\big )h\big (du(\frac{\partial}{\partial x_{l}}),du(\frac{\partial}{\partial x_{i}})\big )\bigg )^{\frac{1}{6}}\nonumber \\
&\quad\times\bigg (\sum_{i,j,l}h\big (du(\frac{\partial}{\partial x_{i}}),du(\frac{\partial}{\partial x_{j}})\big )h\big (du(\frac{\partial}{\partial x_{j}}),du(\frac{\partial}{\partial x_{l}})\big )h\big (du(\frac{\partial}{\partial x_{l}}),du(\frac{\partial}{\partial x_{i}})\big )\bigg )^{\frac{1}{3}}\nonumber\\
&\quad\times\bigg (\sum_{i,j,l}h\big (du(\frac{\partial}{\partial x_{i}}),du(\frac{\partial}{\partial x_{j}})\big )h\big (du(\frac{\partial}{\partial x_{j}}),du(\frac{\partial}{\partial x_{l}})\big )h\big (du(\frac{\partial}{\partial x_{l}}),du(\frac{\partial}{\partial x_{i}})\big )\bigg )^{\frac{1}{3}}\nonumber \\
&=m^{\frac{13}{12}}\bigg (\sum_{\alpha}\big (u^{\alpha}\big )^{2}\bigg )^{\frac{1}{2}}\Vert d_{(3)} u\Vert^{\frac{5}{3}}=m^{\frac{13}{12}}\Vert d_{(3)} u\Vert^{\frac{5}{3}}\bigg(\sum_{\alpha,\beta}h_{\alpha\beta}u^{\alpha}u^{\beta}\bigg )^{\frac{1}{2}},\nonumber
\end{aligned}
\end{equation*}
	where
	\begin{equation*}
	g^{ij}\frac{\partial r}{\partial x_{i}}\frac{\partial r}{\partial x_{j}}=|\nabla r|^{2}=1.
	\end{equation*}
	
	Therefore, by H\"older's inequality, we have
	
	\begin{equation*}
	\begin{aligned}
	& 2\int_{\partial B(R)}g^{is}g^{jk}g^{lr}h_{\alpha\beta}\frac{\partial u^{\alpha}}{\partial x_{i}}\frac{\partial r(x)}{\partial x_{j}}u^{\beta}
	h_{\gamma\delta}\frac{\partial u^{\gamma}}{\partial x_{k}}\frac{\partial u^{\delta}}{\partial x_{l}}h_{\mu\nu}\frac{\partial u^{\mu}}{\partial x_{r}}\frac{\partial u^{\nu}}{\partial x_{s}}ds_{g}\\
	\leq  & \int_{\partial B(R)}2m^{\frac{13}{12}}\Vert d_{(3)} u\Vert^{\frac{5}{3}}\bigg(\sum_{\alpha,\beta}h_{\alpha\beta}u^{\alpha}u^{\beta}\bigg )^{\frac{1}{2}}ds_{g}
	\end{aligned}
	\end{equation*}
	\begin{equation}
	\begin{aligned}
	\leq  & 2m^{\frac{13}{12}}\bigg (\int_{\partial B(R)}\Vert d_{(3)} u\Vert^{2}ds_{g}\bigg )^{\frac{5}{6}}\bigg (\int_{\partial B(R)}\bigg(\sum_{\alpha,\beta}h_{\alpha\beta}u^{\alpha}u^{\beta}\bigg )^{3}ds_{g}\bigg )^{\frac{1}{6}}.
	\end{aligned}
	\label{eq:L12}
	\end{equation}
	By ($\ref{eq:L3}$), we get\\
	\begin{equation}\label{eq:L13}
	\begin{aligned}
	&\int_{B(R)\backslash B(R_{2})}g^{is}g^{jk}g^{lr}\bigg (\frac{\partial h_{\alpha\beta}}{\partial y^{\xi}}u^{\xi}+2h_{\alpha\beta}\bigg )\frac{\partial u^{\alpha}}{\partial x_{i}}\frac{\partial u^{\beta}}{\partial x_{j}}
	 h_{\gamma\delta}\frac{\partial u^{\gamma}}{\partial x_{k}}\frac{\partial u^{\delta}}{\partial x_{l}}h_{\mu\nu}\frac{\partial u^{\mu}}{\partial x_{r}}\frac{\partial u^{\nu}}{\partial x_{s}}dv_{g}\\
	 \geq & \int_{B(R)\backslash B(R_{2})}g^{is}g^{jk}g^{lr}h_{\alpha\beta}\frac{\partial u^{\alpha}}{\partial x_{i}}\frac{\partial u^{\beta}}{\partial x_{j}}
	 h_{\gamma\delta}\frac{\partial u^{\gamma}}{\partial x_{k}}\frac{\partial u^{\delta}}{\partial x_{l}}h_{\mu\nu}\frac{\partial u^{\mu}}{\partial x_{r}}\frac{\partial u^{\nu}}{\partial x_{s}}dv_{g}\\
	 = & \int_{B(R)\backslash B(R_{2})}\Vert d_{(3)} u\Vert^{2}dv_{g}.
	 \end{aligned}
	 \end{equation}
	
	 Denote
	\begin{equation}
	Z(R)=\int_{B(R)\backslash B(R_{2})}\Vert d_{(3)} u\Vert^{2}dv_{g}+D(R_{1}).\label{eq:L14}
	\end{equation}
	Then
	\begin{equation}
	Z'(R)=\int_{\partial B(R)}\Vert d_{(3)} u\Vert^{2}ds_{g}.\label{eq:L15}
	\end{equation}
	According to ($\ref{eq:L9}$)-($\ref{eq:L15}$), we get
	\begin{equation}
	\begin{aligned}
	Z(R)&\leq  C_{1}\bigg (Z'(R)\bigg )^{\frac{5}{6}}\bigg (\int_{\partial B(R)}\bigg(\sum_{\alpha,\beta}h_{\alpha\beta}u^{\alpha}u^{\beta}\bigg )^{3}ds_{g}\bigg )^{\frac{1}{6}},\\
	& \le C_{1}\bigg (Z'(R)\bigg )^{\frac{5}{6}} \eta^{\frac 16}(R)\cdot \operatorname{Vol}_g^{\frac 16} (\partial B(R))\, ,
	\end{aligned}\label{eq:L96}
	\end{equation}
	where $C_{1}=2m^{\frac{13}{12}}$ is a positive constant, and
	\begin{equation}
	\eta(R)=\max_{r(x)=R}\left\{ \bigg(\sum_{\alpha,\beta}h_{\alpha\beta}u^{\alpha}u^{\beta}\bigg )^{3}\right\} .\label{eq:L97}
	\end{equation}
	Since $u(x)\rightarrow0$ as $r(x)\rightarrow\infty,$
	we note $\eta(R)\rightarrow0$
	as $R\rightarrow\infty$, and $\eta (R)$ is nonincreasing for sufficiently large $R\, .$
	Furthermore,
	\begin{equation}
	Z(R)-D(R_{1})=\int_{B(R)\backslash B(R_{2})}\Vert d_{(3)} u\Vert^{2}dv_{g}.\label{eq:L10}
	\end{equation}
	By Proposition \ref{prop:L2}, we have $R_{3}\geq R_{2}$ such
	that $Z(R)>0$ for any $R>R_{3}.$

It follows from (\ref{eq:L96}) that for any $R_{4}\geq R\geq R_{3},$ we have
	\begin{equation*}
	Z^{\frac{6}{5}}(R)\leq  C_{2}Z'(R)\big (\eta(R) \operatorname{Vol}_g(\partial B(R)\big )^{\frac 15},
	\end{equation*}
	and hence
	\begin{equation*}
	\begin{aligned}
	\int_{R}^{R_{4}}\frac{Z'(r)}{Z^{\frac{6}{5}}(r)}dr&\geq\frac{1}{C_{2}}\int_{R}^{R_{4}}\big (\eta(r) \operatorname{Vol}_g(\partial B(r)\big )^{-\frac 15}dr\\
	& \ge \frac{1}{C_{2}}\int_{R}^{R_{4}} \operatorname{Vol}_g^{-\frac 15}(\partial B(r))dr \cdot \eta^{-\frac 15}(R)\\
	&\ge  \frac{1}{C_{2}} R^{-\frac {\zeta}{5}} \eta^{-\frac 15}(R),
	\end{aligned}
    \end{equation*}
	where $C_{2}=C_{1}^{\frac{6}{5}}.$ Letting $R_{4}\rightarrow\infty$, we get
	\begin{equation*}
	\frac{1}{Z^{\frac{1}{5}}(R)}\geq\frac{1}{5C_{2}}R^{-\frac {\zeta}{5}} \eta^{-\frac 15}(R),
	\end{equation*}
	which implies that
	\begin{equation}
	Z(R)\leq  C_{3}R^\zeta \eta(R),\label{eq:L16}
	\end{equation}
	for constant $C_{3}=(5C_{2})^{5}.$ 	
	Hence, according to ($\ref{eq:L10}$) and ($\ref{eq:L16}$), we get
	\begin{equation*}
	\begin{aligned}\int_{B(R)}\frac{\Vert d_{(3)} u\Vert^{2}}{6} dv_{g} &= E_{\Phi_{(3)}}^{R}(u)  =\frac{Z(R)}{6}+\int_{B(R_{2})}\frac{\Vert d_{(3)} u\Vert^{2}}{6} dv_{g}-\frac{D(R_{1})}{6}\\
	&= o(R^{\zeta})\quad \operatorname{as}\quad R \to \infty.
	\end{aligned}
	\end{equation*}
	
\end{proof}

By Proposition $\ref{prop:L3}$ and monotonicy formula \eqref{eq:F28},
we have
\begin{thm}[Liouville Theorem] \label{thm:L4}
	Suppose $u:(M^{m},g)\rightarrow(N^{n},h)$ is a smooth $\Phi_{(3)}$-harmonic
	map and $r(x)$ satisfies the condition (\ref{eq:F7}). If $u(x)\rightarrow p_{0}\in N^{n}$
	as $r(x)\rightarrow\infty$ and
	\begin{equation*}
	\int_{R}^{\infty}\frac{1}{\operatorname{Vol}_{g}^{\frac{1}{5}}\big (\partial B(r)\big )}dr\geq CR^{-\frac{\zeta}{5}}
	\end{equation*}
	for $R$ large enough, then $u$ is a constant map.
\end{thm}
\begin{lem}[\cite{han52,han12,han20,han50,han24,han82,han51,wwdu}]
	\label{lem:LW1} Suppose $(M^{m},g)$ is a complete Riemannian
	manifold with a pole $x_{0}.$ We denote the radial curvature of $M^{m}$
	by $K_{r}$.\\
	$(1)$ If $-\alpha^{2}\leq  K_{r}\leq -\beta^{2}$ with $\alpha\geq\beta\geq0$,
	then
	\begin{equation*}
	\beta\coth(\beta r)[g-dr\otimes dr]\leq  \operatorname{Hess}(r)\leq \alpha\coth(\alpha r)[g-dr\otimes dr].
	\end{equation*}
	$(2)$ If $-\frac{A}{(1+r^{2})^{1+\varepsilon}}\leq  K_{r}\leq \frac{B}{(1+r^{2})^{1+\varepsilon}}$
	with $\varepsilon>0,$ $A\geq0$ and $0\leq  B<2\varepsilon,$
	then	
	\begin{equation*}
	\frac{1-B/2\varepsilon}{r}[g-dr\otimes dr]\leq  \operatorname{Hess}(r)\leq \frac{e^{A/2\varepsilon}}{r}[g-dr\otimes dr].
	\end{equation*}
	$(3)$ If $-\frac{a^{2}}{1+r^{2}}\leq  K_{r}\leq \frac{b^{2}}{1+r^{2}}$
	with $a\geq0$ and $b^{2}\in[0,\frac{1}{4}],$ then	
	\begin{equation*}
	\frac{1+\sqrt{1-4b^{2}}}{2r}[g-dr\otimes dr]\leq  \operatorname{Hess}(r)\leq \frac{1+\sqrt{1+4a^{2}}}{2r}[g-dr\otimes dr].
    \end{equation*}
    \end{lem}
By Lemma $\ref{lem:LW1}$, we have
\begin{lem}
	\label{lem:LW2}Suppose $(M^{m},g)$ is a complete Riemannian manifold
	with a pole $x_{0}.$ We denote the radial curvature of $M^{m}$ by
	$K_{r}$.\\
	$(i)$ If $-\alpha^{2}\leq  K_{r}\leq -\beta^{2}$ with $\alpha\geq\beta\geq0$
	and $(m-1)\beta-6\alpha>0,$ then
	\begin{equation*}
	1+\frac{m-1}{2}\lambda_{\min}-3\max\{2,\lambda_{max}\}\geq m-\frac{6\alpha}{\beta}.
	\end{equation*}
	$(ii)$ If $-\frac{A}{(1+r^{2})^{1+\varepsilon}}\leq  K_{r}\leq \frac{B}{(1+r^{2})^{1+\varepsilon}}$
	with $\varepsilon>0,$ $A\geq0$ and $0\leq  B<2\varepsilon,$
	then
	\begin{equation*}
	1+\frac{m-1}{2}\lambda_{\min}-3\max\{2,\lambda_{max}\}\geq1+(m-1)(1-\frac{B}{2\varepsilon})-6e^{\frac{A}{2\varepsilon}}.
	\end{equation*}
	$(iii)$ If $-\frac{a^{2}}{1+r^{2}}\leq  K_{r}\leq \frac{b^{2}}{1+r^{2}}$
	with $a\geq0$ and $0\leq  b^{2}\leq \frac{1}{4}$$,$
	then
	\begin{equation*}
	1+\frac{m-1}{2}\lambda_{\min}-3\max\{2,\lambda_{max}\}\geq1+(m-1)\frac{1+\sqrt{1-4b^{2}}}{2}-6\frac{1+\sqrt{1+4a^{2}}}{2}.
	\end{equation*}
\end{lem}
\begin{proof}
	If $K_{r}$ satisfies $(i)$, then by Lemma $\ref{lem:LW1}$, we have that
	on $B(r)/\{x_{0}\}$, for every $r>0$,
	\begin{equation*}
	\begin{aligned} & 1+\frac{m-1}{2}\lambda_{\min}-3\max\{2,\lambda_{max}\}\\
	\geq & 1+(m-1)\beta r\coth(\beta r)-6\alpha r\coth(\alpha r)\\
	= & 1+\beta r\coth(\beta r)\bigg( m-1-\frac{6\alpha}{\beta}\frac{\coth(\alpha r)}{\coth(\beta r)}\bigg ) \\
	\geq & m-\frac{6\alpha}{\beta}.
	\end{aligned}
    \end{equation*}
	The last inequality holds due to the increasing function $\beta r\coth(\beta r)\rightarrow1$
	as $r\rightarrow0$ and $\frac{\coth(\alpha r)}{\coth(\beta r)}<1$
	as $0<\beta<\alpha$. Using the same method as $(i)$, the above
	inequality is true for the case $(ii)$ and $(iii)$ on $B(r).$
\end{proof}
We can prove immediately the following theorem by Theorem \ref{thm:L1} and
Lemma $\ref{lem:LW2}$.
\begin{thm} \label{lthm1}
	Suppose $(M^{m},g)$ is a complete Riemannian manifold with a pole
	$x_{0}$ such that the radial curvature $K_{r}$ of $M$ satisfies
	one of the following conditions:\\
	$(i)$ $-\alpha^{2}\leq  K_{r}\leq -\beta^{2}$ with $\alpha\geq\beta\geq0$
	and $(m-1)\beta-6\alpha>0$,\\
	$(ii)$ $-\frac{A}{(1+r^{2})^{1+\varepsilon}}\leq  K_{r}\leq \frac{B}{(1+r^{2})^{1+\varepsilon}}$
	with $\varepsilon>0,$ $A\geq0$, $0\leq  B<2\varepsilon$
	and $1+(m-1)(1-\frac{B}{2\varepsilon})-6e^{\frac{A}{2\varepsilon}}>0,$\\
	$(iii)$ $-\frac{a^{2}}{1+r^{2}}\leq  K_{r}\leq \frac{b^{2}}{1+r^{2}}$
	with $a\geq0$, $0\leq  b^{2}\leq \frac{1}{4}$ and $1+(m-1)\frac{1+\sqrt{1-4b^{2}}}{2}-6\frac{1+\sqrt{1+4a^{2}}}{2}>0.$\\
	Suppose $u:(M^{m},g)\rightarrow(N^{n},h)$ is a smooth $\Phi_{(3)}$-harmonic
	map. If $u(x)\rightarrow p_{0}\in N^{n}$ as $r(x)\rightarrow\infty$
	and
	\begin{equation*}
	\int_{R}^{\infty}\frac{1}{\operatorname{Vol}_{g}^{\frac{1}{5}}\big (\partial B(r)\big )}dr\geq CR^{-\frac{\varLambda}{5}}
	\end{equation*}
	as $R$ large enough, then $u$ is a constant map, where
\begin{equation*}
\varLambda=\begin{cases}
m-\frac{6\alpha}{\beta} & \text{if }K_{r} \text{ satisfies (i)};\\
1+(m-1)(1-\frac{B}{2\varepsilon})-6e^{\frac{A}{2\varepsilon}} & \text{if }K_{r} \text{ satisfies (ii);}\\
1+(m-1)\frac{1+\sqrt{1-4b^{2}}}{2}-6\frac{1+\sqrt{1+4a^{2}}}{2} & \text{if }K_{r} \text{ satisfies (iii).}
\end{cases}
\end{equation*}
\end{thm}

\section{The second variation formula } \label{sec5}

For any smooth map $\Psi:(-\varepsilon,\varepsilon)\times(-\varepsilon,\varepsilon)\times M\rightarrow N$,
we denote $u_{s,t}(x)$ by $\Psi(s,t,x)$, where $\Psi(0,0,x)=u(x).$
Let
\begin{equation*}
V=d\Psi(\frac{\partial}{\partial t})_{\big {|}_{s,t=0}},\quad W=d\Psi(\frac{\partial}{\partial s})_{\big {|}_{s,t=0}}
\end{equation*}
be the variation vector fields of the deformation $u_{s,t}.$
\begin{thm}[The second variation formula] \label{thm:v1} Suppose $u: M^{m}\rightarrow N$ is a $\Phi_{(3)}$-harmonic
	map for the functional $E_{\Phi_{(3)}}$ and $u_{s,t}: M^{m}\rightarrow N, (-\varepsilon<s,t<\varepsilon)$
	is a compactly supported two-parameter variation. Then we have
	\begin{equation*}
	\begin{aligned}I(V,W) & =\frac{\partial^{2}}{\partial s\partial t}E_{\Phi_{(3)}}(u_{s,t})_{\big {|}_{s,t=0}}\\
	& =\int_{M}h\big (R^{N}\big (V,du(e_{i})\big )W, d_{(3)} u(e_{i})\big )dv_{g}\\
	& \quad+\int_{M}\sum_{i,j,k=1}^{m}h\big (\widetilde{\nabla}_{e_{i}}V,\widetilde{\nabla}_{e_{k}}W\big )h\big (du(e_{k}),du(e_{j})\big )h\big (du(e_{i}),du(e_{j})\big )dv_{g}\\
		& \quad+\int_{M}\sum_{i,j,k=1}^{m}h\big (\widetilde{\nabla}_{e_{i}}V,du(e_{k})\big )h\big (\widetilde{\nabla}_{e_{k}}W,du(e_{j})\big )h\big (du(e_{i}),du(e_{j})\big )dv_{g}\\
	& \quad+\int_{M}\sum_{i,j,k=1}^{m}h\big (\widetilde{\nabla}_{e_{i}}V,du(e_{k})\big )h\big (du(e_{k}),\widetilde{\nabla}_{e_{j}}W\big )h\big (du(e_{i}),du(e_{j})\big )dv_{g}\\
	& \quad+\int_{M}\sum_{i,j,k=1}^{m}h\big (\widetilde{\nabla}_{e_{i}}V,du(e_{k})\big )h\big (du(e_{k}),du(e_{j})\big )h\big (\widetilde{\nabla}_{e_{i}}W,du(e_{j})\big )dv_{g}\\
	& \quad+\int_{M}\sum_{i,j,k=1}^{m}h\big (\widetilde{\nabla}_{e_{i}}V,du(e_{k})\big )h\big (du(e_{k}),du(e_{j})\big )h\big (du(e_{i}),\widetilde{\nabla}_{e_{j}}W\big )dv_{g},
	\end{aligned}
	\end{equation*}
	where $R^{N}$ denotes the curvature tensor of $N.$
\end{thm}

	A $\Phi_{(3)}$-harmonic map is called stable if $I(V,V)\geq0$
	for any compactly supported vector field $V$ along $u.$
\begin{proof}
	We still use the symbols $\nabla$ and $\widetilde{\nabla}$ to denote
	the Levi-Civita connection on $(-\varepsilon,\varepsilon)\times(-\varepsilon,\varepsilon)\times M$
	and the induced connection on $\Psi^{-1}TN$ respectively. And we use $\{e_i\}$ to denote a locally orthonormal frame on $M$ and fix any point $x_{0}\in M$ such that
	 $\nabla_{e_{i}} {e_{j}} _{\big {|}_{x_0}} = 0$ for any $i,j$.
	
	From Theorem \ref{Ft1} and Proposition \ref{prop:F1}, at $x_{0},$ we have
	\begin{equation*}
	\begin{aligned}
	&\frac{\partial^{2}}{\partial s\partial t}E_{\Phi_{(3)}}(u_{s,t})_{\big {|}_{s,t=0}}\\
	= & -\frac{\partial}{\partial s}\int_{M}h\big (d\Psi(\frac{\partial}{\partial t}),\tau_{\Phi_{(3)}} {u_{s,t}}\big )_{\big {|}_{s,t=0}}dv_{g}\\
	= & -\sum_{i=1}^{m}\frac{\partial}{\partial s}\int_{M}h\big (d\Psi(\frac{\partial}{\partial t}),\widetilde{\nabla}_{e_{i}}d_{(3)} {u_{s,t}}(e_{i})\big )_{\big {|}_{s,t=0}}dv_{g}
\end{aligned}
\end{equation*}
\begin{equation} \label{eq:FS1}
\begin{aligned}
	= & -\sum_{i=1}^{m}\int_{M}h\big (d\Psi(\frac{\partial}{\partial t}),\widetilde{\nabla}_{\frac{\partial}{\partial s}}\widetilde{\nabla}_{e_{i}}d_{(3)} {u_{s,t}}(e_{i})\big )_{\big {|}_{s,t=0}}dv_{g}\\
	= & -\sum_{i=1}^{m}\int_{M}h\big (d\Psi(\frac{\partial}{\partial t}),\widetilde{\nabla}_{e_{i}}\widetilde{\nabla}_{\frac{\partial}{\partial s}}d_{(3)} {u_{s,t}}(e_{i})\big )_{\big {|}_{s,t=0}}dv_{g}\\
	& +\sum_{i=1}^{m}\int_{M}h\big (d\Psi(\frac{\partial}{\partial t}),R^{N}\big (d\Psi(\frac{\partial}{\partial s}),d\Psi(e_{i})\big )d_{(3)} {u_{s,t}}(e_{i})\big )_{\big {|}_{s,t=0}}dv_{g}\\	
	= & -\sum_{i=1}^{m}\int_{M}e_{i}h\big (d\Psi(\frac{\partial}{\partial t}),\widetilde{\nabla}_{\frac{\partial}{\partial s}}d_{(3)} {u_{s,t}}(e_{i})\big )_{\big {|}_{s,t=0}}dv_{g}\\
	& +\sum_{i=1}^{m}\int_{M}h\big (\widetilde{\nabla}_{e_{i}}d\Psi(\frac{\partial}{\partial t}),\widetilde{\nabla}_{\frac{\partial}{\partial s}}d_{(3)} {u_{s,t}}(e_{i})\big )_{\big {|}_{s,t=0}}dv_{g}\\
	& +\sum_{i=1}^{m}\int_{M}h\big (d\Psi(\frac{\partial}{\partial t}),R^{N}\big (d\Psi(\frac{\partial}{\partial s}),d\Psi(e_{i})\big )d_{(3)} {u_{s,t}}(e_{i})\big )_{\big {|}_{s,t=0}}dv_{g}.
	\end{aligned}
	\end{equation}
	
	We compute the second term in the right hand side of (\ref{eq:FS1}):
	\begin{equation}\label{eq:FS2}
	\begin{aligned}
	& \sum_{i=1}^{m}\int_{M}h\big (\widetilde{\nabla}_{e_{i}}d\Psi(\frac{\partial}{\partial t}),\widetilde{\nabla}_{\frac{\partial}{\partial s}}d_{(3)} {u_{s,t}}(e_{i})\big )_{\big {|}_{s,t=0}}dv_{g}\\
	= & \int_{M}\sum_{i,j,k=1}^{m}h\big (\widetilde{\nabla}_{e_{i}}V,\widetilde{\nabla}_{e_{k}}W\big )h\big (du(e_{i}),du(e_{j})\big )h\big (du(e_{j}),du(e_{k})\big )dv_{g}\\
	& +\int_{M}\sum_{i,j,k=1}^{m}h\big (\widetilde{\nabla}_{e_{i}}V,du(e_{k})\big )h\big (\widetilde{\nabla}_{e_{i}}W,du(e_{j})\big )h\big (du(e_{j}),du(e_{k})\big )dv_{g}\\
	& +\int_{M}\sum_{i,j,k=1}^{m}h\big (\widetilde{\nabla}_{e_{i}}V,du(e_{k})\big )h\big (du(e_{i}),\widetilde{\nabla}_{e_{j}}W\big )h\big (du(e_{j}),du(e_{k})\big )dv_{g}\\
%
	& +\int_{M}\sum_{i,j,k=1}^{m}h\big (\widetilde{\nabla}_{e_{i}}V,du(e_{k})\big )h\big (du(e_{i}),du(e_{j})\big )h\big (\widetilde{\nabla}_{e_{j}}W,du(e_{k})\big )dv_{g}\\
	& +\int_{M}\sum_{i,j,k=1}^{m}h\big (\widetilde{\nabla}_{e_{i}}V,du(e_{k})\big )h\big (du(e_{i}),du(e_{j})\big )h\big (du(e_{j}),\widetilde{\nabla}_{e_{k}}W\big )dv_{g}.
	\end{aligned}
	\end{equation}
	
	The integrand for the first term in the right hand side of (\ref{eq:FS1}) is
	\begin{equation*}
	\begin{aligned}
		& e_{i}h \big (d\Psi(\frac{\partial}{\partial t}),\widetilde{\nabla}_{\frac{\partial}{\partial s}}d_{(3)} {u_{s,t}}(e_{i}) \big )\\
	= & e_{i}h \bigg (d\Psi(\frac{\partial}{\partial t}),\widetilde{\nabla}_{\frac{\partial}{\partial s}}\bigg ( h\big (d\Psi(e_{i}),d\Psi(e_{j})\big )h\big (d\Psi(e_{j}),d\Psi(e_{k})\big )d\Psi(e_{k})\bigg ) \bigg )\\
	= & e_{i}\bigg (h\big(d\Psi(\frac{\partial}{\partial t}),\widetilde{\nabla}_{e_{k}}d\Psi(\frac{\partial}{\partial s})\big)h\big (d\Psi(e_{i}),d\Psi(e_{j})\big )h\big (d\Psi(e_{j}),d\Psi(e_{k})\big )\bigg )\\
	& +e_{i}\bigg ( h\big(d\Psi(\frac{\partial}{\partial t}),d\Psi(e_{k})\big)h\big(\widetilde{\nabla}_{e_{i}}d\Psi(\frac{\partial}{\partial s}),d\Psi(e_{j})\big)h\big (d\Psi(e_{j}),d\Psi(e_{k})\big )\bigg )
\end{aligned}
	\end{equation*}
\begin{equation}
	\begin{aligned}
	& +e_{i}\bigg ( h\big(d\Psi(\frac{\partial}{\partial t}),d\Psi(e_{k})\big)h\big(d\Psi(e_{i}),\widetilde{\nabla}_{e_{j}}d\Psi(\frac{\partial}{\partial s})\big)h\big (d\Psi(e_{j}),d\Psi(e_{k})\big )\bigg )\\
	& +e_{i}\bigg ( h \big(d\Psi(\frac{\partial}{\partial t}),d\Psi(e_{k})\big )h\big (d\Psi(e_{i}),d\Psi(e_{j})\big )h\big (\widetilde{\nabla}_{e_{j}}d\Psi(\frac{\partial}{\partial s}),d\Psi(e_{k})\big )\bigg )\\
	& +e_{i}\bigg ( h\big(d\Psi(\frac{\partial}{\partial t}),d\Psi(e_{k})\big)h\big (d\Psi(e_{i}),d\Psi(e_{j})\big )h\big(d\Psi(e_{j}),\widetilde{\nabla}_{e_{k}}d\Psi(\frac{\partial}{\partial s})\big)\bigg ).
	\end{aligned}
	\label{eq:FS3}
	\end{equation}
	
	Let $X_{1}$, $X_{2},$ $X_{3},$ $X_{4}$ and $X_{5}$ be compactly
	supported vector fields on $M$ and $Y$ be any vector field on $M$.
	We have
	\begin{equation*}
	\begin{aligned}
	& g(X_{1},Y)=h\big(V,\widetilde{\nabla}_{e_{k}}W\big)h\big (du(Y),du(e_{j})\big )h\big (du(e_{j}),du(e_{k})\big ),\\
	& g(X_{2},Y)=h\big (V,du(e_{k})\big )h\big(\widetilde{\nabla}_{Y}W,du(e_{j})\big)h\big (du(e_{j}),du(e_{k})\big ),\\
	& g(X_{3},Y)=h\big (V,du(e_{k})\big )h\big(du(Y),\widetilde{\nabla}_{e_{j}}W\big)h\big (du(e_{j}),du(e_{k})\big ),\\
	& g(X_{4},Y)=h\big (V,du(e_{k})\big )h\big (du(Y),du(e_{j})\big )h\big(\widetilde{\nabla}_{e_{j}}W,du(e_{k})\big),\\
	& g(X_{5},Y)=h\big (V,du(e_{k})\big )h\big (du(Y),du(e_{j})\big )h\big(du(e_{j}),\widetilde{\nabla}_{e_{k}}W\big).
	\end{aligned}
	\end{equation*}
	
	Hence, when $s=0$ and $t=0$, (\ref{eq:FS3}) becomes
	\begin{equation}
	\begin{aligned} & \sum_{i=1}^{m}\bigg ( e_{i}g(X_{1},e_{i})+e_{i}g(X_{2},e_{i})+e_{i}g(X_{3},e_{i})+e_{i}g(X_{4},e_{i})+e_{i}g(X_{5},e_{i})\bigg )\\
	= & \operatorname{div}(X_{1})+\operatorname{div}(X_{2})+\operatorname{div}(X_{3})+\operatorname{div}(X_{4})+\operatorname{div}(X_{5}).
	\end{aligned}
	\label{eq:FS4}
	\end{equation}
	
	The result follows from (\ref{eq:FS1})-(\ref{eq:FS4}).
\end{proof}
\section{Examples of $\Phi_{(3)}$-SSU manifolds} \label{sec6}

Proceeding as in \cite{ww}, we obtain many examples of $\Phi_{(3)}$-SSU manifolds.

\begin{thm}
	\label{thm:SS1}A hypersurface $M$ in Euclidean space is $\Phi_{(3)}$-$\operatorname{SSU}$ if and only if its principal curvatures satisfy
	\begin{equation*}
	0<\lambda_{1}\leq \lambda_{2}\leq \cdots\leq \lambda_{m}<\frac{1}{5}(\lambda_{1}+\cdots+\lambda_{m-1}).
	\end{equation*}
\end{thm}
\begin{proof}
	Similar to the proof of Theorem $5.1$ in \cite{han}, from the
	definition of the $\Phi_{(3)}$-SSU, we have
	\begin{equation*}
	\begin{aligned} & \sum_{i=1}^{m}\big (6\langle B(v,e_i),B(v,e_i)\rangle-\langle B(v,v),B(e_i,e_i)\rangle\big )\\
	\leq  & \lambda_i \big  (6\lambda_{m}-\sum_{i=1}^{m}\lambda_{i}\big  )=\lambda_{i}\big (5\lambda_{m}-\sum_{i=1}^{m-1}\lambda_{i}\big )<0,
	\end{aligned}
	\end{equation*}
	that is, $\lambda_{m}<\frac{1}{5}(\lambda_{1}+\cdots+\lambda_{m-1}).$
	We finish the proof.
\end{proof}
Then we have
\begin{cor}
	The standard sphere $S^{m}$ is $\Phi_{(3)}$-$\operatorname{SSU}$ if and only if
	$m>6.$
\end{cor}
\begin{proof}
	As $S^{m}$ is a compact convex hypersurface in $\mathbb{R}^{m+1},$
	according to Theorem $\ref{thm:SS1}$, its principle curvatures satisfy
	\begin{equation*}
	\lambda_{1}=\lambda_{2}=\cdots=\lambda_{m}=1.
	\end{equation*}
	
	Hence, $m>6.$ We finish the proof.
\end{proof}
\begin{cor}
	The graph of $f(x)=x_{1}^{2}+\cdots+x_{m}^{2},$ $x=(x_{1},\cdots,x_{m})\in\mathbb{R}^{m}$
	is $\Phi_{(3)}$-$\operatorname{SSU}$ if and only if $m>6.$
\end{cor}
\begin{lem}[\cite{wy}] \label{thm:38}An Euclidean hypersurface is $p$-$\operatorname{SSU}$ if
	and only if its principal curvatures satisfy
	\begin{equation*}
	0<\lambda_{1}\leq \lambda_{2}\leq \cdots\leq \lambda_{m}<\frac{1}{p-1}(\lambda_{1}+\cdots+\lambda_{m-1}).
	\end{equation*}
\end{lem}

\begin{thm}	\label{thmp5}
	Every $\Phi_{(3)}$-$\operatorname{SSU}$ manifold $M$ is $p$-$\operatorname{SSU}$ for any $2 \le p \le 6$.
\end{thm}
\begin{proof}
	By Definition, $\Phi_{(3)}$-$\operatorname{SSU}$ manifold enjoys
	\begin{equation}
	\begin{aligned}
	{F_{\Phi_{(3)},x}} (v)=\sum_{i=1}^{m}\big (6\langle B(v,e_{i}),B(v,e_{i})\rangle _{\mathbb R^q}-\langle B(v,v),B(e_{i},e_{i})\rangle_{\mathbb R^q}\big )<0
	\end{aligned}
	\end{equation}
	for all unit tanget vector $v \in T_x(M)$. It follows that
	\begin{equation}
	\begin{aligned}
	F_{p,x}(v)&=(p-2)\langle \mathsf B(v,v), B(v,v)\rangle _{\mathbb R^q} + \langle Q^M_x(v),v\rangle _M\\
	& \le (p-2)\sum^m_{i=1} \bigg (2
	\langle B(v,e_i), B(v, e_i)\rangle _{\mathbb R^q} \bigg ) \\
	& \quad + \sum^m_{i=1} \bigg (2
	\langle B(v,e_i), B(v, e_i)\rangle _{\mathbb R^q} -
	\langle B(v,v), B(e_i, e_i)\rangle _{\mathbb R^q} \bigg )\\
	&\leq\sum_{i=1}^m\big (p\langle B(v, e_i), B(v, e_i)\rangle-\langle B(v,v), B( e_i, e_i)\rangle\big ) \\
	&\leq\sum_{i=1}^m\big (6\langle B(v,e_i), B(v, e_i)\rangle-\langle B(v,v),B( e_i, e_i)\rangle\big )<0,
	\end{aligned}
	\end{equation}
	for $2\leq p\leq 6$.
	So, $M$ is $p$-SSU for any $2\leq p\leq 6$.
\end{proof}





\begin{thm}[Topological Vanishing Theorems] \label{tmp6}
		
	Every compact $\Phi_{(3)}$-SSU manifold $M$ is $6$-connected, i.e.,
	
	\begin{equation}
	\pi_1(M) = \cdots = \pi_{6}(M) = 0.
	\end{equation}
	
\end{thm}

\begin{proof}
	
	Since every compact $p$-SSU manifold is $[p]$-connected (cf. \cite{wwre} Theorem 3.10), and $p=6$ by the previous Theorem, the result follows.
\end{proof}

\begin{thm} \label{tmp7}
		
	The dimension of any compact $\Phi_{(3)}$-SSU manifold $M$ is greater than $6$.
	
\end{thm}

\begin{proof}
	Suppose that $m\leq 6$, then $M$ is not a $6$-SSU manifold (cf. \cite{wy} Theorem 3.10). By the preceeding Theorem, $M$ is not a
	$\Phi_{(3)}$-$\text{SSU}$ manifold. Hence, the dimension of any compact $\Phi_{(3)}$-SSU manifold $M$ is greater than $6$.
\end{proof}

\begin{thm}[Sphere Theorems]
	
	Every compact $\Phi_{(3)}$-$\text{SSU}$ manifold $M$ of dimension $m \le 13$ is homeomorphic to an $m$-sphere.
\end{thm}

\begin{proof}
	In view of Theorem \ref{tmp6}, $
	M$ is 6-connected. By the Hurewicz isomorphism theorem, the 6-connectedness of $M$ implies homology groups $H_1(M)=\cdots=H_6(M)=0$. It follows from Proincare Duality Theorem and the Hurewicz Isomorphism Theorem (cf. E. Spanier \cite{span}) again, $H_{m-6}(M)=\cdots=H_{m-1}(M)=0$, $H_{m}(M)\neq 0$, $m \le 13$ and $M$ is $(m-1)$-connected. Hence, $M$ is a homotopy $m$-sphere, $m \le 13$. Since $M$ is $\Phi_{(3)}$-$\text{SSU}$ manifold, $m\geq 7$. Consequently, a homotopy $m$-sphere $M$ for $m \ge 7$ is homeomorphic to an $m$-sphere
	by a Theorem of S. Smale \cite{smale}.
\end{proof}

\begin{thm}
	\label{thm:SS2}
	Suppose that $\widetilde{M}$ is a compact convex hypersurface of $\mathbb{R}^{q}$ and the principal curvatures of $\widetilde{M}$ satisfy
	\begin{equation*}
	0<\lambda_{1}\leq \lambda_{2}\leq \cdots\leq \lambda_{q-1}.
	\end{equation*}
	If $\langle Ric^{M}(v),v\rangle>\frac{5}{6}k\lambda_{q-1}^{2}$,
	where $M$ is a compact connected minimal $k$-submanifold of $\widetilde{M}$
	and $v$ is any unit tangent vector to $M$, then $M$ is $\Phi_{(3)}$-$\operatorname{SSU}$.
\end{thm}
\begin{proof}
	We denote the second fundamental form of $M$ in $\mathbb{R}^{q}$,
	$M$ in $\widetilde{M}$ and $\widetilde{M}$ in $\mathbb{R}^{q}$
	by $B$, $B_{1}$ and $\widetilde{B}$. According to Gauss equation,
	we get
	\begin{equation}
	B(X,Y)=B_{1}(X,Y)+\widetilde{B}(X,Y)\vartheta,\label{eq:ESS0}
	\end{equation}
	where $\vartheta$ is the unit normal field of $\widetilde{M}$ in $\mathbb{R}^{q}$.
	By the definition of minimal submanifold, we have
	\begin{equation}
	 \sum_{i=1}^{k}B(e_{i},e_{i})=\sum_{i=1}^{k}B_{1}(e_{i},e_{i})+\sum_{i=1}^{k}\widetilde{B}(e_{i},e_{i})\vartheta=\sum_{\alpha=1}^{k}\widetilde{B}(e_{i},e_{i})\vartheta,\label{eq:ESS1}
	\end{equation}
	where $\{e_{i}\}_{i=1}^{k}$ is a local orthonormal
	frame on $M$. Denote $\widetilde{B}(e_{i},e_{j})=\lambda_{i}\delta_{i j}.$
	
	Hence,
	\begin{equation*}
	\begin{aligned}
	&\sum_{i=1}^{k}\big (6\langle B(v,e_i),B(v,e_i)\rangle-\langle B(v,v),B(e_i,e_i)\rangle\big )\\
		= & -6\langle Ric^{M}(v),v\rangle+5\sum_{i=1}^{k}\langle B(v,v),B(e_i,e_i)\rangle\\
%
	= & -6\langle Ric^{M}(v),v\rangle+5\sum_{i=1}^{k}\widetilde{B}(v,v)\widetilde{B}(e_i,e_i)\\
	\leq  & -6\langle Ric^{M}(v),v\rangle+5\sum_{i=1}^{k}\lambda_{i}\lambda_{q-1}\\
	\leq  & -6\langle Ric^{M}(v),v\rangle+5k\lambda_{q-1}^{2}<0.
	\end{aligned}
	\end{equation*}
	The first equality holds due to Gauss equation. We have the desired
	result.
\end{proof}
The following lemma will be used in our later proof. For the second
fundamental form of an Ellipsoid in $\mathbb{R}^{m+1}$, we have
\begin{lem}[\cite{hw,sp,wwll}]
    Let $\{\lambda_{i}\}_{i=1}^{m}$ be a family of principal curvatures of $E^{m}$ in $\mathbb{R}^{m+1}$ with $0<\lambda_{1}\leq \lambda_{2}\leq \cdots\leq \lambda_{m}$. Then
\begin{equation*}
\frac{\min(a_{i})}{\big (\max(a_{i})\big )^{2}}\leq  \lambda_{1}\leq \lambda_{2}\leq \cdots\leq \lambda_{m}\leq \frac{\max(a_{i})}{\big (\min(a_{i})\big )^{2}},
\end{equation*}
where
\begin{equation*}
E^{m}=\left\{ (x_{1},\cdots,x_{m+1})\in\mathbb{R}^{m+1}:\frac{x_{1}^{2}}{a_{1}^{2}}+\cdots+\frac{x_{m+1}^{2}}{a_{m+1}^{2}}=1, a_{i}>0, 1\leq i\leq m+1\right\}.
\end{equation*}	
\end{lem}
Similarly, we can prove the following results by using Theorem $\ref{thm:SS2}$.
As the proofs are similar, we omit the details.
\begin{thm}
	Suppose that $M$ is a compact minimal $k$-submanifold of an ellipsoid
	$E^{q-1}$ in $\mathbb{R}^{q}$ and $v$ is any unit tangent vector
	to $M.$ Then $M$ is $\Phi_{(3)}$-$\operatorname{SSU}$ when
	\begin{equation*}
	\langle Ric^{M}(v),v\rangle>\frac{5}{6}\frac{\big (\max_{1\leq  i\leq  q}(a_{i})\big )^{2}}{\big (\min_{1\leq  i\leq  q}(a_{i})\big )^{4}}k.
	\end{equation*}
\end{thm}
\begin{cor}
	Suppose that $M$ is a compact minimal $k$-submanifold of the unit
	sphere $S^{q-1}$ and $v$ is any unit tangent vector to $M.$ Then
	$M$ is $\Phi_{(3)}$-$\operatorname{SSU}$ when $\langle Ric^{M}(v),v\rangle>\frac{5}{6}k.$
\end{cor}
\begin{thm}
	Suppose that $M$ is a compact $k$-submanifold of the unit
	sphere $S^{q-1}$ and $B_{1}$ is the second fundamental form of $M$ in $S^{q-1}$.
	Then $M$ is $\Phi_{(3)}$-$\operatorname{SSU}$ when
	\begin{equation*}
	\Vert B_{1}\Vert^{2}<\frac{k-6}{\sqrt{k}+6},
	\end{equation*}
	where $k>6.$
\end{thm}
\begin{proof}
	By Theorem $\ref{thm:SS2}$ and Cauchy-Schwarz inequality, we get
	\begin{equation*}
	\begin{aligned} & \sum_{i=1}^{k}\big (6\langle B(v,e_{i}),B(v,e_{i})\rangle-\langle B(v,v),B(e_{i},e_{i})\rangle\big )\\
	= & \sum_{i=1}^{k}\big (6\langle B_{1}(v,e_{i}),B_{1}(v,e_{i})\rangle-\langle B_{1}(v,v),B_{1}(e_{i},e_{i})\rangle\big )-(k-6)\\
%
	\leq  & \sum_{i=1}^{k}6\langle B_{1}(v,e_{i}),B_{1}(v,e_{i})\rangle+|B_{1}(v,v)| \big (\big| \sum_{i=1}^{k}B_{1}(e_{i},e_{i})\big|^{2}\big ) ^{\frac{1}{2}}-(k-6)\\
	\leq  & 6\Vert B_{1}\Vert^{2}+\sqrt{k}|B_{1}(v,v)|\big( \sum_{i=1}^{k}|B_{1}(e_{i},e_{i})|^{2}\big ) ^{\frac{1}{2}}-(k-6)\\
	\leq  & (6+\sqrt{k})\Vert B_{1}\Vert^{2}-(k-6)<0.
	\end{aligned}
	\end{equation*}
	
	Hence, by the definition of the $\Phi_{(3)}$-SSU, we get the desired
	result.
\end{proof}

\section{Stable $\Phi_{(3)}$-harmonic maps from $\Phi_{(3)}$-SSU manifolds} \label{sec7}

We recall some definitions and facts of submanifolds which will be used in the following results, see {\cite{hw}}.

Let $M^{m}$ be isometrically immersed in the Euclidean space $\mathbb{R}^{q}$
and $B$ be the second fundamental form $M$ in $\mathbb{R}^{q}$. We denote
the standard flat connection of $\mathbb{R}^{q}$ and the Riemannian connection
on $M$ by $\overline{\nabla}$ and $\nabla$. These are related by
\begin{equation*}
\overline{\nabla}_{X}Y=\nabla_{X}Y+B(X,Y),
\end{equation*}
where $X,Y$ are smooth vector fields on $M.$ The tensors $A$ and
$B$ are related by
\begin{equation}
\langle A^{\eta}X,Y\rangle=\langle B(X,Y),\eta\rangle,\label{eq:ssu}
\end{equation}
where $A^{\eta}X$ is the Weingarten map with the normal vector field $\eta \in T^{\bot}M$.

For each $x\in M$, we denote an orthonormal basis for the normal
space $T_{x}^{\bot}M$ to $M$ at $x$ by $
\{e_{m+1},\cdots, e_{q}\}.$
Let $v \in T_{x}M$. The Ricci tensor $Ric^{M}:T_{x}M\rightarrow T_{x}M$
is defined by
\begin{equation*}
Ric^{M}(v)=\sum_{i=1}^{m}R(v,e_{i})e_{i}.
\end{equation*}

From the Gauss curvature equation, we have
\begin{equation}
Ric^{M}=\sum_{\alpha=m+1}^{q}\mathrm{trace}(A^{e_\alpha})A^{e_\alpha}-\sum_{\alpha=m+1}^{q}A^{e_\alpha}A^{e_\alpha}.\label{eq:SS1}
\end{equation}

Then we have
\begin{thm} \label{thm:fss1}
	Let $(M^{m},g)$ be a compact $\Phi_{(3)}$-$\operatorname{SSU}$ manifold and $(N,h)$
	be a compact Riemannian manifold. Then every stable $\Phi_{(3)}$-harmonic
	map $u:(M^{m},g)\rightarrow(N,h)$ is constant.
\end{thm}

\begin{proof}
	Let $\{v_{\ell}^\top\}$ be the tangential projection of an orthonormal frame field $\{v_{\ell}\}_{\ell=1}^q$ in $\mathbb R^q$ onto $M$.
	For convenience, we choose $\{v_1,\cdots,v_m\}=\{e_1,\cdots,e_m\}$ to be tangential to $M$, $\{v_{m+1},\cdots,v_q\}=\{e_{m+1},\cdots,e_{q}\}$ to be normal to $M$, and $\nabla ^{\Psi}e_i = 0$ at a point in $M$.
	%
	Since $v_{\ell}^\top=v_{\ell}-v_{\ell}^\bot$ and $v_{\ell}$ are parallel in $\mathbb R^q$, we have
	\begin{equation}
	\begin{aligned}
	\nabla ^u_{e_i} du(v_{\ell}^\top)&=(\nabla ^u_{e_i} du)(v_{\ell}^\top) + du (\nabla ^M_{e_i}v_{\ell}^\top)= (\nabla ^u_{e_i} du)(v_{\ell}^\top) + du \bigg (\left(\nabla ^{\mathbb R ^q}_{e_i}(v_{\ell}-v_{\ell}^\bot)\right)^\top\bigg )\\
	&= (\nabla ^u_{e_i} du)(v_{\ell}^\top) + du \left(A^{v_{\ell}^\bot}(e_i)\right )\, .
	\end{aligned}\label{3.18}
	\end{equation}
	
	%
	
	In view of (\ref{eq:ssu}), we have
	
	\begin{equation}
	\nabla ^u_{e_i} du(v_{\ell}^\top)=\sum_{k=1}^{m}\nabla^u_{e_{i}}du(e_{k}) + \sum_{\ell=m+1}^{q}\sum_{k=1}^{m}B_{ik}^{\ell}du(e_{k}),
	\end{equation}
	where  $B_{ij}^{\ell} = \langle B(e_i, e_j), e_{\ell}\rangle$ is a
	components of $B.$
	
	%
	%
	
	According to Proposition \ref{prop:F1}, we have
	\begin{eqnarray}
	&& \sum_{\ell=1}^q\int_{M}\langle(\Delta du)(v_{\ell}), d_{(3)} u(v_{\ell})\rangle dv_{g}=\int_{M}\sum_{i,j,\ell} \delta_{ij}\langle(\Delta du)(e_{i}), d_{(3)} u(e_{j})\rangle dv_{g}\nonumber\\
	&& =\int_{M}\sum_{i=1}^m\langle(\Delta du)(e_{i}), d_{(3)} u(e_{i})\rangle dv_{g}=\int_{M}\langle(\Delta du), d_{(3)} u\rangle dv_{g}=\int_{M}\langle\delta du,\delta d_{(3)} u\rangle dv_{g}\nonumber\\
	&& =-\int_{M}\langle\delta du,\tau_{\Phi_{(3)}} (u)\rangle dv_{g}=0.
	\label{eq:SS0}
	\end{eqnarray}
	
	By using the Weitzenb\"ock formula, we have
	\begin{equation*}
	-\sum_{k}R^{N}\big (du(X),du(e_{k})\big )du(e_{k})+du\big (Ric^{M}(X)\big )=(\Delta du)(X)+(\nabla^{2}du)(X),
	\end{equation*}
	where $X$ is a smooth vector field in $M$. We assume $i,j,k,\hbar,\imath \in\{1,\cdots,m\},$ $\alpha,\beta \in\{m+1,\cdots,q\},$ $\ell \in\{1, \cdots
	m, \cdots q\}.$
	Hence,
	\begin{equation} \label{eq:SS2}
	\begin{aligned}
	& \sum_{\ell=1}^qI\big (du(v_{\ell}^\top),du(v_{\ell}^\top)\big )\\
	= & -\int_{M}\sum_{i=1}^mh\big (du(Ric^{M}(e_{i})), d_{(3)} u(e_{i})\big )dv_{g}+\int_{M}\sum_{i=1}^mh\big ((\nabla^{2}du)(e_{i}), d_{(3)} u(e_{i})\big )dv_{g}\\
	& +  \int_{M}\sum_{i,j,k,\ell}h\big (\widetilde{\nabla}_{e_{i}}du(v_{\ell}^\top),\widetilde{\nabla}_{e_{k}}du(v_{\ell})\big )h\big (du(e_{k}),du(e_{j})\big )h\big (du(e_{i}),du(e_{j})\big )dv_{g}\\
	& +  \int_{M}\sum_{i,j,k,\ell}h\big (\widetilde{\nabla}_{e_{i}}du(v_{\ell}^\top),du(e_{k})\big )h\big (\widetilde{\nabla}_{e_{k}}du(v_{\ell}),du(e_{j})\big )h\big (du(e_{i}),du(e_{j})\big )dv_{g}\\
	& + \int_{M}\sum_{i,j,k,\ell}h\big (\widetilde{\nabla}_{e_{i}}du(v_{\ell}^\top),du(e_{k})\big )h\big (du(e_{k}),\widetilde{\nabla}_{e_{j}}du(v_{\ell})\big )h\big (du(e_{i}),du(e_{j})\big )dv_{g}\\
	& +  \int_{M}\sum_{i,j,k,\ell}h\big (\widetilde{\nabla}_{e_{i}}du(v_{\ell}^\top),du(e_{k})\big )h\big (du(e_{k}),du(e_{j})\big )h\big (\widetilde{\nabla}_{e_{i}}du(v_{\ell}^\top),du(e_{j})\big )dv_{g}\\
	& + \int_{M}\sum_{i,j,k,\ell}h\big (\widetilde{\nabla}_{e_{i}}du(v_{\ell}^\top),du(e_{k})\big )h\big (du(e_{k}),du(e_{j})\big )h\big (du(e_{i}),\widetilde{\nabla}_{e_{j}}du(v_{\ell}^\top)\big )dv_{g}.
	\end{aligned}
	\end{equation}

	We compute at $x_0$. Via (\ref{eq:SS1}), The first integrand in (\ref{eq:SS2}) is

	\begin{equation}
	\begin{aligned}
	& \sum_{i=1}^m h\big (du(Ric^{M}(e_{i})), d_{(3)} u(e_{i})\big )\\
	= & \sum_{i,\alpha}h\big (du(\mathrm{trace}(A^{e_\alpha})A^{e_\alpha}(e_{i})), d_{(3)} u(e_{i})\big )-\sum_{i,\alpha}h\big (du(A^{e_\alpha}A^{e_\alpha}(e_{i})), d_{(3)} u(e_{i})\big ).
	\end{aligned}
	\label{eq:SS3}
	\end{equation}
	
	The second integrand in (\ref{eq:SS2}) is
	\begin{equation*}
	\begin{aligned}
	& \sum_{i=1}^m h\big ((\nabla^{2}du)(e_{i}), d_{(3)} u(e_{i})\big )\\
	= & \sum_{i,j,k}h\big ((\nabla^{2}du)(e_{i}),du(e_{k})\big )h\big (du(e_{i}),du(e_{j})\big )h\big (du(e_{j}),du(e_{k})\big )\\
	= & \sum_{i,j,k,\hbar}e_{\hbar}\bigg (h\big (\nabla_{e_{\hbar}}du(e_{i}),du(e_{k})\big )h\big (du(e_{i}),du(e_{j})\big )h\big (du(e_{j}),du(e_{k})\big )\bigg )\\
	& - \sum_{i,j,k,\hbar}h\big (\nabla_{e_{\hbar}}du(e_{i}),\nabla_{e_{\hbar}}du(e_{k})\big )h\big (du(e_{i}),du(e_{j})\big )h\big (du(e_{j}),du(e_{k})\big )\\
	& - \sum_{i,j,k,\hbar}h\big (\nabla_{e_{\hbar}}du(e_{i}),du(e_{k})\big )h\big (\nabla_{e_{\hbar}}du(e_{i}),du(e_{j})\big )h\big (du(e_{j}),du(e_{k})\big )\\
	& - \sum_{i,j,k,\hbar}h\big (\nabla_{e_{\hbar}}du(e_{i}),du(e_{k})\big )h\big (du(e_{i}),\nabla_{e_{\hbar}}du(e_{j})\big )h\big (du(e_{j}),du(e_{k})\big )\\
	\end{aligned}
\end{equation*}
\begin{equation}
\begin{aligned}
	& - \sum_{i,j,k,\hbar}h\big (\nabla_{e_{\hbar}}du(e_{i}),du(e_{k})\big )h\big (du(e_{i}),du(e_{j})\big )h\big (\nabla_{e_{\hbar}}du(e_{j}),du(e_{k})\big )\\
	& - \sum_{i,j,k,\hbar}h\big (\nabla_{e_{\hbar}}du(e_{i}),du(e_{k})\big )h\big (du(e_{i}),du(e_{j})\big )h\big (du(e_{j}),\nabla_{e_{\hbar}}du(e_{k})\big ).
	\end{aligned}
	\label{eq:SS4}
	\end{equation}

	The third integrand in (\ref{eq:SS2}) is
	\begin{equation} \label{eq:SS5}
	\begin{aligned}
	&\sum_{i,j,k,\ell}h\big (\widetilde{\nabla}_{e_{i}}du(v_{\ell}^\top),\widetilde{\nabla}_{e_{k}}du(v_{\ell}^\top)\big)h\big (du(e_{k}),du(e_{j})\big )h\big (du(e_{i}),du(e_{j})\big )\\
	= & \sum_{i,j,k,\hbar,\imath,\alpha,\beta}h\bigg (B_{i\hbar}^{\alpha}du(e_{\hbar})+\widetilde{\nabla}_{e_{i}}du(e_{\hbar}),B_{k\imath}^{\beta}du(e_{\imath})\\
	&  +\widetilde{\nabla}_{e_{k}}du(e_{\imath})\bigg )h\big (du(e_{k}),du(e_{j})\big )h\big (du(e_{i}),du(e_{j})\big )\\
	= & \sum_{i,j,k,\hbar,\imath,\alpha}B_{i\hbar}^{\alpha}B_{k\imath}^{\alpha}h\big (du(e_{\hbar}),du(e_{\imath})\big )h\big (du(e_{k}),du(e_{j})\big )h\big (du(e_{i}),du(e_{j})\big )\\
	& + \sum_{i,j,k,\hbar}h\big((\nabla_{e_{\hbar}}du)(e_{i}),(\nabla_{e_{\hbar}}du)(e_{k})\big)h\big (du(e_{k}),du(e_{j})\big )h\big (du(e_{i}),du(e_{j})\big )\\
	= & \sum_{i,j,k.\alpha}h\big(du(A^{e_\alpha}(e_{i})),du(A^{e_\alpha}(e_{k}))\big)h\big (du(e_{k}),du(e_{j})\big )h\big (du(e_{i}),du(e_{j})\big )\\
	& + \sum_{i,j,k,\hbar}h\big((\nabla_{e_{\hbar}}du)(e_{i}),(\nabla_{e_{\hbar}}du)(e_{k})\big)h\big (du(e_{k}),du(e_{j})\big )h\big (du(e_{i}),du(e_{j})\big ).
	\end{aligned}
	\end{equation}
	
	The fourth integrand in (\ref{eq:SS2}) is
	\begin{equation*}
	\begin{aligned}
	&\sum_{i,j,k,\ell}h\big (\widetilde{\nabla}_{e_{i}}du(v_{\ell}^\top),du(e_{k})\big)h\big(\widetilde{\nabla}_{e_{k}}du(v_{\ell}^\top),du(e_{j})\big)h\big (du(e_{i}),du(e_{j})\big )\\
	= & \sum_{i,j,k,\hbar.\imath,\alpha,\beta}h\big(B_{i\hbar}^{\alpha}du(e_{\hbar})+\widetilde{\nabla}_{e_{i}}du(e_{\hbar}),du(e_{k})\big)\\
	& \cdot h\big(B_{k\imath}^{\beta}du(e_{\imath})+\widetilde{\nabla}_{e_{k}}du(e_{\imath}),du(e_{j})\big)h\big (du(e_{i}),du(e_{j})\big )
\end{aligned}
	\end{equation*}
\begin{equation}
	\begin{aligned}
	= & \sum_{i,j,k,\hbar,\imath,\alpha}B_{i\hbar}^{\alpha}B_{k\imath}^{\alpha}h\big (du(e_{\hbar}),du(e_{k})\big )h\big (du(e_{\imath}),du(e_{j})\big )h\big (du(e_{i}),du(e_{j})\big )\\
	& + \sum_{i,j,k,\hbar}h\big(\widetilde{\nabla}_{e_{i}}du(e_{\hbar}),du(e_{k})\big)h\big(\widetilde{\nabla}_{e_{k}}du(e_{\hbar}),du(e_{j})\big)h\big (du(e_{i}),du(e_{j})\big )\\
	= & \sum_{i,j,k.\alpha}h\big (du(A^{e_\alpha}(e_{i})),du(e_{k})\big )h\big (du(A^{e_\alpha}(e_{k})),du(e_{j})\big )h\big (du(e_{i}),du(e_{j})\big )\\
	& + \sum_{i,j,k,\hbar}h\big((\nabla_{e_{\hbar}}du)(e_{i}),du(e_{k})\big)h\big((\nabla_{e_{\hbar}}du)(e_{k}),du(e_{j})\big)h\big (du(e_{i}),du(e_{j})\big ).
	\end{aligned}
	\label{eq:SS6}
	\end{equation}
	
	The fifth integrand in (\ref{eq:SS2}) is
	\begin{equation*}
	\begin{aligned}
	& \sum_{i,j,k,\ell}h\big (\widetilde{\nabla}_{e_{i}}du(v_{\ell}^\top),du(e_{k})\big)h\big(du(e_{k}),\widetilde{\nabla}_{e_{j}}du(v_{\ell}^\top)\big)h\big (du(e_{i}),du(e_{j})\big )
	\end{aligned}
	\end{equation*}
	\begin{equation} 	\label{eq:SS7}
	\begin{aligned}
	= & \sum_{i,j,k,\hbar,\imath,\alpha,\beta}h\big(B_{i\hbar}^{\alpha}du(e_{\hbar})+\widetilde{\nabla}_{e_{i}}du(e_{\hbar}),du(e_{k})\big)\\
	& \cdot h\big (du(e_{k}),B_{j\imath}^{\beta}du(e_{\imath})+\widetilde{\nabla}_{e_{j}}du(e_{\imath})\big )h\big (du(e_{i}),du(e_{j})\big )\\
	= & \sum_{i,j,k,\hbar, \imath,\alpha}B_{i\hbar}^{\alpha} B_{j\imath}^{\alpha}h\big (du(e_{\hbar}),du(e_{k})\big )h\big (du(e_{k}),du(e_{\imath})\big ) h\big ( du(e_{i}),du(e_{j})\big )\\
	& + \sum_{i,j,k,\hbar}h\big(\widetilde{\nabla}_{e_{i}}du(e_{\hbar}),du(e_{k})\big)h\big(du(e_{k}),\widetilde{\nabla}_{e_{j}}du(e_{\hbar})\big)h\big (du(e_{i}),du(e_{j})\big )\\
	= & \sum_{i,j,k.\alpha}h\big (du(A^{e_\alpha}(e_{i})),du(e_{k})\big )h\big (du(e_{k}),du(A^{e_\alpha}(e_{j}))\big )h\big (du(e_{i}),du(e_{j})\big )\\
	& + \sum_{i,j,k,\hbar}h\big((\nabla_{e_{\hbar}}du)(e_{i}),du(e_{k})\big)h\big(du(e_{k}),(\nabla_{e_{\hbar}}du)(e_{j})\big)h\big (du(e_{i}),du(e_{j})\big ).
	\end{aligned}
	\end{equation}
	
	The sixth integrand in (\ref{eq:SS2}) is	
	\begin{equation}
	\begin{aligned}
	& \sum_{i,j,k,\ell}h\big (\widetilde{\nabla}_{e_{i}}du(v_{\ell}^\top),du(e_{k})\big)h\big (du(e_{k}),du(e_{j})\big )h\big(\widetilde{\nabla}_{e_{i}}du(v_{\ell}^\top),du(e_{j})\big)\\
	= & \sum_{i,j,k,\hbar,\alpha}h\big(B_{i\hbar}^{\alpha}du(e_{\hbar})+\widetilde{\nabla}_{e_{i}}du(e_{\hbar}),du(e_{k})\big)h\big (du(e_{k}),du(e_{j})\big )\\
	& \cdot h\big(B_{i\hbar}^{\alpha}du(e_{\hbar})+\widetilde{\nabla}_{e_{i}}du(e_{\hbar}),du(e_{j})\big)\\
	= & \sum_{i,j,k,\hbar.\alpha}B_{i\hbar}^{\alpha}B_{i\hbar}^{\alpha}h\big (du(e_{\hbar}),du(e_{k})\big )h\big (du(e_{k}),du(e_{j})\big )h\big (du(e_{\hbar}),du(e_{j})\big )\\
	& + \sum_{i,j,k,\hbar}h\big(\widetilde{\nabla}_{e_{i}}du(e_{\hbar}),du(e_{k})\big)h\big (du(e_{k}),du(e_{j})\big )h\big(\widetilde{\nabla}_{e_{i}}du(e_{\hbar}),du(e_{j})\big)\\
	%
	= & \sum_{i,j,k.\alpha}h\big (du(A^{e_\alpha}(e_{i})),du(e_{k})\big )h\big (du(e_{k}),du(e_{j})\big )h\big (du(A^{e_\alpha}(e_{i})),du(e_{j})\big )\\
	& + \sum_{i,j,k,\hbar}h\big((\nabla_{e_{\hbar}}du)(e_{i}),du(e_{k})\big)h\big (du(e_{k}),du(e_{j})\big )h\big((\nabla_{e_{\hbar}}du)(e_{i}),du(e_{j})\big).
	\end{aligned}
	\label{eq:SS8}
	\end{equation}
	
	The seventh integrand in (\ref{eq:SS2}) is
	\begin{equation}
	\begin{aligned}
	&\sum_{i,j,k,\ell}h\big(\widetilde{\nabla}_{e_{i}}du(v_{\ell}^\top),du(e_{k})\big)h \big (du(e_{k}),du(e_{j})\big )h\big(du(e_{i}),\widetilde{\nabla}_{e_{j}}du(v_{\ell}^\top)\big)\\
	= & \sum_{i,j,k,\hbar,\imath,\alpha,\beta}h\big(B_{i\hbar}^{\alpha}du(e_{\hbar})+\widetilde{\nabla}_{e_{i}}du(e_{\hbar}),du(e_{k})\big)h\big (du(e_{k}),du(e_{j})\big )\\
	& \cdot  h\big(du(e_{i}),B_{j\imath}^{\beta}du(e_{\imath})+\widetilde{\nabla}_{e_{j}}du(e_{\imath})\big)\\
	= & \sum_{i,j,k,\hbar.\imath,\alpha}B_{i\hbar}^{\alpha}B_{j\imath}^{\alpha}h\big (du(e_{\hbar}),du(e_{k})\big )h\big (du(e_{k}),du(e_{j})\big )h\big (du(e_{i}),du(e_{\imath})\big ) \\
	& + \sum_{i,j,k,\hbar}h\big(\widetilde{\nabla}_{e_{i}}du(e_{\hbar}),du(e_{k})\big)h\big (du(e_{k}),du(e_{j})\big )h\big(du(e_{i}),\widetilde{\nabla}_{e_{j}}du(e_{\hbar})\big)\\
	= & \sum_{i,j,k.\alpha}h\big (du(A^{e_\alpha}(e_{i})),du(e_{k})\big )h\big (du(e_{k}),du(e_{j})\big )h\big (du(e_{i}),du(A^{e_\alpha}(e_{j}))\big )\\
	& + \sum_{i,j,k,\hbar}h\big((\nabla_{e_{\hbar}}du)(e_{i}),du(e_{k})\big)h\big (du(e_{k}),du(e_{j})\big )h\big(du(e_{i}),(\nabla_{e_{\hbar}}du)(e_{j})\big).
	\end{aligned}
	\label{eq:SS9}
	\end{equation}
	
	Combining (\ref{eq:SS2})-($\ref{eq:SS9}$), we have
	\begin{eqnarray} \label{eq:SS10}
	\begin{aligned}
	&\sum_{\ell=1}^qI\big (du(v_{\ell}^\top),du(v_{\ell}^\top)\big )\\
	= & -\int_{M}\sum_{i,\alpha}h\big (du(\mathrm{trace}(A^{e_\alpha})A^{e_\alpha}(e_{i})), d_{(3)} u(e_{i})\big )dv_{g}\\
	& + \int_{M}\sum_{i,\alpha}h\big (du(A^{e_\alpha}A^{e_\alpha}(e_{i})), d_{(3)} u(e_{i})\big )dv_{g}\\
	& + \int_{M}\sum_{i,j,k.\alpha}h\big (du(A^{e_\alpha}(e_{i})),du(A^{e_\alpha}(e_{k}))\big )h\big (du(e_{k}),du(e_{j})\big )h\big (du(e_{i}),du(e_{j})\big )dv_{g}\\
	& + \int_{M}\sum_{i,j,k.\alpha}h\big (du(A^{e_\alpha}(e_{i})),du(e_{k})\big )h\big (du(A^{e_\alpha}(e_{k})),du(e_{j})\big )h\big (du(e_{i}),du(e_{j})\big )dv_{g}\\
	& + \int_{M}\sum_{i,j,k.\alpha}h\big (du(A^{e_\alpha}(e_{i})),du(e_{k})\big )h\big (du(e_{k}),du(A^{e_\alpha}(e_{j}))\big )h\big (du(e_{i}),du(e_{j})\big )dv_{g}\\
	& + \int_{M}\sum_{i,j,k.\alpha}h\big (du(A^{e_\alpha}(e_{i})),du(e_{k})\big )h\big (du(e_{k}),du(e_{j})\big )h\big (du(A^{e_\alpha}(e_{i}),du(e_{j}))\big )dv_{g}\\	
	&+ \int_{M}\sum_{i,j,k.\alpha}h\big (du(A^{e_\alpha}(e_{i})),du(e_{k})\big )h\big (du(e_{k}),du(e_{j})\big )h\big (du(e_{i}),du(A^{e_\alpha}(e_{j}))\big )dv_{g}.
	\end{aligned}
	\end{eqnarray}

	As the matrix $h\big (du(e_{i}),du(e_{j})\big )$ is symmetric, we take a local orthonormal frame
	$\{e_{i}\}_{i=1}^{m}$ such that
	\begin{equation}
	h\big (du(e_{i}),du(e_{j})\big )=\lambda_{i}^{2}\delta_{ij},\qquad i,j=1,\cdots,m.\label{eq:SS11}
	\end{equation}
	
	The first integrand in ($\ref{eq:SS10}$) is
	\begin{eqnarray}\label{eq:SS12}
	&&\quad\sum_{i,\alpha}h\big (du\big (\mathrm{trace}(A^{e_\alpha})A^{e_\alpha}(e_{i})\big ), d_{(3)} u(e_{i})\big )\nonumber\\
	&&=\sum_{i,j,k,\alpha}h\big (du\big (\mathrm{trace}(A^{e_\alpha})A^{e_\alpha}(e_{i})\big ),du(e_{k})\big )h\big (du(e_{i}),du(e_{j})\big )h\big (du(e_{j}),du(e_{k})\big )\nonumber\\
	&&=\sum_{i,j,k,\hbar,\imath,\alpha}\langle A^{e_\alpha}(e_{\hbar}),e_{\hbar}\rangle\langle A^{e_\alpha}(e_{i}),e_{\imath}\rangle h\big (du(e_{\imath}),du(e_{k})\big )h\big (du(e_{i}),du(e_{j})\big )h\big (du(e_{j}),du(e_{k})\big )\nonumber\\
	&&=\sum_{i,j,k,\hbar,\imath}\langle B(e_{\hbar},e_{\hbar}),B(e_{i},e_{\imath})\rangle\lambda_{\imath}^{2}\delta_{\imath k}\lambda_{i}^{2}\delta_{ij}\lambda_{j}^{2}\delta_{jk}\\
	&&=\sum_{i,j}\lambda_{i}^{6}\langle B(e_{i},e_{i}),B(e_{j},e_{j})\rangle.\nonumber
	\end{eqnarray}
	
	The second integrand in ($\ref{eq:SS10}$) is
	\begin{eqnarray}\label{eq:SS13}
	&&~~~~\sum_{i,\alpha}h\big (du\big (A^{e_\alpha}A^{e_\alpha}(e_{i})\big ), d_{(3)} u(e_{i})\big )\nonumber\\
	&&=\sum_{i,j,k,\alpha}h\big (du\big (A^{e_\alpha}A^{e_\alpha}(e_{i})\big ),du(e_{k})\big )h\big (du(e_{i}),du(e_{j})\big )h\big (du(e_{j}),du(e_{k})\big )\nonumber\\
	&&=\sum_{i,j,k,\hbar,\imath,\alpha}\langle A^{e_\alpha}(e_{i}),e_{\imath}\rangle\langle A^{e_\alpha}(e_{\imath}),e_{\hbar}\rangle h\big (du(e_{\hbar}),du(e_{k})\big )h\big (du(e_{i}),du(e_{j})\big )h\big (du(e_{j}),du(e_{k})\big )\nonumber\\
	&&=\sum_{i,j,k,\hbar,\imath}\langle B(e_{i},e_{\imath}),B(e_{\imath},e_{\hbar})\rangle\lambda_{\hbar}^{2}\delta_{\hbar k}\lambda_{i}^{2}\delta_{ij}\lambda_{j}^{2}\delta_{jk}\\
	&&=\sum_{i,j}\lambda_{i}^{6}\langle B(e_{i},e_{j}),B(e_{i},e_{j})\rangle\nonumber.
	\end{eqnarray}
	
	The third integrand in ($\ref{eq:SS10}$) is
	\begin{eqnarray}\label{eq:SS14}
	&&~~~~\sum_{i,j,k.\alpha}h\big (du\big (A^{e_\alpha}(e_{i})\big ),du\big (A^{e_\alpha}(e_{k})\big )\big )h\big (du(e_{k}),du(e_{j})\big )h\big (du(e_{i}),du(e_{j})\big )\nonumber\\
	&&=\sum_{i,j,k,\hbar,\imath,\alpha}\langle A^{e_\alpha}(e_{i}),e_{l}\rangle\langle A^{e_\alpha}(e_{k}),e_{r}\rangle h\big (du(e_{\hbar}),du(e_{\imath})\big )h\big (du(e_{k}),du(e_{j})\big )h\big (du(e_{i}),du(e_{j})\big )\nonumber\\
	&&=\sum_{i,j,k,\hbar,r}\langle B(e_{i},e_{\hbar}),B(e_{k},e_{\imath})\rangle\lambda_{\hbar}^{2}\delta_{\hbar r}\lambda_{k}^{2}\delta_{kj}\lambda_{i}^{2}\delta_{ij}\nonumber\\
	&&=\sum_{i,j}\lambda_{i}^{4}\lambda_{j}^{2}\langle B(e_{i},e_{j}),B(e_{i},e_{j})\rangle\\
	&&\leq\sum_{i,j} (\frac{2}{3}\lambda_{i}^{6}+\frac{1}{3}\lambda_{j}^{6}\big )\langle B(e_{i},e_{j}),B(e_{i},e_{j})\rangle\nonumber\\
	&&=\sum_{i,j}\lambda_{i}^{6}\langle B(e_{i},e_{j}),B(e_{i},e_{j})\rangle.\nonumber
	\end{eqnarray}
	
	The fourth integrand in ($\ref{eq:SS10}$) is
	\begin{eqnarray}\label{eq:SS15}
	&&~~~~\sum_{i,j,k.\alpha}h\big (du\big (A^{e_\alpha}(e_{i})\big ),du(e_{k})\big )h\big (du(A^{e_\alpha}(e_{k})),du(e_{j})\big )h\big (du(e_{i}),du(e_{j})\big )\nonumber\\
	%
	&&=\sum_{i,j,k,\hbar,\imath,\alpha}\langle A^{e_\alpha}(e_{i}),e_{\hbar}\rangle\langle A^{e_\alpha}(e_{k}),e_{\imath}\rangle h\big (du(e_{\hbar}),du(e_{k})\big )h\big (du(e_{\imath}),du(e_{j})\big )h\big (du(e_{i}),du(e_{j})\big )\nonumber\\
	&&=\sum_{i,j,k,\hbar,\imath}\langle B(e_{i},e_{\hbar}),B(e_{k},e_{\imath})\rangle\lambda_{\hbar}^{2}\delta_{\hbar k}\lambda_{\imath}^{2}\delta_{rj}\lambda_{i}^{2}\delta_{ij}\nonumber\\
	&&=\sum_{i,j}\lambda_{i}^{4}\lambda_{j}^{2}\langle B(e_{i},e_{j}),B(e_{i},e_{j})\rangle\\
	&&\leq\sum_{i,j}\big (\frac{2}{3}\lambda_{i}^{6}+\frac{1}{3}\lambda_{j}^{6}\big)\langle B(e_{i},e_{j}),B(e_{i},e_{j})\rangle\nonumber\\
	&&=\sum_{i,j}\lambda_{i}^{6}\langle B(e_{i},e_{j}),B(e_{i},e_{j})\rangle.\nonumber
	\end{eqnarray}
	
	The fifth integrand in ($\ref{eq:SS10}$) is
	\begin{eqnarray}\label{eq:SS16}
	&&~~~~\sum_{i,j,k.\alpha}h\big (du\big (A^{e_\alpha}(e_{i})\big ),du(e_{k})\big )h\big (du(e_{k}),du\big (A^{e_\alpha}(e_{j})\big )\big )h\big (du(e_{i}),du(e_{j})\big )\nonumber\\
	&&=\sum_{i,j,k,\hbar,\imath,\alpha}\langle A^{e_\alpha}(e_{i}),e_{\hbar}\rangle\langle A^{e_\alpha}(e_{j}),e_{\imath}\rangle h\big (du(e_{\hbar}),du(e_{k})\big )h\big (du(e_{k}),du(e_{\imath})\big )h\big (du(e_{i}),du(e_{j})\big )\nonumber\\
	&&=\sum_{i,j,k,\hbar,\imath}\langle B(e_{i},e_{\hbar}),B(e_{j},e_{\imath})\rangle\lambda_{\hbar}^{2}\delta_{\hbar k}\lambda_{k}^{2}\delta_{k\imath}\lambda_{i}^{2}\delta_{ij}\nonumber\\
	&&=\sum_{i,j}\lambda_{i}^{2}\lambda_{j}^{4}\langle B(e_{i},e_{j}),B(e_{i},e_{j})\rangle\nonumber\\
	&&\leq\sum_{i,j}\big (\frac{1}{3}\lambda_{i}^{6}+\frac{2}{3}\lambda_{j}^{6}\big )\langle B(e_{i},e_{j}),B(e_{i},e_{j})\rangle\\
	&&=\sum_{i,j}\lambda_{i}^{6}\langle B(e_{i},e_{j}),B(e_{i},e_{j})\rangle.\nonumber
	\end{eqnarray}
	
	The sixth integrand in ($\ref{eq:SS10}$) is
	\begin{eqnarray}\label{eq:SS17}
	&&~~~~\sum_{i,j,k.\alpha}h\big ( du\big (A^{e_\alpha}(e_{i})\big ),du(e_{k})\big )h\big (du(e_{k}),du(e_{j})\big )h\big (du\big (A^{e_\alpha}(e_{i})\big ),du(e_{j})\big )\nonumber\\	
	&&=\sum_{i,j,k,\hbar,\imath,\alpha}\langle A^{e_\alpha}(e_{i}),e_{\hbar}\rangle\langle A^{e_\alpha}(e_{i}),e_{\imath}\rangle h\big (du(e_{\hbar}),du(e_{k})\big )h\big (du(e_{k}),du(e_{j})\big )h\big (du(e_{\imath}),du(e_{j})\big )\nonumber\\
	&&=\sum_{i,j,k,\hbar,\imath}\langle B(e_{i},e_{\hbar}),B(e_{i},e_{\imath})\rangle\lambda_{\hbar}^{2}\delta_{\hbar k}\lambda_{k}^{2}\delta_{kj}\lambda_{\imath}^{2}\delta_{\imath j}\\
	&&=\sum_{i,j}\lambda_{i}^{6}\langle B(e_{i},e_{j}),B(e_{i},e_{j})\rangle.\nonumber
	\end{eqnarray}
	
	The seventh integrand in ($\ref{eq:SS10}$) is
	\begin{eqnarray}\label{eq:SS18}
	&&~~~~\sum_{i,j,k.\alpha}h\big (du\big (A^{e_\alpha}(e_{i})\big ),du(e_{k})\big )h\big (du(e_{k}),du(e_{j})\big )h\big (du(e_{i}),du\big (A^{e_\alpha}(e_{j})\big )\big )\nonumber\\
	&&=\sum_{i,j,k,\hbar,\imath,\alpha}\langle A^{e_\alpha}(e_{i}),e_{\hbar}\rangle\langle A^{e_\alpha}(e_{j}),e_{\imath}\rangle h\big (du(e_{\hbar}),du(e_{k})\big )h\big (du(e_{k}),du(e_{j})\big )h\big (du(e_{i}),du(e_{\imath})\big )\nonumber\\
	&&=\sum_{i,j,k,\hbar,\imath}\langle B(e_{i},e_{l}),B(e_{j},e_{\imath})\rangle\lambda_{l}^{2}\delta_{\hbar k}\lambda_{k}^{2}\delta_{kj}\lambda_{i}^{2}\delta_{i\imath}\nonumber\\
	%
	&&=\sum_{i,j}\lambda_{i}^{2}\lambda_{j}^{4}\langle B(e_{i},e_{j}),B(e_{i},e_{j})\rangle\\
	&&\leq\sum_{i,j}\big (\frac{1}{3}\lambda_{i}^{6}+\frac{2}{3}\lambda_{j}^{6}\big)\langle B(e_{i},e_{j}),B(e_{i},e_{j})\rangle\nonumber\\
	&&=\sum_{i,j}\lambda_{i}^{6}\langle B(e_{i},e_{j}),B(e_{i},e_{j})\rangle.\nonumber
	\end{eqnarray}
	
	Combining ($\ref{eq:SS10}$)-($\ref{eq:SS18}$), we have
	\begin{equation}
	\begin{aligned} & \sum_{\ell=1}^q I\big (du(v_{\ell}^\top),du(v_{\ell}^\top)\big )\\
	\leq  & \int_{M}\sum_{i,j}\lambda_{i}^{6}\big (6\langle B(e_{i},e_{j}),B(e_{i},e_{j})\rangle-\langle B(e_{i},e_{i}),B(e_{j},e_{j})\rangle\big )dv_{g}.
	\end{aligned}
	\label{eq:SS19}
	\end{equation}
	
	If $u$ is not constant on $M,$ then
	\begin{equation*}
	\sum_{\ell=1}^q I\big (du(v_{\ell}^\top),du(v_{\ell}^\top)\big ) < 0,
	\end{equation*}
	that is, there exists a variational vector field $du(v_{\ell}^\top)$ along which
	decreases $E_{\Phi_{(3)}}$-energy for some $1\leq  \ell\leq  q$. Therefore $u$ is not be a stable $\Phi_{(3)}$-harmonic map. This contradiction proves that $u$ is constant.
\end{proof}

By the results of Section \ref{sec6} and Theorem \ref{thm:fss1}, we have
\begin{cor}
	If $M^m$ is a compact manifold and satisfies the condition of examples in Section \ref{sec6}, then every stable $\Phi_{(3)}$-harmonic
	map $u:(M^m,g)\rightarrow(N,h)$ is a constant map from $(M^m,g)$ into any compact Riemannian manifold $N$.
\end{cor}

\section{Stable $\Phi_{(3)}$-harmonic maps into $\Phi_{(3)}$-SSU manifolds} \label{sec8}

In this section, we prove the following theorem
\begin{thm} \label{thm:iss1}
	Suppose $(N^{n},g)$ is a compact $\Phi_{(3)}$-$\operatorname{SSU}$ manifold and $(M^{m},g)$
	is any compact manifold. Then every stable $\Phi_{(3)}$-harmonic
	map $u:(M^{m},g)\rightarrow(N^{n},h)$ is constant.
\end{thm}
\begin{proof}
	We choose a local orthonormal frame field $\{e_{1},\cdots,e_{m}\}$
	on $M$. Let $\mathsf v, \mathsf v^\top, \mathsf v^{\perp}$ denote a unit vector in $\mathbb R^q$ the tangential projection of $\mathsf v$ onto $N$, and the normal projection of $\mathsf v$ onto $N$ respectively. We can choose an adopted orthonormal basis $\{\mathsf v_{\ell}\}_{{\ell}=1}^q$ in $\mathbb R^q$ such that $\{\mathsf v_{\ell}\}_{{\ell}=1}^n$ is tangent to $N$, and $\{\mathsf v_{\ell}\}_{{\ell}=n+1}^q$ is normal to $N$ at a point in $N$. Denote by $\mathsf f_t^{\mathsf v_{\ell}^\top}$ the flow generated by $\mathsf v_{\ell}^\top$.  Set $du(e_i)=\sum_{\alpha=1}^nu_{i\alpha} \mathsf e_{\alpha}\quad \text{for}\quad 1 \le i \le m\, .$ As $\mathsf v_{\ell}$ is parallel in $\mathbb R^q$, we have
	\begin{equation}\label{3.9}
	\begin{aligned}
	\nabla ^u_{e_i}\mathsf v_{\ell}^\top&=\nabla ^N_{du(e_i)}\mathsf v_{\ell}^\top=\left(\nabla ^{\mathbb R ^q}_{du(e_i)}\mathsf v_{\ell}^\top\right)^\top
	=\left(\nabla ^{\mathbb R ^q}_{du(e_i)}(\mathsf v_{\ell}-\mathsf v_{\ell}^\bot)\right)^\top\\
	&=\mathsf A^{\mathsf v_{\ell}^\bot}(du(e_i))\, .
	\end{aligned}
	\end{equation}
	Hence, if $\mathsf v_{\ell}=  \mathsf e_\nu$ for some $\nu \ge n+1$ at a point in $N$, then \begin{equation}\nabla ^u_{e_i}\mathsf v_{\ell}^\top=\sum_{\alpha=1}^nu_{i\alpha} \mathsf B_{\alpha\beta} ^{\nu},\end{equation} where $\mathsf B_{\alpha\beta} ^{\nu}=\langle  \mathsf B (\mathsf e_{\alpha} , \mathsf e_{\beta} ), \mathsf e_\nu\rangle$ are the components of the 2nd fundamental form of $N$ in $\mathbb R^q\, .$
	
	%
	We note
	the matrix
	\begin{equation}
	\bigg (\sum_{i=1}^{m}u_{i\alpha}u_{i\beta}\bigg )_{\alpha,\beta=1,\cdots,n} = \bigg (u_{i\alpha}\bigg )^T \cdot \bigg (u_{i\alpha}\bigg )\label{eq:SS20}
	\end{equation}
	is symmetric. We take  local orthonormal frame fields $\{e_1,\cdots,e_m\}$ on $M$ and $\{\mathsf e_1,\cdots,\mathsf e_n\}$ on $N$ so that the symmrtic matrix is diagonalizable. Namely,
	\begin{equation}
	\sum_{i=1}^{m}u_{i\alpha}u_{i\beta}=\lambda_{\alpha}^{2}\delta_{\alpha\beta}.\label{eq:SSS20}
	\end{equation}
	
	Suppose that $i,j,k\in\{1,\cdots,m\}$, $\alpha,\beta,\gamma,\eta,\iota,\kappa,\sigma,\tau\in\{1,\cdots,n\}$, $\nu\in\{n+1,\cdots,q\}$, and $\ell \in \{1, \cdots, q\}.$ Using the second variation formula, and the extrinsic average variation method, we have
	\begin{equation}
	\begin{aligned}
	& \sum_{\ell=1}^qI\big (\mathsf v_{\ell}^\top,\mathsf v_{\ell}^\top\big )\\
	= & \int_{M}\sum_{i,\ell}h\big (R^{N}(\mathsf v_{\ell}^\top,du(e_{i}))\mathsf v_{\ell}^\top, d_{(3)} u(e_{i})\big )dv_{g}\\
	& +\int_{M}\sum_{i,j,k,\ell}h\big (\widetilde{\nabla}_{e_{i}}\mathsf v_{\ell}^\top,\widetilde{\nabla}_{e_{k}}\mathsf v_{\ell}^\top\big)h\big (du(e_{k}),du(e_{j})\big )h\big (du(e_{i}),du(e_{j})\big )dv_{g}\\
	& +\int_{M}\sum_{i,j,k,\ell}h\big(\widetilde{\nabla}_{e_{i}}\mathsf v_{\ell}^\top,du(e_{k})\big)h\big(\widetilde{\nabla}_{e_{k}}\mathsf v_{\ell}^\top,du(e_{j})\big)h\big (du(e_{i}),du(e_{j})\big )dv_{g}\\
	& +\int_{M}\sum_{i,j,k,\ell}h\big(\widetilde{\nabla}_{e_{i}}\mathsf v_{\ell}^\top,du(e_{k})\big )h\big(du(e_{k}),\widetilde{\nabla}_{e_{j}}\mathsf v_{\ell}^\top\big )h\big (du(e_{i}),du(e_{j})\big )dv_{g}\\
	& +\int_{M}\sum_{i,j,k,\ell}h\big (\widetilde{\nabla}_{e_{i}}\mathsf v_{\ell}^\top,du(e_{k})\big)h\big (du(e_{k}),du(e_{j})\big )h\big(\widetilde{\nabla}_{e_{i}}\mathsf v_{\ell}^\top,du(e_{j})\big)dv_{g}\\
	& +\int_{M}\sum_{i,j,k,\ell}h\big(\widetilde{\nabla}_{e_{i}}\mathsf v_{\ell}^\top,du(e_{k})\big)h\big (du(e_{k}),du(e_{j})\big )h\big(du(e_{i}),\widetilde{\nabla}_{e_{j}}\mathsf v_{\ell}^\top\big)dv_{g}.
	\end{aligned}
	\label{eq:SS21}
	\end{equation}

	Hence, at $x_0,$ we give estimates on every term in (\ref{eq:SS21}):\\
	\begin{equation} \label{eq:SS22}
	\begin{aligned}
	&\sum_{i,\ell}h\big (R^{N}\big (\mathsf v_{\ell}^\top,du(e_{i})\big )\mathsf v_{\ell}^\top, d_{(3)} u(e_{i})\big )\\	
	= & \sum_{i,j,k,\ell}h\big (R^{N}\big (\mathsf v_{\ell}^\top,du(e_{i})\big )\mathsf v_{\ell}^\top,du(e_{k})\big )h\big (du(e_{i}),du(e_{j})\big )h\big (du(e_{j}),du(e_{k})\big )\\
	= & \sum_{i,j,k,\alpha,\beta,\gamma,\eta,\iota,\kappa,\sigma,\tau} u_{i\beta}u_{k\kappa}u_{i\tau}u_{j\sigma}u_{j\eta}u_{k\gamma}h\big (R^{N}(\mathsf e_{\alpha},\mathsf e_{\beta})\mathsf e_ {\iota},\mathsf e_ {\kappa}\big )h(\mathsf e_ {\tau},\mathsf e_ {\sigma})h(\mathsf e_{\eta},\mathsf e_{\gamma})\\
	= & \sum_{\alpha,\beta,\gamma,\eta,\kappa,\sigma}\lambda_{\beta}^{2}\delta_{\beta \tau}\lambda_{\sigma}^{2}\delta_{\sigma\eta}\lambda_{\kappa}^{2}\delta_{\kappa\gamma}\delta_{\tau\sigma}\delta_{\eta\gamma}h\big (R^{N}(\mathsf e_{\alpha},\mathsf e_{\beta})\mathsf e_{\alpha},\mathsf e_{\kappa}\big )\\
	= & \sum_{\alpha,\beta}\lambda_{\alpha}^{6}h\big (R^{N}(\mathsf e_{\alpha},\mathsf e_{\beta})\mathsf e_{\alpha},\mathsf e_{\beta}\big )\\
	= & \sum_{\alpha,\beta}\lambda_{\alpha}^{6}\big (\langle \mathsf B(\mathsf e_{\alpha},\mathsf e_{\beta}),\mathsf B(\mathsf e_{\alpha},\mathsf e_{\beta})\rangle-\langle \mathsf B(\mathsf e_{\alpha},\mathsf e_{\alpha}),\mathsf B(\mathsf e_{\beta},\mathsf e_{\beta})\rangle\big ),
	\end{aligned}
	\end{equation}
	and
	\begin{equation}\label{eq:SS23}
	\begin{aligned}
	&\sum_{i,j,k,\ell}h\big(\widetilde{\nabla}_{e_{i}}\mathsf v_{\ell}^\top,\widetilde{\nabla}_{e_{k}}\mathsf v_{\ell}^\top\big)h\big (du(e_{k}),du(e_{j})\big )h\big (du(e_{i}),du(e_{j})\big )\\
	= & \sum_{i,j,k,\alpha,\beta,\ell}u_{i\alpha}u_{k\beta}h\big(\nabla_{\mathsf e_{\alpha}}v_{\ell}^\top,\nabla_{\mathsf e_{\beta}}v_{\ell}^\top\big)h\big (du(e_{k}),du(e_{j})\big )h\big (du(e_{i}),du(e_{j})\big )\\
	= & \sum_{i,j,k,\alpha,\beta,\gamma,\eta,\iota,\kappa,\sigma,\tau,\nu}u_{i\alpha}u_{k\beta} \mathsf B_{\alpha \gamma}^{\nu}\mathsf B_{\beta \eta}^{\nu}u_{k}^{\iota}u_{j\kappa}u_{i\sigma}u_{j\tau}h(\mathsf e_{\gamma},\mathsf e_{\eta})h(\mathsf e_{\iota},\mathsf e_{\kappa})h(\mathsf e_{\sigma},\mathsf e_{\tau})\\
	= & \sum_{\alpha,\beta,\gamma,\eta,\iota,\kappa,\sigma,\tau,\nu}\lambda_{\alpha}^{2}\lambda_{\kappa}^{2}\lambda_{\beta}^{2}\delta_{\beta\iota}\delta_{\kappa\tau}\delta_{\alpha\sigma}\delta_{\gamma\eta}\delta_{\iota\kappa}\delta_{\sigma\tau}\mathsf B_{\alpha\gamma}^{\nu}\mathsf B_{\beta\eta}^{\nu}\\
	= & \sum_{\alpha,\gamma,\nu}\lambda_{\alpha}^{6}\mathsf B_{\alpha\gamma}^{\nu}\mathsf B_{\alpha\gamma}^{\nu}\\
	= & \sum_{\alpha,\beta}\lambda_{\alpha}^{6}\langle \mathsf B(\mathsf e_{\alpha},\mathsf e_{\beta}),\mathsf B(\mathsf e_{\alpha},\mathsf e_{\beta})\rangle,
\end{aligned}
\end{equation}
and
	\begin{equation*}
     \begin{aligned}
	& \sum_{i,j,k,\ell}h\big ( \widetilde{\nabla}_{e_{i}}\mathsf v_{\ell}^\top,du(e_{k})\big )h\big (\widetilde{\nabla}_{e_{k}}\mathsf v_{\ell}^\top,du(e_{j})\big )h\big (du(e_{i}),du(e_{j})\big )\\
	= & \sum_{i,j,k,\alpha,\beta,\ell}u_{i\alpha}u_{k\beta}h\big (\nabla_{\mathsf e_{\alpha}}\mathsf v_{\ell}^\top,du(e_{k})\big )h\big (\nabla_{\mathsf e_{\beta}}v_{\ell}^\top,du(e_{j})\big )h\big (du(e_{i}),du(e_{j})\big )\\
	= & \sum_{i,j,k,\alpha,\beta,\gamma,\eta,\iota,\kappa,\sigma,\tau,\nu}u_{i\alpha}u_{j\beta}u_{k\eta}u_{k\iota}u_{i\tau}u_{j\sigma}\mathsf B_{\alpha\gamma}^{\nu}\mathsf B_{\beta\sigma}^{\nu}h(\mathsf e_{\gamma},\mathsf e_{\eta})h(\mathsf e_{\tau},\mathsf e_{\iota})h(\mathsf e_{\tau},\mathsf e_{\kappa})\\
	= & \sum_{\alpha,\beta,\gamma,\eta,\iota,\kappa,\sigma,\tau,\nu}\lambda_{\alpha}^{2}\lambda_{\iota}^{2}\lambda_{\beta}^{2}\delta_{\beta\tau}\delta_{\iota\kappa}\delta_{\beta\eta}\delta_{\gamma\eta}\delta_{\sigma\iota}\delta_{\kappa\tau}\mathsf B_{\alpha\gamma}^{\alpha}\mathsf B_{\beta\sigma}^{\nu}\\
	= & \sum_{\alpha,\beta,\nu}\lambda_{\alpha}^{4}\lambda_{\beta}^{2}B_{\alpha\beta}^{\nu}B_{\alpha\beta}^{\nu}\\
	= & \sum_{\alpha,\beta}\lambda_{\alpha}^{4}\lambda_{\beta}^{2}\langle \mathsf B(\mathsf e_{\alpha},\mathsf e_{\beta}),\mathsf B(\mathsf e_{\alpha},\mathsf e_{\beta})\rangle
\end{aligned}
	\end{equation*}
\begin{equation}
     \begin{aligned}
	\leq  & \sum_{\alpha,\beta}\big (\frac{2}{3}\lambda_{\alpha}^{6}+\frac{1}{3}\lambda_{\beta}^{6}\big)\langle \mathsf B(\mathsf e_{\alpha},\mathsf e_{\beta}),\mathsf B(\mathsf e_{\alpha},\mathsf e_{\beta})\rangle\\
	= & \sum_{\alpha,\beta}\lambda_{\alpha}^{6}\langle \mathsf B(\mathsf e_{\alpha},\mathsf e_{\beta}),\mathsf B(\mathsf e_{\alpha},\mathsf e_{\beta})\rangle,
	\end{aligned}
	\label{eq:SS24}
	\end{equation}
and	
	\begin{equation}
	\begin{aligned}
	& \sum_{i,j,k,\ell}h\big ( \widetilde{\nabla}_{e_{i}}\mathsf v_{\ell}^\top,du(e_{k})\big )h\big (du(e_{k}),\widetilde{\nabla}_{e_{j}}\mathsf v_{\ell}^\top\big )h(du(e_{i}),du(e_{j}))\\
	= & \sum_{i,j,k,\alpha,\beta,\ell}u_{i\alpha}u_{j}^{\beta}h\big (\nabla_{\mathsf e_{\alpha}}\mathsf v_{\ell}^\top,du(e_{k})\big )h\big (du(e_{k}),\nabla_{\mathsf e_{\beta}}\mathsf v_{\ell}^\top\big )h\big (du(e_{i}),du(e_{j})\big )\\
	= & \sum_{i,j,k,\alpha,\beta,\gamma,\eta,\iota,\kappa,\sigma,\tau,\nu}u_{i\alpha}u_{j\beta}u_{k\eta}u_{k\iota}u_{i\kappa}u_{j\sigma}\mathsf B_{\beta\gamma}^{\nu}\mathsf B_{\beta\tau}^{\nu}h(\mathsf e_{\gamma},\mathsf e_{\eta})h(\mathsf e_{\iota},\mathsf e_{\tau})h(\mathsf e_{\kappa},\mathsf e_{\sigma})\\
	= & \sum_{\alpha,\beta,\gamma,\eta,\iota,\kappa,\sigma,\tau,\nu}\lambda_{\alpha}^{2}\lambda_{\beta}^{2}\lambda_{\eta}^{2}\delta_{\alpha\kappa}\delta_{\beta\sigma}\delta_{\eta\iota}\delta_{\gamma\eta}\delta_{\iota\tau}\delta_{\kappa\sigma}\mathsf B_{\alpha\gamma}^{\nu}\mathsf B_{\beta\tau}^{\nu}\\
	= & \sum_{\alpha,\beta,\nu}\lambda_{\alpha}^{4}\lambda_{\beta}^{2}\mathsf B_{\alpha\beta}^{\nu}\mathsf B_{\alpha\beta}^{\nu}\\
	%
	= & \sum_{\alpha,\beta,\nu}\lambda_{\alpha}^{4}\lambda_{\beta}^{2}\langle \mathsf B(\mathsf e_{\alpha},\mathsf e_{\beta}),\mathsf B(\mathsf e_{\alpha},\mathsf e_{\beta})\rangle\\
	\leq  & \sum_{\alpha,\beta}\big(\frac{2}{3}\lambda_{\alpha}^{6}+\frac{1}{3}\lambda_{\beta}^{6}\big )\langle \mathsf B(\mathsf e_{\alpha},\mathsf e_{\beta}),\mathsf B(\mathsf e_{\alpha},\mathsf e_{\beta})\rangle\\
	= & \sum_{\alpha,\beta}\lambda_{\alpha}^{6}\langle \mathsf B(\mathsf e_{\alpha},\mathsf e_{\beta}),\mathsf B(\mathsf e_{\alpha},\mathsf e_{\beta})\rangle,
	\end{aligned}\label{eq:SS25}
	\end{equation}	
and	
	\begin{equation}
	\begin{aligned}
	&\sum_{i,j,k,\ell}h\big (\widetilde{\nabla}_{e_{i}}\mathsf \mathsf v_{\ell}^\top,du(e_{k})\big )h\big (du(e_{k}),du(e_{j})\big )h\big (\widetilde{\nabla}_{e_{i}}\mathsf \mathsf v_{\ell}^\top,du(e_{j})\big )\\
	= & \sum_{i,j,k,\alpha,\beta,\ell}u_{i\alpha}u_{i\beta}h\big (\nabla_{\mathsf e_{\alpha}}\mathsf v_{\ell}^\top,du(e_{k})\big )h\big (du(e_{k}),du(e_{j})\big )h\big (\nabla_{\mathsf e_{\beta}}\mathsf v_{\ell}^\top,du(e_{j})\big )\\
	= & \sum_{i,j,k,\alpha,\beta,\gamma,\eta,\iota,\kappa,\sigma,\tau,\nu}u_{i\alpha}u_{i\beta}u_{k\eta}u_{k\iota}u_{j\tau}u_{j\sigma}\mathsf B_{\alpha \gamma}^{\nu}\mathsf B_{\beta\kappa}^{\nu}h(\mathsf e_{\gamma},\mathsf e_{\eta})h(\mathsf e_{\iota},\mathsf e_{\tau})h(\mathsf e_{\kappa},\mathsf e_{\sigma})\\
	= & \sum_{\alpha,\beta,\gamma,\eta,\iota,\kappa,\sigma,\tau,\nu}\lambda_{\alpha}^{2}\lambda_{\eta}^{2}\lambda_{\tau}^{2}\delta_{\alpha\beta}\delta_{\eta\iota}\delta_{\tau\sigma}\delta_{\gamma\eta}\delta_{\iota\tau}\delta_{\kappa\sigma}\mathsf B_{\alpha\gamma}^{\nu}\mathsf B_{\alpha\kappa}^{\nu}\\
	= & \sum_{\alpha,\beta,\nu}\lambda_{\alpha}^{2}\lambda_{\beta}^{4}\mathsf B_{\alpha,\beta}^{\nu}\mathsf B_{\alpha,\beta}^{\nu}\\
	\leq  & \sum_{\alpha,\beta}\big (\frac{1}{3}\lambda_{\alpha}^{6}+\frac{2}{3}\lambda_{\beta}^{6}\big )\langle \mathsf B(\mathsf e_{\alpha},\mathsf e_{\beta}),\mathsf B(\mathsf e_{\alpha},\mathsf e_{\beta})\rangle\\
	= & \sum_{\alpha,\beta}\lambda_{\alpha}^{6}\langle \mathsf B(\mathsf e_{\alpha},\mathsf e_{\beta}),\mathsf B(\mathsf e_{\alpha},\mathsf e_{\beta})\rangle,
	\end{aligned}
	\label{eq:SS26}
	\end{equation}
	and
	\begin{equation}\label{eq:SS27}
	\begin{aligned}
	& \sum_{i,j,k,\ell}h\big ( \widetilde{\nabla}_{e_{i}}\mathsf v_{\ell}^\top,du(e_{k})\big )h\big (du(e_{k}),du(e_{j})\big )h\big (du(e_{i}),\widetilde{\nabla}_{e_{j}}\mathsf v_{\ell}^\top\big )\\
	= & \sum_{i,j,k,\alpha,\beta,\ell}u_{i\alpha}u_{j\beta}h\big (\nabla_{\mathsf e_{\alpha}} \mathsf v_{\ell}^\top,du(e_{k})\big )h\big (du(e_{k}),du(e_{j})\big )h\big (du(e_{i}),\nabla_{\mathsf e_{\beta}}\mathsf v_{\ell}^\top\big )\\
	= & \sum_{i,j,k,\alpha,\beta,\gamma,\eta,\iota,\kappa,\sigma,\tau,\nu}u_{i\alpha
}u_{j\beta}u_{k\eta}u_{k\iota}u_{j\tau}u_{i\kappa}\mathsf B_{\alpha\gamma}^{\nu}\mathsf B_{\beta\sigma}^{\nu}h(\mathsf e_{\gamma},\mathsf e_{\eta})h(\mathsf e_{\iota},\mathsf e_{\tau})h(\mathsf e_{\kappa},\mathsf e_{\sigma})\\
= & \sum_{\alpha,\beta,\gamma,\eta,\iota,\kappa,\sigma,\tau,\nu}\lambda_{\alpha}^{2}\lambda_{\beta}^{2}\lambda_{\eta}^{2}\delta_{\alpha\kappa}\delta_{\beta\tau}\delta_{\eta\iota}\delta_{\gamma\eta}\delta_{\iota\tau}\delta_{\kappa\sigma}\mathsf B_{\alpha\gamma}^{\nu}\mathsf B_{\beta\sigma}^{\nu}\\
= & \sum_{\alpha,\beta,\nu}\lambda_{\alpha}^{2}\lambda_{\beta}^{4}\mathsf B_{\alpha\beta}^{\nu}\mathsf B_{\alpha\beta}^{\nu}\\
\leq  & \sum_{\alpha,\beta}\big (\frac{1}{3}\lambda_{\alpha}^{6}+\frac{2}{3}\lambda_{\beta}^{6}\big )\langle B(\mathsf e_{\alpha},\mathsf e_{\beta}),B(\mathsf e_{\alpha},\mathsf e_{\beta})\rangle\\
= & \sum_{\alpha,\beta}\lambda_{\alpha}^{6}\langle \mathsf B(\mathsf e_{\alpha},\mathsf e_{\beta}),\mathsf B(\mathsf e_{\alpha},\mathsf e_{\beta})\rangle.
\end{aligned}
\end{equation}

Combining ($\ref{eq:SS20}$)-($\ref{eq:SS27}$), we have

\begin{equation}
\begin{aligned}  \sum_{\ell=1}^q I\big (\mathsf v_{\ell}^\top,\mathsf v_{\ell}^\top\big )
\leq \int_{M}\sum_{\alpha,\beta}\lambda_{\alpha}^{6}\big (6\langle \mathsf B(\mathsf e_{\alpha},\mathsf e_{\beta}),\mathsf B(\mathsf e_{\alpha},\mathsf e_{\beta})\rangle-\langle \mathsf B(\mathsf e_{\alpha},\mathsf e_{\alpha}),\mathsf B(\mathsf e_{\beta},\mathsf e_{\beta})\rangle\big )dv_{g}.
\end{aligned}
\label{eq:SS28}
\end{equation}

As $N$ is a $\Phi_{(3)}$-SSU manifold, we know that if $u$ is a
nonconstant, then
\begin{equation*}
\sum_{\ell=1}^qI\big (\mathsf v_{\ell}^\top,\mathsf v_{\ell}^\top\big ) < 0.
\end{equation*}

That is, making a variation of $u$ along a vector field $\mathsf v_{\ell}^\top$ decreases $E_{\Phi_{(3)}}$-energy
for some $1\leq  \ell\leq  q$. Hence, $u$ is
not a stable $\Phi_{(3)}$-harmonic map, which leads to a contradiction. Consequently, $u$ is a constant.
\end{proof}

By Section \ref{sec6} and Theorem \ref{thm:iss1}, we have
\begin{cor}
If $N^n$ is a compact manifold and satisfies the condition of examples in Section \ref{sec6}, $M^m$ is any compact Riemannian manifold, then every stable $\Phi_{(3)}$-harmonic
map $u:(M^m,g)\rightarrow(N^n,h)$ is a constant map.
\end{cor}

\section{The Infimum of $\Phi_{(3)}$-energy in the homotopic class of maps into $\Phi_{(3)}$-SSU manifolds}
\label{sec9}

\begin{lem}\label{L:4.3}
If $N$ is a compact $\Phi_{(3)}$-$\operatorname{SSU}$ manifold, then there is a number $0<\rho<1$ such that for any compact manifold $M$ and any map $u:M\rightarrow N$ there is a map $u_1:M\rightarrow N$ homotopic to $u$ with $E_{\Phi_{(3)}}(u_1)\leq \rho E_{\Phi_{(3)}}(u)$.
\end{lem}
\begin{proof}
Let $\widehat{T}_yN$ be the space of the unit tangent vectors to $N$ at the point $y \in N$. Since $N$ is $\Phi_{(3)}$-SSU and $\widehat{T}_yN$ is compact, by \eqref{eq:SSU}, there exists $\kappa > 0$ such that for every $y \in N\, $ and every $ x \in \widehat{T}_yN\, ,$
\begin{equation}\label{4.23}
{F_{\Phi_{(3)}}}_x (v) < - q \kappa.
\end{equation}

Similar to the method of \cite{hw}, let $v$ be a unit vector in  $\mathbb{R}^{q}$,

It follows from ($\ref{eq:SSU}$), ($\ref{eq:SS28}$) and \eqref{4.23}  that
\begin{equation}
\sum_{\ell=1}^q\frac{d^2}{dt^2}E_{\Phi_{(3)}}(f_t^{v_{\ell}^\top}\circ u)_{\big |_{t=0}}
\le  -6q\kappa\int_M e_{\Phi_{(3)}} (u)\, dv_{g}=-6q\kappa E_{\Phi_{(3)}}(u).\label{4.24}
\end{equation}
We now proceed in steps.

{\bf{Step 1}}. There is a number $\xi\geq 6\kappa>0$ such that for $1\leq \ell\leq q $, $|t|\leq 1$ and all $ X, Y, Z\in \Gamma(TN)$,
\begin{equation}
\begin{aligned}
&\bigg | \frac{d^3}{dt^3} \bigg (\langle d f_t^{v_{\ell}^\top}( X),d f_t^{ v_{\ell}^\top}( Y)\rangle \langle d f_t^{v_{\ell}^\top}( Y),d f_t^{ v_{\ell}^\top}( Z)\rangle \langle d f_t^{v_{\ell}^\top}( Z),d f_t^{ v_{\ell}^\top}( X)\rangle \bigg )\bigg |\\
\le &\xi \bigg |\langle X,Y\rangle \langle Y, Z\rangle \langle  Z,X\rangle \bigg |.\label{4.25}
\end{aligned}
\end{equation}
\begin{proof}
	Let $SN$ be the unit sphere bundle of $N$. Then the function defined on the compact set $[-1,1]\times SN\times SN \times SN$ by
	\begin{eqnarray}\label{4.26}
	\begin{aligned}
	(t,  x,  y, z)\mapsto \max_{1\leq \ell\leq q}\bigg | \frac{d^3}{dt^3} \big (\langle d f_t^{v_{\ell}^\top}( x),d f_t^{ v_{\ell}^\top}( y)\rangle \langle d f_t^{v_{\ell}^\top}( y),d f_t^{ v_{\ell}^\top}( z)\rangle \langle d f_t^{v_{\ell}^\top}( z),d f_t^{ v_{\ell}^\top}( x)\rangle \big )\bigg |
	\end{aligned}
	\end{eqnarray}
	is continuous and thus has a maximum. Let ${\xi}_0$ be this maximum and $\xi=\max\{5\kappa, \xi_0\}$. Then \eqref{4.25} follows by homogeneity.
\end{proof}
{\bf{Step 2}}. There is a smooth vector field $V$ on $N$ such that if $\xi$ is as in Step 1,
then we have
\begin{equation}
\frac{d}{dt}E_{\Phi_{(3)}}(f_t^{V}\circ u)_{\big |_{t=0}}\leq 0,\label{4.27}
\end{equation}
\begin{equation}
\frac{d^2}{dt^2}E_{\Phi_{(3)}}(f_t^{ V}\circ u)_{\big |_{t=0}}\leq - 6\kappa E_{\Phi_{(3)}}(u),\label{4.28}
\end{equation}
and
\begin{equation}
\left|\frac{d^3}{dt^3}E_{\Phi_{(3)}}(f_t^{V}\circ u)\right| \leq \xi E_{\Phi_{(3)}}(u),\quad{\text{for}}\quad |t|\leq 1.\label{4.29}
\end{equation}
\begin{proof}
	From \eqref{4.24} it is seen that
	\begin{equation}\label{4.30}
	\frac{d^2}{dt^2}E_{\Phi_{(3)}}(f_t^{ v_{\ell}^\top}\circ u)_{\big |_{t=0}}\leq - 6\kappa E_{\Phi_{(3)}}(u),
	\end{equation}
	for some $1 \le \ell \le q$. Or we would have $\sum_{\ell=1}^q\frac{d^2}{dt^2}E_{\Phi_{(3)}}(f_t^{v_{\ell}^\top}\circ u)_{\big |_{t=0}}  > -6q\kappa E_{\Phi_{(3)}}(u)\, ,$ contradicting \eqref{4.24}. If $\frac{d}{dt}E_{\Phi_{(3)}}( f_t^{ v_{\ell}^\top}\circ u)_{\big |_{t=0}}\leq 0,$ set $V=  v_{\ell}^\top$; otherwise, set $ V=- v_{\ell}^\top$. Then \eqref{4.27} and \eqref{4.28} hold. From \eqref{4.25}, we have
	\begin{eqnarray}\label{4.31}
	&& \quad  \left|\frac{d^3}{dt^3}E_{\Phi_{(3)}}(f_t^{V}\circ u)\right| \nonumber\\
	&& = \frac 16 \int_M\sum_{i,j,k=1}^m \frac{d^3}{dt^3} \bigg (\langle df_t^{v_{\ell}^\top}\big (du(e_i)\big ),df_t^{v_{\ell}^\top}\big (du(e_j)\big )\rangle \langle df_t^{v_{\ell}^\top}\big (du(e_j)\big ),df_t^{v_{\ell}^\top}\big (du(e_k)\big )\rangle \nonumber\\
	&& \quad \times \langle df_t^{v_{\ell}^\top}\big (du(e_k)\big ),df_t^{v_{\ell}^\top}\big (du(e_i)\big )\rangle\bigg )\, dv_{g}\nonumber\\
	&&\leq \frac {\xi}{6} \int_M\sum_{i,j,k=1}^m\langle du(e_i), du(e_j)\rangle \langle du(e_j), du(e_k) \rangle \langle du(e_k), du(e_i)\rangle \, dv_{g}\\
	&&= \xi E_{\Phi_{(3)}}(u).\nonumber
	\end{eqnarray}
	So \eqref{4.29} is right.
\end{proof}
{\bf{Step 3}}. Let $\zeta=\frac{5\kappa}{\xi}$ ($\zeta\leq 1$, as $5\kappa\leq \xi$), $\rho=1-\frac{\kappa \zeta^2}{2}$, and $V$ be as in Step 2. Then $0<\rho<1$ and
\begin{align}\label{4.32}
E_{\Phi_{(3)}}( f_{\zeta}^{ V}\circ u)\leq \rho E_{\Phi_{(3)}}(u).
\end{align}
\begin{proof}
	Let $E_{\Phi_{(3)}}(t)=E_{\Phi_{(3)}}( f_t^ V\circ u)$. Then by Step 2 for $0\leq t\leq \zeta$, we have
	\begin{align*}
	&E_{\Phi_{(3)}}^{\prime \prime}(t)=E_{\Phi_{(3)}}^{\prime \prime}(0)+\int_0^tE_{\Phi_{(3)}}^{\prime \prime \prime}(s)ds\leq -6\kappa E_{\Phi_{(3)}}(u) + \xi \zeta E_{\Phi_{(3)}}(u)=-\kappa E_{\Phi_{(3)}}(u).
	\end{align*}
	Thus
	\begin{align}\label{4.33}
	E_{\Phi_{(3)}}^{\prime}(t)=E_{\Phi_{(3)}}^{\prime}(0)+\int_0^tE_{\Phi_{(3)}}^{\prime \prime}(s)ds\leq -\kappa tE_{\Phi_{(3)}}(u)
	\end{align}
	and
	\begin{align}\label{4.34}
	E_{\Phi_{(3)}}(\zeta)=E_{\Phi_{(3)}}(0)+\int_0^{\zeta}E_{\Phi_{(3)}}^{\prime}(s)ds\leq \bigg(1-\frac{\kappa\zeta^2}{2}\bigg )E_{\Phi_{(3)}}(u)=\rho E_{\Phi_{(3)}}(u).
	\end{align}
	As both $E_{\Phi_{(3)}}(\zeta)$ and $E_{\Phi_{(3)}}(u)$ are positive, this inequality implies $\rho$ is positive. Let $u_1= f_{\zeta}^{ V}\circ u$. Then $u_1$ is homotopic to $u$ and we have just shown $E_{\Phi_{(3)}}(u_1)\leq \rho E_{\Phi_{(3)}}(u)$.
\end{proof}
From Step 1, Step 2 and Step 3, we know this lemma is right.
\end{proof}

\begin{thm}
\noindent
If $N$ is a compact $\Phi_{(3)}$-$\operatorname{SSU}$ manifold, then for every compact Riemannian manifold $M$, the homotopic class of any map from $M$ into $N$ contains elements of arbitrarily small $\Phi_{(3)}$-energy.\label{T: 9.2}
\end{thm}

\begin{proof}
Let $u:M\rightarrow N$ be any smooth map from $M$ to $N$. By using Lemma \ref{L:4.3}, we can find a map $u_1:M\rightarrow N$ which is homotopic to $u$ with $E_{\Phi_{(3)}}(u_1)\leq \rho E_{\Phi_{(3)}}(u)$.
Another application of the lemma gives an $u_2$ homotopic to $u_1$ with $E_{\Phi_{(3)}}(u_2)\leq \rho E_{\Phi_{(3)}}(u_1)\leq \rho^2E_{\Phi_{(3)}}(u)$. By induction, there is $u_{\ell}$ $(\ell=1,2,\cdots)$ homotopic to $u$ with $E_{\Phi_{(3)}}(u_\ell)\leq \rho^\ell E_{\Phi_{(3)}}(u)$.
But $0<\rho<1$ whence $\lim_{\ell\rightarrow\infty}E_{\Phi_{(3)}}(u_{\ell})=0$ as required.
\end{proof}

\begin{cor}\label{C:4.1}
If $N$ is a compact $\Phi_{(3)}$-$\operatorname{SSU}$ manifold, then the infimum of the $\Phi_{(3)}$-energy $E_{\Phi_{(3)}}$ is zero among maps homotopic to the identity map on $N\, .$
\end{cor}
\begin{proof} This follows at once from Theorem \ref{T: 9.2}  by choosing $M=N$ and the smooth map to be the identity map on $N\, .$
\end{proof}

\section{The Infimum of $\Phi_{(3)}$-energy in the homotopic class of maps from $\Phi_{(3)}$-SSU manifolds}
\label{sec10}

\begin{lem}\label{L: 4.4}
	If $N$ is a compact Riemannian manifold such that the infimum of the $\Phi_{(3)}$-energy $E_{\Phi_{(3)}}$  is zero among maps homotopic to the identity and if $M$ is a compact Riemannian manifold, then
	the infimum of the $\Phi_{(3)}$-energy $E_{\Phi_{(3)}}$ is zero in each homotopy class of maps from $N$ to $M$.
\end{lem}
\begin{proof}
	
	Let $K, M, N$ be three Riemannian manifolds of dimensions $k, m$ and $n$ respectively, and $M$ be a compact manifold.
	Let $K \xrightarrow {\psi} N \xrightarrow {u} M$ be smooth maps.
	Denote a symmetric $n$ $\times$$n$-matrix $U$ with the $i$-$j$ entry $U_{ij}$ given by $U_{ij} = \langle du(e_i), du(e_j)\rangle_M\, $, and a symmetric $k$$\times$$k$-matrix $\Psi$ with the $i$-$j$ entry
	$\Psi_{ij}$ given by $\Psi_{ij} = \langle d\psi(e_i), d\psi(e_j)\rangle_N\, .$ Then the $\Phi_{(3)}$-energy density $e_{\Phi_{(3)}}(u)$ of $u$ satisfies
	\begin{equation}
	\begin{aligned}
	e_{\Phi_{(3)}}(u) & =\frac 16 \sum_{i,j,k=1}^n\langle du(e_i),du(e_j)\rangle \langle du(e_j),du(e_k)\rangle \langle du(e_k),du(e_i)\rangle\, \\
	& = \frac 16 \text {trace} (U \cdot U \cdot U)\, .
	\end{aligned}
	\end{equation}
	Similarly, $e_{\Phi_{(3)}}(\psi) =\frac 16 \text {trace} (\Psi^3)\, ,$ and
	\begin{equation}
	\begin{aligned}e_{\Phi_{(3)}}(u \circ \psi ) &=\frac 16 \text {trace}  (U^3 \cdot \Psi^3)\\
	&\le \max_ {x\in N} |U_{i,j}(x)|^3 \cdot \frac 16 \text {trace}  \Psi ^3\\
	&\le \max_ {x\in N} |U_{i,j}(x)|^3 \cdot e_{\Phi_{(3)}} (\Psi)\, .
	\end{aligned}
	\end{equation}
	
	Thus, if $u: N \to M$ is any smooth map, its composition with $\psi ^{\ell} : N \to N$ homotopic to the identity map on $N$ and $E_{\Phi_{(3)}} (\psi^{\ell}) \to 0\, ,$ as $\ell \to \infty$, then
	\begin{equation}
	\begin{aligned}E_{\Phi_{(3)}}(u \circ \psi ^{\ell}) &  \le \max_ {x\in N} |U_{i,j}(x)|^3 \cdot E_{\Phi_{(3)}} (\Psi^\ell) \\
	& \to 0\quad \text{as}\quad  \ell \to \infty\, .
	\end{aligned}
	\end{equation}
\end{proof}

\noindent
\begin{thm} \label{T:10.2}
	If $N$ is compact $\Phi_{(3)}$-$\operatorname{SSU}$ manifold, then for every compact Riemannian manifold $M$, the homotopic class of any map from $N$ into $M$ contains elements of arbitrarily small $\Phi_{(3)}$-energy $E_{\Phi_{(3)}}$.
\end{thm}

\begin{proof} This follows at once from Corollary \ref{C:4.1} and Lemma \ref{L: 4.4}.
\end{proof}

By virtue of Theorems \ref{thm:fss1}, \ref{thm:iss1}, \ref{T: 9.2}, \ref{T:10.2}, Definition \ref{D: 1.7},
\cite [p.131-132, Definition and Theorem A]{ww}, and \cite[p.5, Theorem 1.2]{hw}, we have the following result.

\begin{thm}\label{T: 10.3}
For $i=1,2,3$, every compact $\Phi_{(i)}$-$\operatorname{SSU}$ manifold is $\Phi_{(i)}$-$\operatorname{SU}$, and hence is $\Phi_{(i)}$-$\operatorname{U}.$
\end{thm}

This generalizes Proposition \ref{P:1.8}.\bigskip

{\bf{Acknowledgements}}: The authors wish to thank Professors Sun-Yung  Alice Chang and Paul C. Yang for references and communications, the editor for the editorship, and the referee for helpful comments, suggestions and remarks which make the present form of the paper possible.
Work supported in part by the National Natural Science Foundation of China (Grant No.11971415, 11771456), the NSF (DMS-1447008), and Nanhu Scholars Program for Young
Scholars of Xinyang Normal University.

\hspace*{1mm}{Shuxiang Feng and Yingbo Han \\}

\hspace*{1mm}{School of Mathematics and Statistics\\}
\hspace*{7mm}{Xinyang Normal University\\}
\hspace*{7mm}{Xinyang, 464000, Henan, P. R. China\\}
\hspace*{7mm}{shuxiangfeng78@163.com (Shuxiang Feng)}\\
\hspace*{7mm}{yingbohan@163.com (Yingbo Han)}\\

 \hspace*{1mm}{Kaige Jiang \\}

\hspace*{1mm}{School of Mathematics and Statistics\\}
\hspace*{7mm}{Shangqiu Normal University\\}
\hspace*{7mm}{Shangqiu, 476000, Henan, P. R. China\\}
\hspace*{7mm}{kaigejiang@163.com}\\

\hspace*{1mm}{Shihshu Walter Wei\\}

\hspace*{1mm}{Department of Mathematics\\}
\hspace*{7mm}{University of Oklahoma\\}
\hspace*{7mm}{Norman, Oklahoma 73019-0315, U.S.A.\\}
\hspace*{7mm}{wwei@ou.edu}\\

\end{document}